\documentclass{lmcs}


\usepackage{latexsym}
\usepackage{amssymb}
\usepackage{amsmath}
\usepackage{amsthm}
\usepackage{mathrsfs} 
\usepackage{accents} 
\usepackage{graphicx}
\usepackage{mathtools}
\usepackage{hyperref}

\newcommand{\adjoins}[1]{\mathrel{{{|}\mspace{-3.5mu}{\rightarrow}}_{\!#1}}}

\newcommand{\Aka}{A.k.a.}
\newcommand{\allabove}[2]{#1\mathbin{\uparrow}#2}
\newcommand{\allbelow}[2]{#1\mathbin{\downarrow}#2}
\newcommand{\anon}{\mathord{\tinysp\rule[0.5ex]{0.45em}{0.5pt}\tinysp}}
\newcommand{\apr}{a'}
\newcommand{\arr}[1]{\left\larr\arrmarg #1 \arrmarg\right\rarr}
\newcommand{\Asc}{\mathtextop{Asc}\nolimits}
\newcommand{\astbimpl}
  {\mathbin{\accentset{\scriptstyle*\mspace{7mu}}{\smash{\Rightarrow}\rule{0pt}{0.8ex}}}}

\newcommand{\bigset}[1]{\bigl\{ #1 \bigr\}}
\newcommand{\bigsuchthat}{\bigm|}
\newcommand{\bimpl}{\mathbin{\Rightarrow}}
\newcommand{\cl}{\mathtextop{cl$\tinysp$}\nolimits}
\newcommand{\Cl}{\mathtextop{Cl$\tinysp$}\nolimits}
\newcommand{\clop}{\mathtextop{cl$\dtinysp$op$\tinysp$}\nolimits}
\newcommand{\ClRul}{\mathtextop{ClRul}\nolimits}
\newcommand{\clsys}{\mathtextop{cl$\dtinysp$sys$\tinysp$}\nolimits}
\newcommand{\ClSys}{\mathtextop{ClSys}\nolimits}
\newcommand{\ClTheor}{\mathtextop{ClTheor}\nolimits}
\newcommand{\coll}{\mathcal}
\newcommand{\compose}{\circ}
\newcommand{\Dc}{\mathtextop{Dc}\nolimits}
\newcommand{\dcclsys}{\mathtextop{dc$\dtinysp$cl$\dtinysp$sys$\tinysp$}\nolimits}
\newcommand{\DcClSys}{\mathtextop{DcClSys}\nolimits}
\newcommand{\DcNucSys}{\mathtextop{DcNucSys}\nolimits}
\newcommand{\DcSpos}{\mathtextop{DcSpo}\nolimits}
\newcommand{\ddown}{\mathop{\rotatebox[origin=c]{-90}{\ensuremath{\twoheadrightarrow}}}\nolimits\negtinysp}
\newcommand{\defeq}{:=}
\newcommand{\Defeq}{\:\defeq\:}
\newcommand{\deltapr}{\delta'}
\newcommand{\dflt}{\mathtextopfont{df}}
\newcommand{\dircl}{directed-closed}
\newcommand{\dirinacc}{directed-inaccessible}
\newcommand{\D}{\mathscr{D}}
\newcommand{\DirSubsets}{\D}
\renewcommand{\dj}{\mathtextop{dj}\nolimits}
\newcommand{\dotbimpl}
  {\mathbin{\accentset{\textstyle\text{.}\mspace{5mu}}{\smash{\Rightarrow}\rule{0pt}{0.95ex}}}}
\newcommand{\dottJoin}{\mathop{\setbox0=\hbox{${\Join}{\cdot}$}\setbox1=\hbox{${\Join}$}%
	\makebox[\wd1][l]{${\Join}\hspace{-.5\wd0}\raisebox{.8ex}{${\cdot}$}$}}\nolimits}
\newcommand{\dtinysp}{\mspace{2mu}}

\renewcommand{\emph}[1]{\textit{#1\/}}
\renewcommand{\emptyset}{\mathord{\varnothing}}
\newcommand{\Eq}{\:=\:}

\newcommand{\fit}{{\rightt{\mathtextopfont{fit}}}}
\newcommand{\fitit}{\phi}
\newcommand{\fitnuc}{\Delta\tinysp}
\newcommand{\FitNuc}{\Nuc_\fit}
\newcommand{\FitNucSys}{\NucSys_\fit}
\newcommand{\fix}{\mathtextop{fix}\nolimits}
\newcommand{\fpr}{f'}
\newcommand{\Fun}{\mathtextop{Fun}\nolimits}
\newcommand{\galoisconn}{\rightleftharpoons}
\newcommand{\Galoisconn}{\xrightleftharpoons{\quad}}
\newcommand{\gammapr}{\gamma^{\tinysp\prime}}
\newcommand{\gencl}[1]{\overlinesymb{#1}}
\newcommand{\Heytalg}{Heyting algebra}
\newcommand{\HMJ}{Hofmann--Mislove--Johnstone}
\newcommand{\hpr}{h'}
\newcommand{\id}{\mathtextopfont{id}}
\newcommand{\ie}{i.e.}
\renewcommand{\iff}{if~and only~if}
\newcommand{\impl}{\mathrel{\Rightarrow}}
\renewcommand{\implies}{\mathrel{\Longrightarrow}}
\newcommand{\Implies}{\:\implies\:}
\newcommand{\In}{\:\in\:}
\newcommand{\Inc}{\mathtextop{Inc}\nolimits}

\newcommand{\inter}{{}^\circ}
\newcommand{\inters}{\cap}

\newcommand{\interskip}{}
\newcommand{\interthmskip}{\vspace{-.5ex}}
\newcommand{\isequiv}{\mathrel{\Longleftrightarrow}}
\newcommand{\Isequiv}{\:\isequiv\:}
\newcommand{\isimplied}{\mathrel{\Longleftarrow}}
\newcommand{\J}{\mathscr{J}}
\newcommand{\join}{\vee}
\renewcommand{\Join}{\bigvee}
\newcommand{\Joinable}{\!\J\negdtinysp}
\newcommand{\knownas}[1]{\textit{#1\/}}
\newcommand{\ldown}{\mathop{\smash{\downarrow}}\nolimits}
\newcommand{\leftdt}[1]{\negdtinysp#1\dtinysp}
\newcommand{\leftt}[1]{\negtinysp#1\tinysp}
\newcommand{\Leq}{\:\leq\:}
\newcommand{\longto}{\longrightarrow}
\newcommand{\lup}{\mathop{\smash{\uparrow}}\nolimits}
\newcommand{\mathtextopfont}[1]{\text{\rm#1}}
\newcommand{\mathtextop}[1]{\mathop{\smash{\mathtextopfont{#1}}}}
\newcommand{\meet}{\wedge}
\newcommand{\Meet}{\bigwedge}
\newcommand{\narr}[1]{\!#1\!}
\newcommand{\narrdt}[1]{\negdtinysp#1\negdtinysp}
\newcommand{\narrt}[1]{\negtinysp#1\negtinysp}
\newcommand{\negdisplayhalfskip}{}
\newcommand{\negdisplayshortskip}{}
\newcommand{\negdisplayskip}{}
\newcommand{\negdtinysp}{\mspace{-2mu}}

\newcommand{\negtinysp}{\mspace{-1mu}}
\newcommand{\NN}{\mathbb{N}}    
\newcommand{\notion}[1]{{\bfseries #1}}
\newcommand{\nuc}{\mathtextop{nuc}\nolimits}
\newcommand{\Nuc}{\mathtextop{Nuc}\nolimits}
\newcommand{\nuccore}{\mathtextop{nuc}\nolimits}
\newcommand{\nucfilt}{\mathtextop{nuc$\dtinysp$filt}\nolimits}
\newcommand{\NucFilt}{\mathtextop{NucFilt}\nolimits}
\newcommand{\nucsys}{\mathtextop{nuc$\dtinysp$sys}\nolimits}
\newcommand{\NucSys}{\mathtextop{NucSys}\nolimits}
\newcommand{\oneker}{\mathord{\nabla}\tinysp}
\newcommand{\op}{\mathtextopfont{op}}
\newcommand{\open}{{}^\circ}
\newcommand{\overlinepart}[2][1]%
  {{\setbox0=\hbox{$#2$}\accentset{\raisebox{0.0618ex}{\rule{#1\wd0}{0.09ex}}}{#2}}}
\newcommand{\overlinesymb}[1]{\overlinepart[0.8]{#1}}
\newcommand{\overbar}{\overlinesymb}
\newcommand{\pair}{\arr}
\renewcommand{\P}{\mathscr{P}}
\newcommand{\Pow}{\P}
\newcommand{\PowA}{\Pow\!A}
\newcommand{\PowE}{\Pow\negtinysp E}

\newcommand{\PowL}{\Pow\negtinysp L}
\newcommand{\PowP}{\Pow\negtinysp P}
\newcommand{\PowPowE}{\Pow\negdtinysp\Pow\negtinysp E}
\newcommand{\PowX}{\Pow\negtinysp X}
\newcommand{\Precl}{\mathtextop{Precl}\nolimits}
\newcommand{\Prenuc}{\mathtextop{Prenuc}\nolimits}
\newcommand{\regnuc}[1]{#1^\natural}
\newcommand{\relnot}{\lnot\tinysp}
\newcommand{\rightdt}[1]{\dtinysp#1\negdtinysp}
\newcommand{\rightt}[1]{\tinysp#1\negtinysp}
\newcommand{\Rulesdflt}{\mathrm{R}_{\dtinysp\dflt}}
\newcommand{\Rulesnuc}{\mathtextopfont{R}_{\nuc}}
\newcommand{\Sc}{\mathtextop{Sc}\nolimits}
\newcommand{\ScCl}{\mathtextop{ScCl}\nolimits}
\newcommand{\ScClSys}{\mathtextop{ScClSys}\nolimits}
\newcommand{\sccore}{\mathtextop{sc}\nolimits}
\newcommand{\ScNuc}{\mathtextop{ScNuc}\nolimits}
\newcommand{\ScNucSys}{\mathtextop{ScNucSys}\nolimits}
\newcommand{\Scottcl}{Scott-closed}
\newcommand{\Scottcont}{Scott-continuous}
\newcommand{\Scottop}{Scott-open}
\newcommand{\ScPrecl}{\mathtextop{ScPrecl}\nolimits}
\newcommand{\set}[1]{\{ #1 \}}
\newcommand{\setdiff}{\narrdt\setminus\negtinysp}
\newcommand{\stJoin}{\mathop{\smash{\textstyle\bigvee}}\nolimits}
\newcommand{\stMeet}{\mathop{\smash{\textstyle\bigwedge}}\nolimits}
\newcommand{\Subseteq}{\:\subseteq\:}
\newcommand{\suchthat}{\mid}
\newcommand{\theor}{\mathtextop{theor$\tinysp$}\nolimits}
\newcommand{\thmskip}{\medskip}
\newcommand{\tInters}{\mathop{\textstyle\bigcap}\nolimits}
\newcommand{\tinysp}{\mspace{1mu}}
\newcommand{\tJoin}{\mathop{\textstyle\bigvee}\nolimits}
\newcommand{\tMeet}{\mathop{\textstyle\bigwedge}\nolimits}
\newcommand{\tUnion}{\mathop{\textstyle\bigcup}\nolimits}
\newcommand{\txtskip}{\bigskip}
\newcommand{\union}{\cup}
\newcommand{\waybelow}{\ll}
\newcommand{\wide}[1]{\,#1\,}
\newcommand{\widedt}[1]{\dtinysp#1\dtinysp}
\newcommand{\wider}[1]{\:#1\:}
\newcommand{\widet}[1]{\tinysp#1\tinysp}

\newcommand{\nothing}{}
\newcommand{\setparenarr}{\let\larr(\let\rarr)\let\arrmarg\nothing}
\newcommand{\setanglearr}{\let\larr\langle\let\rarr\rangle\let\arrmarg\tinysp}
\setanglearr

\renewcommand{\leq}{\leqslant}

\renewcommand{\geq}{\geqslant}

\newcommand{\tJoinfof}{\tJoin\negdtinysp\negtinysp f}
\newcommand{\dottJoinF}{\dottJoin\negtinysp\!F}
\newcommand{\tJoinF}{\tJoin\negtinysp\!F}
\newcommand{\tJoinFof}{\tJoin\negdtinysp F}
\newcommand{\dottJoinG}{\dottJoin\!G}
\newcommand{\tJoinG}{\tJoin\!G}
\newcommand{\tJoinGof}{\tJoin\negtinysp G}
\newcommand{\tJoinA}{\tJoin\negdtinysp\!A}
\newcommand{\tJoinD}{\tJoin\negdtinysp\!D}
\newcommand{\tJoinY}{\tJoin\negdtinysp Y}

\newcommand{\tJoinS}{\tJoin\negtinysp\!S}

\newenvironment{myproof}[1][Proof.]{\begin{proof}[\normalfont\bfseries#1]}{\end{proof}}
\newcommand{\mycaption}[1]{{\scshape\caption{\normalfont\small#1}}}

\begin{document}

\title{Closure operators on dcpos}

\author[F.~Dacar]{France Dacar}
\address{Department of Knowledge Technologies,
	Jo\v{z}ef~Stefan Institute,
	Jamova 39,
	1000 Ljubljana,
	Slovenia}
\email{france.dacar@ijs.si}
\urladdr{\url{https://kt.ijs.si/france_dacar/}}

\begin{abstract}
We examine collective properties of closure operators on posets that are at least dcpos.
The first theorem sets the tone of the paper:
it tells how a~set of preclosure maps on a~dcpo determines the least closure operator above~it,
and pronounces the related induction principle and its sibling, the obverse induction principle.
Using this theorem we prove that the poset of closure operators on a~dcpo is a~complete lattice,
and then provide a~constructive proof of the Tarski's theorem for dcpos.
We go on to construct the joins
in the complete lattice of \Scottcont\ closure operators on a~dcpo,
and to prove that the complete lattice of nuclei on a~preframe is a~frame,
giving some constructions in the special case of the frame of all nuclei on a~frame.
In~the rather drawn-out proof of the \HMJ\ theorem
we show off the utility of the obverse induction, applying it in the proof of the clinching lemma.
After that we shift a~viewpoint and prove some results, analogous to the results about dcpos,
for posets in which certain special subposets have enough maximal elements;
these results specialize to dcpos, but at the price of using the axiom of choice.
We conclude by pointing out two convex geometries associated with closure operators on a~dcpo.
\end{abstract}

\maketitle

\section{Overview of the paper}
\label{sec:overview}

The central theme of the paper are the properties
of the poset $\Cl(\negtinysp P)$ of all closure operators on a~poset $P$ that is at least a~dcpo,
and of subposets of $\Cl(\negtinysp P)$ consisting of some special kind of closure operators,
for instance of the subposet $\ScCl(\negtinysp P)$ of all \Scottcont\ closure operators.
The~results of the paper are, among other things,
subsuming and/or extending several known fixed point theorems for~dcpos.

\txtskip

In section~\ref{sec:prelims} we present notation and terminology, and state some basic facts.

\txtskip

In section~\ref{sec:completlatt-of-clopers-on-dcpo}
	we establish that, for a~dcpo $P$, the poset $\Cl(\negtinysp P)$ is a~complete lattice.
More is true, actually: the complete lattice $\ClSys(\negtinysp P)$ of all closure systems in~$P$,
	which is antiisomorphic to the complete lattice $\Cl(\negtinysp P)$,
is a~closure system in the powerset lattice~$\PowP\tinysp$
(which means that the intersection of any~set of closure systems in~$P$
	is a~closure system~in~$P$);
also, the fixpoint set of any preclosure map on~$P$ is a~closure system.
These results are consequences of Theorem~\ref{thm:in-dcpo-cl-generd-by-precls-&-induct}
which describes how a~set~$G$ of preclosure maps on a~dcpo~$P$
	determines the least closure operator~$\gencl{G}$ above~it.

Theorem~\ref{thm:in-dcpo-cl-generd-by-precls-&-induct}
	also pronounces the \knownas{induction principle}:
if a~subset of~$P$ is closed under directed joins in~$P$ and is closed under~$G$,
then it is closed under~$\gencl{G}$.
The induction principle has a~sort of dual, the \knownas{obverse induction principle},
but the passage from the former to the latter involves the~law of~excluded middle;
since the one application of the obverse induction principle in the paper
is in a~proof which intentionally avoids using the~law of~excluded middle,
the obverse induction principle is stated and proved~on~its~own.

Using Theorem~\ref{thm:in-dcpo-cl-generd-by-precls-&-induct} and its corollaries
we provide a~constructive proof of the version of Tarski's fixed point theorem for dcpos.

\txtskip

In Section~\ref{sec:Scottcont-clopers-on-dcpo}
we consider the poset $\ScCl(\negtinysp P)$ of all \Scottcont\ closure operators on a~dcpo~$P$.
The main result of the section
	is~Theorem~\ref{thm:in-dcpo-Scottcont-preclmaps-gener-Scottcont-clop}:
if $G$ is a~set of \Scottcont\ preclosure maps on $P$,
then the least closure operator $\gencl{G}$ above~$G$
is the directed pointwise join of the composition monoid~$G^*$ generated~by~$G$,
and $\gencl{G}$ is~\Scottcont.
Consequently, the poset $\ScCl(\negtinysp P)$ is an interior system in the complete lattice~$\Cl(\negtinysp P)$
	and so is itself a~complete lattice with the joins inherited from~$\Cl(\negtinysp P)$,
while the poset $\DcClSys(\negtinysp P)$ of all \dircl\ closure systems in~$P$,
	which is in fact the complete lattice
		of the closure systems associated with the \Scottcont\ closure operators,
is~a~closure system in~$\ClSys(\negtinysp P)$ as~well~as~in~$\PowP$.

For~every closure operator~$\gamma$ on~$P$
there exists the greatest \Scottcont\ closure operator $\widedt{\sccore(\gamma)}$ below~$\gamma\tinysp$,
called the \knownas{\Scottcont\ core} of the closure operator~$\gamma\tinysp$.
Correspondingly,~for every closure system~$C$ in~$P$
there exists the least \dircl\ closure system~$\dtinysp\dcclsys(C)$ that includes~$C$.
Not much can be said about $\widedt{\sccore(\gamma)}$ and~$\dtinysp\dcclsys(C)$
	for a~general dcpo~$P$,
but if $P$ is a~domain (a~continuous dcpo), then both can~be constructed.
Proposition~\ref{prop:sc(gamma)-on-domain}
	gives the construction of the \Scottcont\ core of~a~closure operator on a~domain~$P$,
while Proposition~\ref{prop:dcclsys-in-domain}
	has the construction of~$\dtinysp\dcclsys(C)$ for~a~closure system~$C$ in a~domain~$P$.

\txtskip

In section~\ref{sec:frame-of-nuclei-on-preframe} we carry out the project
that is only sketched in~\cite{Escardo-JFN-2003}:
we prove, in the fullness of time, that the poset of all nuclei on a~preframe is a~frame.

A~\knownas{frame} is a~complete lattice in which binary meets distribute over arbitrary joins.
A~\knownas{preframe} is a~meet-semilattice that is also a~dcpo
	and in which binary meets distribute over directed joins.
Given a~meet-semilattice~$P$,
a~\knownas{nucleus} (\knownas{prenucleus}) on~$P$
	is a~closure operator (a~preclosure map) on~$P$ that preserves binary meets;
the fixpoint set of a~nucleus on $P$ is called a~\knownas{nuclear system} in~$P$.
We denote by $\Nuc(\negtinysp P)$ the poset of all nuclei on~$P$
and by $\NucSys(\negtinysp P)$ the poset of all nuclear systems in~$P$.

The~main result of the section is Theorem~\ref{thm:in-preframe-prenuclei-generate-nucleus}
which states that for any set $\Gamma$ of~prenuclei in a~preframe~$P$
the least closure operator $\gencl{\Gamma}$ above the set of preclosure maps~$\Gamma$ is~a~nucleus.
It~follows that for a~preframe~$P$
the poset $\Nuc(\negtinysp P)$ is closed under all joins in the complete lattice $\Cl(\negtinysp P)$
	and is thus a~complete lattice;
besides this, $\Nuc(\negtinysp P)$ is closed under the (pointwise calculated) binary meets in~$\Cl(\negtinysp P)$.
For every closure operator~$\gamma$ on a~preframe~$P$
there exists the largest nucleus $\widedt{\nuccore(\gamma)}$ below~$\gamma\tinysp$,
the \knownas{nuclear core} of~$\gamma\tinysp$.
Since $\Nuc(\negtinysp P)$ is an interior system in $\Cl(\negtinysp P)$,
$\NucSys(\negtinysp P)$ is a~closure system in $\ClSys(\negtinysp P)$ and~hence~in~$\PowP$.

Theorem~\ref{thm:frame-of-nuclei-on-preframe} tells~us
that on a~preframe $P$ the complete lattice $\Nuc(\negtinysp P)$ is in~fact~a~frame.
Both Theorem~\ref{thm:in-preframe-prenuclei-generate-nucleus}
	and Theorem~\ref{thm:frame-of-nuclei-on-preframe}
are proved using the induction principle.

We conclude the section	by taking a~look at \Scottcont\ nuclei on~$P$.

\txtskip

In section~\ref{sec:nucs-on-frames} we consider the nuclei on a~frame~$L\tinysp$.
Since a~frame is a~special preframe,
all results for preframes specialize to the frame~$L\tinysp$.
But, since a~frame is a~\emph{very} special preframe,
we can say much more about the frame of nuclei $\Nuc(L)$ on the frame $L$
than about the frame of nuclei on a~mere preframe.
For~instance, by~Proposition~\ref{prop:frame-L-Nuc(L)-sub-(completlatt)-of-Cl(L)},
the subset~$\Nuc(L)$ of~$\Cl(L)$ is closed not only under arbitrary joins in~$\Cl(L)$
but also under arbitrary (not~just~binary) meets in~$\Cl(L)\dtinysp$;
that is, $\Nuc(L)$~is a~sub\nobreakdash-(complete lattice) of the complete lattice~$\Cl(L)$.

A~frame is relatively pseudocomplemented, that is, it is a~complete \Heytalg.
Corollary~\ref{cor:in-frame-nucC(y)=Meet(x-in-X)((y=>gamma(x))=>gamma(x))}
gives a~formula, which uses the operation $\impl$ of relative pseudocomplementation,
for the nuclear core $\widedt{\nuccore(\gamma)}$ of a~closure operator $\gamma$ on~$L\tinysp$.

\txtskip

In~\cite{Escardo-JFN-2003} the author demonstrates the utility of join induction
	in a~proof of the \HMJ\ theorem.
In~section~\ref{sec:HMJ-theorem} we prove this theorem in a~way
that let us observe the obverse induction principle in action.
Our~proof of the \HMJ\ theorem is spread through proofs of three lemmas,
	with parts of it reasoned out in the connecting text;
the short concluding reasoning then ties everything together.
The~obverse induction principle is used in the proof of Lemma~\ref{lem:scott-open-filter-is-nuclear}.

\txtskip

For a~dcpo~$P$, the complete lattice $\ClSys(\negtinysp P)$ of all closure systems in~$P$
	is a~closure system in the powerset lattice $\PowP$
and so it is determined by a~set of closure rules on~$P$,
say by the set of all closure rules obeyed by~$\ClSys(\negtinysp P)$.
In section~\ref{sec:do-it-with-max-elems}
we prove that $\ClSys(\negtinysp P)$ is already determined
	by the set of all~\knownas{default closure rules associated with~$P$},
which are the closure rules on~$P$ of the form $B\negdtinysp\adjoins{}c\tinysp$,
where $c\in P$ is a~maximal lower bound of $B\subseteq P$.
However, the proof of this result requires the~axiom of~choice.

In section~\ref{sec:do-it-with-max-elems}
	we actually develop a~little theory which operates with maximal elements.
We prove several assertions of the following form:
if every subset of~$P$ of some special kind \knownas{has a~ceiling},
meaning that in the subset every element has a~maximal element above it,
then $P$~has a~certain property.
For example, if $P$ is \knownas{default-enabled},
	which means that every lower bound of any~subset of~$P$
		is below some maximal lower bound of the subset
	(that~is, the sets of all lower bounds of arbitrary subsets of~$P$ have ceilings),
then the closure systems in~$P$ are determined by the set of all default closure rules;
in other words, if $P$ has enough default closure rules associated with~it,
then the closure system $\ClSys(\negtinysp P)$ in~$\PowP$
	is determined by the default closure rules.
Similarly, if we require that in a~default-enabled meet-semilattice
	some special subsets have ceilings,
	then there are enough \knownas{nuclear} closure rules,
		besides the default closure rules,
so that $\Nuc(\negtinysp P)$ is an interior system in~$\Cl(\negtinysp P)\tinysp$;
if in addition to this we require existence of still more subsets that have ceilings,
	then we are able to prove that $\Nuc(\negtinysp P)$~is~a~frame.
Interestingly, we can prove all this without ever invoking the~axiom of~choice.
These results mimic the results
	in sections~\ref{sec:completlatt-of-clopers-on-dcpo}
		and~\ref{sec:frame-of-nuclei-on-preframe},
and they in fact imply them by specialization,
though at~the~price of being compelled to involve the~axiom of~choice.

The class of default-enabled posets is strictly larger than the class of~dcpos.

Every~default-enabled poset~$P$, and in particular every dcpo, has the following two properties:
the fixpoint set of every~preclosure map on~$P$ is a~closure system~in~$P\tinysp$;
the set of all closure systems in~$P$ is a~closure system in~$\PowP$.
Here is a~project that may turn out to be more of an adventure than it appears
	at first sight:
characterize, in structural terms,
	the posets that have the one, or the other, or both of these properties.

\txtskip

In section~\ref{sec:conv-geoms-assoc-with-dcpo} we prove that for every dcpo~$P$
the closure operator~$\widedt{\clsys_P}$ on~$\PowP$
    (which for each $X\subseteq P$
	yields the least closure system in~$P$ that includes~$X$)
and the closure operator $\widedt{\dcclsys_P}$ on~$\PowP$
    (which for each~$X\subseteq P$
	yields the least \dircl\ closure system in $P$ that includes~$X$)
are convex,
meaning that they satisfy the anti-exchange axiom.
Since all that we need to obtain these two results
is the property of a~dcpo~$P$ that $\ClSys(\negtinysp P)$ is a~closure system in~$\PowP$,
both~results are valid also for every default-enabled poset~$P$.

\section{Preliminaries}
\label{sec:prelims}

For any set $X$ we denote the set (and complete lattice) of all subsets of~$X$
	by~$\PowX$.

For a~subset~$A$ of a~monoid~$M$
	we denote by $A^*$ the submonoid of $M$ generated~by~$A\tinysp$.
An \notion{ordered monoid} is a~monoid $M$ which is partially ordered by $\leq$
	and whose operation is increasing in both operands,
meaning that $x_1\leq y_1$ and $x_2\leq y_2$ imply $x_1x_2\leq y_1y_2$
	for all $x_1,\tinysp x_2,\tinysp y_1,\tinysp y_2\in M$;
this condition is satisfied iff $x\leq y$ implies $zx\leq zy$ and $xz\leq yz$
	for all $x,\tinysp y,\tinysp z \in M$.

If $X$ and $Y$ are sets,
	then $\Fun(X,Y)$ denotes the set of all functions from~$X$ to~$Y$.
We write $\Fun(X,X)$ as $\Fun(X)$.
The set $\Fun(X)$ carries the structure of a~monoid,
with the monoid operation the composition of endomaps on~$X$
and with the neutral element the identity function~$\id_X$ on~$X$.
A~\notion{composition monoid} is a~submonoid of the monoid~$\Fun(X)$ for some set~$X$.

If~$X$ and~$Y$ are~sets, $A\subseteq X$, and $F\subseteq\Fun(X,Y)$,
then we denote by $F(A)$ the subset $\set{f(a)\suchthat \text{$f\in F$ and $a\in A$}}$ of~$B$.
We write $\set{f}(A)$ as $f(A)$ and $F\bigl(\set{a}\bigr)$ as~$F(a)$.

\txtskip

Let $X$ be a set and $F$ a set of endofunctions on~$X$.

A~\notion{fixed point}, or~\notion{fixpoint}, \notion{of~$F$}
	is an~$x\in X$ such that $f(x) = x$ for every~$f\in\nolinebreak F$.
We~denote the set of all fixed~points~of~$F$~by~$\dtinysp\fix(F)$,
	and call it the~\notion{fixpoint~set~of~$F$}.
We~write $\dtinysp\fix(\set{f})$ as~$\dtinysp\fix(f)\tinysp$;
always $\dtinysp\fix(f)\narrt\subseteq f(X)$,
where $\dtinysp\fix(f) = f(X)$ iff $f$~is an idempotent endofunction.
Note that $\dtinysp\fix(F) = \tInters_{f\in F}\fix(f)\tinysp$;
in particular $\dtinysp\fix(\emptyset) = X$.%
\footnote{\,The `intersection' of the empty set of subsets of~$X$
	is by convention its meet $X$ in~$\PowX$.}

Let $A\subseteq X$.
The~subset $A$ is said to be \notion{closed under~$F$} if $F(A)\subseteq A\dtinysp$;
as a~special case, $a\in X$ is a fixed point of $F$ iff $\set{a}$ is closed under~$F$.
The~subset~$A$~is said to be \notion{inversely closed under~$F$}
	if~$f^{-1}(A)\subseteq A$ for every $f\in F$.
In the presence of the~law of~excluded middle,
the subset $A$~is closed under~$F$
iff its complement $X\setdiff A$ is~inversely~closed under~$F\tinysp$:
the statement
    $(\forall\negtinysp f\narr\in F)(\forall x\narr\in X)
		(x\narrdt\in A\widet\impl f(x)\narrdt\in A)$,
which says that $A$~is~closed under~$F$,
is~equivalent
to the statement
    $(\forall\negtinysp f\narr\in F)(\forall x\narr\in X)
		(f(x)\narrdt\notin A \widet\impl x\narrdt\notin A)$,
which says that $X\setdiff A$ is~inversely closed under~$F$.

\txtskip

Let $P$ be a~poset (ordered by $\leq\dtinysp$).

Of $x,\tinysp y\in P$ such that $x\leq y$
we shall say that $x$~is \notion{below}~$y$ and that $y$~is \notion{above}~$x\tinysp$.

If the poset~$P$ has a~least (a~greatest) element, we shall write it
	as $\bot \narrdt= \bot_{\tinysp P}\tinysp$ ($\top \narrdt= \top_{\!P}$).

Let $x\in P$ and $A\subseteq P$.
We shall denote by $\allabove{A}{x}$ the set $\set{y\narrt\in A\negtinysp\suchthat y\geq x}$
	of all~ele\-ments of~$A$ above~$x\tinysp$,
and by $\allbelow{A}{x}$ the set $\set{y\narrt\in A\negtinysp\suchthat y\leq x}$
	of all elements of~$A$ below~$x\tinysp$.
In~particular, $\lup x \defeq \allabove{P}{x}$
	is the \notion{principal filter of $P$} generated by~$x\tinysp$,
and $\ldown x \defeq \allbelow{P\negtinysp}{x}$
	is the \notion{principal ideal of $P$} generated by~$x\tinysp$.
We~shall~write~$x\geq A$ ($x\leq A$),
	and say that $x$~is \notion{above~$A$} (\notion{below}~$A$),
to mean that $x\geq a$ ($x\leq a$) for every~$a\in A\tinysp$,
	in~other words, that~$x$~is~an upper (lower) bound~of~$A\tinysp$.
A~\notion{lower set of~$P$} is a~subset $A$ of~$P$
	such that $x\in A\tinysp$, $y\in P$, and $x\geq y$ imply $y\in A\tinysp$.
Dually, an~\notion{upper set of~$P$} is a~subset~$A$ of~$P$
	such that $x\in A\tinysp$, $y\in P$, and $x\leq y$ imply $y\in A\tinysp$.
In the presence of the~law of~excluded middle,
a~subset~$A$ of the~poset~$P$ is a~lower set of~$P$
	iff its complement $P\setdiff A$ is an upper set of~$P$.

We assume as known the notions, and the basic properties, of
	a~\knownas{join} (least upper bound),
	a~\knownas{meet} (greatest upper bound),
	a~\knownas{join/meet-semilattice},
	a~\knownas{lattice},
	a~\knownas{complete lattice},
	and of~a~\knownas{(complete) meet/join-irreducible/prime element}.
Sometimes we will observe meets and joins in several posets,
	possibly in several subposets of the same poset;
on such occasions we shall avoid confusion by writing the symbols for joins and meets
		in each poset~$P$ under observation
	as $\join^P$, $\meet^{\negtinysp P}$, $\stJoin^P$, and $\stMeet^{\!P}$.

We shall say that a~subset~$X$ of the poset~$P$
	is \notion{joinable in~$P$} if it has a~join in~$P$. 
The~set of all joinable subsets of the poset~$P$ shall be denoted by $\Joinable\negtinysp P$.
Let $A\subseteq B$ be subposets of a~poset~$P$, and let $X\subseteq A\tinysp$.
Mark that $X$~being joinable in~$A$ is in general unrelated to it being joinable in~$B\tinysp$:
$X$ may have joins in both subposets, have a~join in only one of them, or in none.
Therefore, neither of the two inclusions
	$\Joinable{A}\subseteq\Joinable{B}$, $\Joinable{A}\supseteq\Joinable{B}$
holds in general.
When~$X$ does~have both joins $\stJoin^A\!X$ and~$\stJoin^B\!X$,
	then $\stJoin^A\!X\geq\stJoin^B\!X\tinysp$.

A~subset $D$ of the poset~$P$ is said to be \notion{directed}
if every finite (possibly empty) subset of $D$ has an~upper bound in~$D$,
or equivalently,
if it is nonempty and any two elements of $D$ have an upper bound in~$D$.
Dually, a~subset $F$ of the poset $P$ is said to be \notion{filtered}
if every finite subset of $F$ has a~lower bound in~$F$.%
\footnote{\,In this we follow~\cite{Gierz-etal-CL&D-2003}.
Some authors use ``filtered'' as a~synonym for ``directed''.}
%
We denote by $\DirSubsets P$
	the set of all directed subsets of the poset~$P$.
If~$A\subseteq P$,
then a~subset~$D$ of~$A$ is directed in the poset~$P$ iff it is directed in the subposet~$A\dtinysp$:
$\DirSubsets\negtinysp{A} = \PowA\widet\inters\DirSubsets{P}\tinysp$.
We denote by $\Joinable\negtinysp\DirSubsets{P}$ 
	the set of all joinable directed subsets of the poset~$P\tinysp$:
$\Joinable\negtinysp\DirSubsets{P} = \Joinable{P} \inters \DirSubsets{P}$.

\pagebreak[3]

A~subset~$A$ of the poset~$P$ is said to be \notion{closed under (existing) directed~joins},
	or~shorter, \notion{\dircl},
if for every joinable directed subset~$D$ of~$P$
	its join in~$P$ belongs~to~$A\tinysp$.
A~subset~$A$ of~$P$ is said to be \notion{inaccessible by (existing) directed~joins},
	or~shorter, \notion{\dirinacc},
if every joinable directed subset~$D$ of~$P$ whose join in~$P$ lies in~$A$
	has some element in~$A\tinysp$. 
In~the presence of~the~law of~excluded middle,
a~subset~$A$ of~$P$ is \dircl~iff 
	its complement $P\setdiff A$ is \dirinacc;
indeed, the~statement
	$(\forall D\narr\in\Joinable\negtinysp\DirSubsets\negtinysp{P})
			(D\narrdt\subseteq A \widet\impl \tJoin\!D\narrdt\in A)$,
which says that $A$~is~\dircl,
is equivalent to the statement
	$(\forall D\narr\in\Joinable\negtinysp\DirSubsets\negtinysp{P})
			(\tJoin\!D\narrdt\notin A \widet\impl D\narrdt\nsubseteq A)$,
which says that $P\setdiff A$~is~\dirinacc.

A~subset of the poset~$P$
	is said to be~\notion{\Scottcl} if it is a~\dircl\ lower set
	and it is said to be~\notion{\Scottop} if it is a~\dirinacc\ upper set.
In the presence of the law of excluded middle,
	a~subset of $P$ is \Scottop\ iff its complement is~\Scottcl.


A~\notion{directed-complete poset}, or a~\notion{dcpo} for short,
is a~poset in which every directed subset has a~join;
mark that the empty poset is a~dcpo.
A~dcpo is said to be \notion{pointed} if it has a~least element.
A~\dircl\ subposet of a~dcpo~$P$ is a~dcpo;
it is a~\notion{sub-dcpo~of~$P$},
	meaning that its directed joins are inherited from the dcpo~$P$.

\txtskip

A~function $f$ from a~poset~$P$ to a~poset~$Q$ is \notion{increasing} (\notion{decreasing})
if~$x\leq y$ in~$P$ always implies $f(x)\leq f(y)$ ($f(x)\geq f(y)$)~in~$Q$.
If~$f\colon P\to Q$ and $g\colon Q\to R$ are increasing functions between posets,
then the composite $gf = g\narrt\compose f\colon P\to R$ is increasing.
The set of all increasing functions $P\to Q$ shall be denoted by $\Inc(P,Q)$,
and $\Inc(P,P)$ will be shortened to $\Inc(\negtinysp P)$.
If $P$ is a~poset ordered by $\leq\dtinysp$,
then the poset on the same set of elements, but ordered by~$\geq\dtinysp$,
is said to be the \notion{opposite}, or \notion{dual}, \notion{of~$P$}
	and is denoted by $P^{\tinysp\op}$.
If~$P$~and~$Q$ are posets,
then an increasing function $P\to Q$ is also an increasing function $P^{\tinysp\op}\to Q^\op$,
and a~decreasing function $P\to Q$
is~the~same~thing as an increasing function $P\to Q^\op$
	or an increasing function $P^{\tinysp\op}\to Q\tinysp$.

For any set $X$ and any poset $P$ the set $\Fun(X,P)$ of all functions $X\to P$
	is made into a~poset by endowing it with the pointwise partial ordering:
if $f,\tinysp g\colon X\to P$, then we let $f\leq g$ iff $f(x)\leq g(x)$ for every $x\in X$.
For any posets $P$ and~$Q$ the set $\Inc(P,Q)$ is~made~into a~poset as a~subposet of $\Fun(P,Q)$.
Let~$P$,~$Q$,~$R$ be~posets.
The composition
\begin{equation*}
\Inc(P,Q)\narrt\times\Inc(Q,R) \wide\longto \Inc(P,R) \wide: \pair{f,g} \wide\longmapsto gf
\end{equation*}
is~increasing in both operands:
if $f_1\leq f_2$ in $\Inc(P,Q)$ and $g\in\Inc(Q,R)\tinysp$,
	then $gf_1\leq gf_2$
	because $\Inc(P,Q)$ and $\Inc(P,R)$ are ordered pointwise and $g$ is increasing,
and if we have $f\in\Inc(P,Q)$ and $g_1\leq g_2$ in $\Inc(Q,R)$,
	then $g_1f\leq g_2f$
	because $\Inc(Q,R)$ and $\Inc(P,R)$ are ordered pointwise.
In particular, $\Inc(\negtinysp P)$ is an ordered composition monoid.

\txtskip

Let $P$ be a~poset.

Let $f$ be an endofunction on the poset~$P$.
We say that $f$ \notion{ascends} (\notion{descends}) on $x\in P$
	if $x\leq f(x)$ ($x\geq f(x)$).
The function $f$ is \notion{ascending} (\notion{descending})
	if $f$ ascends (descends) on every element of~$P$.
An~ascending and increasing function $f$ is called a~\notion{preclosure~map on~$P$}.%
\footnote{\Aka\ an~\knownas{inflationary map}, as in~\cite{Escardo-JFN-2003}.
	This kind of map is known under many names.}
The set $\Asc\tinysp(\negtinysp P)$ of all ascending maps on the~poset~$P$ is a~composition monoid;
it is also a~poset (ordered pointwise),
	but~it~is~not, in general, an ordered monoid.
However, the set $\Precl(\negtinysp P)$ of all preclosure maps on~$P$
is a~submonoid, and hence a~sub-(ordered monoid),%
\footnote{\,Every submonoid of an ordered monoid $M$ becomes a~sub-(ordered monoid)
when we equip it with the partial order induced from~$M$.}
of~the~ordered composition monoid~$\Inc\tinysp(\negtinysp P)$.


A~\notion{closure operator on~$P$} is an idempotent preclosure map on~$P$,
that is, it is an endofunction on $P$ that is ascending, increasing, and idempotent;
the poset of all closure operators on $P$, ordered pointwise, is denoted by~$\Cl(\negtinysp P)$.
Dually, an \notion{interior operator on~$P$} is a~closure operator on~$P^{\tinysp\op}$,
that is, it is an endofunction on~$P$ that is descending, increasing, and idempotent.

A~\notion{closure system in~$P$}
	is a~fixpoint set of some closure operator~$\gamma$ on~$P\tinysp$;
note that $\dtinysp\fix(\gamma)=\gamma(\negtinysp P)$ since $\gamma$ is idempotent.
A~subset~$C$ of~$P$ is a~closure system in~$P$
	iff for every $x\in P$ the set $\allabove{C}{x}$ has a~least element.
If $C$ is a~closure system in~$P$,
then there is a~unique closure operator on~$P$, denoted~by~$\cl_C\tinysp$,
	whose fixpoint set is~$C\tinysp$:
it sends each $x\narrt\in P$ to the least element of $\allabove{C}{x}$.
We~de\-note~by~$\ClSys(\negtinysp P)$ the poset of all closure systems in~$P$ ordered by inclusion.
The~map $\Cl(\negtinysp P)\to\ClSys(\negtinysp P) : \gamma\mapsto\fix(\gamma)$
	is an antiisomorphism of posets,
with the inverse $\ClSys(\negtinysp P)\to\Cl(\negtinysp P) : C\mapsto\cl_C\tinysp$.
The~notion of an~\notion{interior system} is dual to that of a~closure system;
that is, an interior system in~$P$ is a~fixpoint set of some interior operator on~$P$.


Let $\gamma$ be a~closure operator on~$P$.
If~a~subset $X$ of the closure system $\gamma(\negtinysp P)$
	has a~join in the subposet~$\gamma(\negtinysp P)$ of~$P$,
we write it $\tJoin^{\tinysp\gamma}\!X$
	{\large(}instead of $\stJoin^{\tinysp\gamma(\negtinysp P)}\!X${\large)}.
If~a~subset~$S$ of~$P$ has a~join $\tJoinS$ in the poset~$P$,
then the~subset~$\gamma(S)$ of~$\gamma(P)$ has in the subposet~$\gamma(P)$ the~join
	$\tJoin^{\tinysp\gamma}\!\gamma(S) = \gamma(\tJoinS)\dtinysp$;
in~particular, if~$S\subseteq\gamma(P)$ has a~join $\tJoin\negdtinysp\!S$ in~$P$,
then~$S = \gamma(S)$ has the~join $\tJoin^{\tinysp\gamma}\negtinysp\!S = \gamma(\tJoinS)$ in~$\gamma(P)$.%
\footnote{\,However, a~subset of $\gamma(P)$ may have a~join in~$\gamma(P)$
	without having a~join in~$P\tinysp$.}
Restricting the codomain of the closure operator~$\gamma$ to its image~$\gamma(P)$
we obtain the surjective map $\gammapr\colon P\to \gamma(P) : x\mapsto\gamma(x)$
	which preserves all existing~joins.
How~about the meets in the subposet~$\gamma(P)\dtinysp$?
Let $X\subseteq \gamma(P)$.
If~$X$~has a~meet in~$P$ then this is also its meet in~$\gamma(P)$,
and if~$X$~has a~meet in~$\gamma(P)$ then this is also its~meet~in~$P$.


If~$g$ is a~preclosure map on $P$ and $h$ is a~closure operator on $P$,
then $g\leq h$ iff~$gh=h$, iff~$hg=h$, iff~$\dtinysp\fix(g)\supseteq\fix(h)\tinysp$;
this holds, in particular, if both~$g$ and~$h$ are closure operators.


Let $\beta$ and $\gamma$ be closure operators on a~poset~$P$
	such that $\beta\gamma\leq\gamma\beta\dtinysp$;
then $\gamma\beta$ is a~closure operator on~$P$,
and it is the join of the closure operators $\gamma$ and $\beta$ in the poset~$\Cl(\negtinysp P)$.


\txtskip

At this point we could go on, for quite a while, about Galois connections and their basic properties,
	but instead we shall assume that the reader knows about all this.

\txtskip

Let $L$ be a~complete lattice.

A~subset of $L$ is a~closure system \iff\ it is closed under arbitrary meets.
The~set~$\ClSys(L)$ of all closure systems in~$L$ is closed under arbitrary intersections,
and is therefore a~closure system in the complete lattice~$\PowL\dtinysp$;
consequently, the poset $\ClSys(L)$ is a~complete lattice.
The~poset~$\Cl(L)$, which is antiisomorphic to the poset~$\ClSys(L)$,
	is~likewisee a~complete lattice.
Let $\Gamma\subseteq\Cl(L)$.
The~meet of~$\Gamma$ in~$\Cl(L)$ is calculated pointwise,
    $\bigl(\tMeet\negdtinysp\Gamma\bigr)(x)
	= \tMeet\negtinysp\Gamma(x)$
    for $x\in L\tinysp$,
and the fixpoint set of the meet is obtained~as
    $\dtinysp\fix(\tMeet\negdtinysp\Gamma)
	= \set{\tinysp\tMeet_{\gamma\in\Gamma}z_\gamma
		\narrdt\suchthat\text{$z_\gamma\in\fix(\gamma)$ for $\gamma\in\Gamma$}}$.
The~join~$\tJoin\negtinysp\Gamma$ in~$\Cl(L)$ is the closure operator
whose fixpoint set is~$\dtinysp\fix(\Gamma)$.

Let $X$ be a~subset of $L\tinysp$.
The subset $X$ generates the closure system $\widedt{\clsys(X)}$,
	the~least closure system in $L$ that includes~$X$.
The closure system $\widedt{\clsys(X)}$ is the set of the meets of all subsets of~$X$.
If $\Gamma\subseteq\Cl(L)$,
then the fixpoint set of $\tMeet\negdtinysp\Gamma$
	is $\clsys(\tUnion_{\gamma\in\Gamma}\fix(\gamma))$,
which gives us the expression for $\dtinysp\fix(\tMeet\negdtinysp\Gamma)$ in the preceding paragraph.
The closure operator~$\cl_X$ on~$L\tinysp$, whose fixpoint set is $\widedt{\clsys(X)}$, is~given,
for $y\in L$,
by~$\dtinysp\cl_X(y)
	= \tMeet(\allabove{\clsys(X)}{y})
	= \tMeet(\allabove{X\negdtinysp}{y})$.

\txtskip

We now turn to closure operators on, and closure systems in,
	the complete lattice~$\PowE$ of all subsets of some set~$E\tinysp$.


A~\notion{closure rule on $E$}
is a~pair $\pair{B,c}\in\ClRul(E)\defeq(\PowE)\narrdt\times E\tinysp$,
which will be~written as~$B\negdtinysp\adjoins{}c\tinysp$,
with the subset~$B$ of~$E$ called the~\notion{body}
	and the element~$c$ of~$E$ called the~\notion{head}
of the closure rule.
For $R$ a~set of closure rules on~$E$
we~write~$\pair{B,c}\in R$~as~$R\colon B\negdtinysp\adjoins{}c\dtinysp$.
When the `ambient set' $E$ is known and we write $B\negdtinysp\adjoins{}c\tinysp$,
it is understood that $B$~is a~subset of~$E$ and $c$~is an element of~$E$.

A~\notion{closure theory on $E$} is a~set $T$ of closure rules on~$E$ which is
\begin{itemize}[topsep=1ex,itemsep=1ex,leftmargin=3cm]
\item[\notion{reflexive}\,:\,]
	for all $B\subseteq E$ and all $b\in E\tinysp$,\\
	if $b\in B$ then $T\colon B\negdtinysp\adjoins{}b\dtinysp$;
\item[\notion{transitive}\,:\,]
	for all $B,\tinysp C\subseteq E$ and all $d\in E\tinysp$,\\
	if $T\colon B\negdtinysp\adjoins{}c$ for every $c\in C$
		and $T\colon C\negtinysp\adjoins{}d\tinysp$,
	then $T\colon B\negdtinysp\adjoins{}d\tinysp$.
\end{itemize}
We denote by $\ClTheor(E)$ the poset of all closure theories on $E$
	ordered by inclusion.%
\footnote{A set of closure rules is also known as~\knownas{implicational system},
while a~closure theory is referred to as a~\knownas{complete implicational system}.}

We say that a subset $X$ of $E$ \notion{obeys} a~closure rule $B\negdtinysp\adjoins{}c$ on $E$,
or that the closure rule $B\negdtinysp\adjoins{}c$ \notion{is obeyed by} the subset~$X$,
if $B\subseteq X$ implies $c\in X$.
If $\coll{X}$ is a~set of subsets of $E$ and $R$ is a~set of closure rules on~$E$,
then we say that $\coll{X}$ \notion{obeys}~$R$, ot that $R$~\notion{is obeyed by}~$\coll{X}$,
if every set in $\coll{X}$ obeys every rule in~$R\tinysp$.

The relation ``$\negtinysp X$ obeys $B\negdtinysp\adjoins{}c\tinysp$''
	between a~subset $X$ of $E$ and a~closure rule $B\negdtinysp\adjoins{}c$ on~$E$
gives rise to a~Galois connection
\begin{equation*}
\pair{\sigma,\varrho}
	\wide\colon \Pow\negtinysp\ClRul(E) \wide\Galoisconn (\PowPowE)^\op~,
\end{equation*}
where for every $R\subseteq\ClRul(E)$,
	$\sigma(R)$ is the set of all subsets of $E$ that obey~$R\tinysp$,
and for every $\coll{X}\subseteq\PowE$,
	$\varrho\tinysp(\coll{X})$ is the set of all closure rules on $E$ obeyed by~$\coll{X}$.

It is easy to verify that $\sigma(R)$~is a~closure system in the complete lattice~$\PowE$
	for every~$R\subseteq\ClRul(E)$,
and that $\varrho\tinysp(\coll{X})$~is a~closure theory on~$E$
	for every~$\coll{X}\subseteq\PowE\tinysp$.
The converse is also true:
every closure system $\coll{C}$ in $\PowE$ is of the form $\sigma(R)$
	for some $R\subseteq\ClRul(E)$
	(in particular for $R=\varrho(\coll{C})$),
and every closure theory $T$ on $E$ is of the form $\varrho\tinysp(\coll{X})$
	for some $\coll{X}\subseteq\nolinebreak\PowE$
	(in particular for $\coll{X}=\sigma(T)$).
In short: $\sigma\bigl(\Pow\negtinysp\ClRul(E)\bigr) = \ClSys(\PowE)$
	and $\varrho\bigl(\PowPowE\bigr) = \ClTheor(E)\tinysp$.

Given a~set $R$ of closure rules,
we shall say that the closure system $\sigma(R)$ in~$\PowE$ is~\notion{determined by~$R\tinysp$},
and also that the closure operator on $\PowE$ which has $\dtinysp\fix(\gamma)=\sigma(R)$
	is~\notion{determined by~$R\tinysp$}.

The restriction $\ClTheor(E)\to\ClSys(\PowE)$ of~$\sigma$
	is an antiisomorphism of complete lattices
whose inverse $\ClSys(\PowE)\to\ClTheor(E)$
	is the restriction of~$\varrho\tinysp$.

The isomorphism $\Cl(\PowE) \to \ClSys(\PowE)^\op : \gamma \mapsto \fix(\gamma)$
composes with the~iso\-morphism~$\ClSys(\PowE)^\op \to \ClTheor(E)$
to yield the isomorphim of complete lattices $\Cl(\PowE) \to \ClTheor(E)$
which sends a~closure operator $\gamma$ on $\PowE$
to the closure theory
	$\dtinysp\theor(\gamma)
	    \defeq \bigset{B\negdtinysp\adjoins{}c\bigsuchthat c\in\negtinysp\gamma(B)}$
		on~$E\tinysp$.
The~isomorphism $\ClTheor(E) \to \ClSys(\PowE)^\op$
composes with the isomorphism
	$\ClSys(\PowE)^\op \to \Cl(\PowE) : \coll{C}\mapsto\cl_{\coll{C}}$,
yieding the isomorphism of complete lattices $\ClTheor(E) \to \Cl(\PowE)$,
which is the inverse of the isomorphism $\dtinysp\theor\colon\Cl(\PowE) \to \ClTheor(E)\tinysp$,
and which sends a~closure theory~$T$ on~$E$ to the closure operator~$\dtinysp\clop(T)$ on~$\PowE$
given by $\dtinysp\clop(T)(X)
		= \bigset{y\in E\bigsuchthat T\colon\narrdt X\adjoins{}y}$ for $X\subseteq E\tinysp$.

\section{The complete lattice of closure operators on a~dcpo}
\label{sec:completlatt-of-clopers-on-dcpo}

Let~$P$ be a~dcpo,
and let $M$~be the ordered monoid of all preclosure maps on~$P\negtinysp$.
In~$M$~all~directed joins exist, and they are calculated pointwise:
if~$F$ is a~directed subset of~$M$,
then at each $x\in P$ the set $F(x)$ is directed,
thus the map $\varphi\colon P\to P : x\mapsto\tJoin\!F(x)$ is well defined,
and one easily verifies that it is a~preclosure~map;
it follows that $\varphi = \tJoin\!F$ in the poset~$M$.
The ordered monoid~$M$ is therefore a~dcpo;
moreover, $M$~is a~pointed~dcpo since the identity map $\id_P$ is its least element.
Mark that every submonoid of~$M$ is a~directed subset~of~$M$
	because $f,\tinysp g\leq fg$ for any $f,\tinysp g\in M$:
$f\narrt\leq fg$ because $g$~is ascending and $f$~is increasing,
	and~$g\narrt\leq fg$ because $f$~is ascending.

\txtskip

The following theorem describes
how a set of preclosure maps on a~dcpo determines the least~closure operator above~it.

\thmskip

\begin{thm}\label{thm:in-dcpo-cl-generd-by-precls-&-induct}
Let\/ $P$ be a~dcpo,
and let\/ $G$ be a~set of preclosure maps on\/~$P\negtinysp$.
Then\/ $\dtinysp\fix(\leftt G)$ is a~closure system in\/~$P\negtinysp$,
and the corresponding closure operator\/ $\gencl{G}$ on\/~$P\negtinysp$,
	which has\/ $\dtinysp\fix(\gencl{G}\tinysp) = \fix(\leftt G)$,
is the least of all closure operators on\/~$P$ that are above\/~$\nolinebreak G$.

The closure operator\/ $\gencl{G}$ satisfies the~{\bfseries induction principle}:
if a~subset of\/~$P$ is closed under directed joins 
	and is~closed~under\/~$G$,
then~it~is closed under\/~$\gencl{G}$.

Moreover, the~{\bfseries obverse induction principle} for\/~$\gencl{G}$ holds:
if~a~subset of\/ $P$ is inaccessible by directed joins and is inversely closed under\/~$G$,
then it is inversely closed~under\/~$\gencl{G}$.
\end{thm}

\interskip

\begin{myproof}
In the ordered monoid $M\defeq\Precl(\negtinysp P)$
let $H$ be the intersection of all \dircl\ submonoids that include~$G$;
$H$~is the least such submonoid.
Since~$H$ is a~directed subset of the dcpo~$M$ the join $h = \tJoin\!H$ in~$M$ exists,
and $h\in H$ because $H$~is \dircl\ in~$M$,
thus $h$~is the~greatest element~of~$H$.
Since~$H$~is a~submonoid of~$M$, we have $h\tinysp h\in H$, hence~$h\tinysp h\leq h$,
so $h$~is a~closure operator.
We~have~$h\geq G$ because~$h$ is the greatest element of~$H\supseteq G$.
Let~$k\geq G$ be a~closure operator; then~$k\in M$.
The set $K = \allbelow{M\negdtinysp}{k}$ is a submonoid of~$M$
	since~$\id_P\in K$
	and since~$f,\tinysp \fpr \in K$ implies $f\negtinysp\fpr \leq k\tinysp k = k\dtinysp$;
also, $K$~is closed under directed joins in~$M$
	(because it is closed under all existing joins~in~$M$)
	and it~includes~$G$,
thus it includes~$H$, whence~$h\leq k$.

\emph{The induction principle.}
Let $A$ be a subset of~$P$ that is \dircl\ and is closed under~$G$.
Let~$F$ be the set of all $f\in M$ such that $f(A)\subseteq A\tinysp$.
Then $F$~is a~submonoid~of~$M$ and includes~$G\tinysp$.
Also, $F$~is \dircl\ in~$M$
	because~$A$ is \dircl\ in~$P$
	and because the directed joins in~$M$ are calculated pointwise.
It~follows that $H\subseteq F$, hence~$h\in F$, that~is, $h(A)\subseteq A\tinysp$.
This proves the induction principle for the closure operator~$h\tinysp$.

If $g\in G$ then $\dtinysp\fix(g)\supseteq\fix(h)$
	because $g\leq h$ with~$g$ a~preclosure map and~$h$ a~closure operator;
this yields the inclusion $\dtinysp\fix(\leftt G)\supseteq\fix(h)$.
Conversely, if~$a$~is a~fixed point~of~$G$,
then $\set{a}$~is closed under~$G$ and is evidently \dircl,
so it is closed under~$h$ by the induction principle,
therefore~$h(a) = a \in \fix(h)$.

The closure operator $\gencl{G} \defeq h$ has the properties
	stated in the first assertion of the proposition.

\emph{The obverse induction principle.}
Suppose that a~\dirinacc\ subset $A$ of~$P$ is inversely closed under~$G$.
Let $F$ be the set of all $f\in M$ such that $f^{-1}(A)\subseteq A\dtinysp$;
$F$~includes~$G$ and it is a~submonoid of~$M$.
Let~$E$ be a~directed subset~of~$F$; we shall show that $\tJoin\!E\in F$.
Let~$x\in P\negtinysp$,
and suppose that $\bigl(\tJoin\!E\bigr)(x) = \tJoin\negdtinysp E(x) \in\nolinebreak A\dtinysp$;
since~$A$ is \dirinacc\ there exists $e\in E$ with $e(x)\in A\tinysp$,
and we have $x\in A$ because $A$ is inversely closed under~$e\tinysp$.
It~follows that $A$ is inversely closed under~$\tJoin\!E$.
We see that~$F$ is \dircl,
so~$F$~includes~$H$ and therefore~contains the closure operator~$h=\gencl{G}$,
whence $A$~is inversely closed under~$\gencl{G}$.
\end{myproof}

\thmskip

In the classical logic, which uses the~law of~excluded middle with abandon,
the obverse induction principle for a~subset~$A$ of~$P$
is just a~rephrasing of the induction principle for the complement $P\setdiff A\tinysp$.
Since we want to apply the obverse induction principle in situations
	where the~law of~excluded middle is not admissible,
we~proved~it on~its~own.

\txtskip

Let $P\negtinysp$, $G$, and $\gencl{G}$ be as in Theorem~\ref{thm:in-dcpo-cl-generd-by-precls-&-induct}.
We~shall say that the closure operator~$\gencl{G}$
	is~\notion{generated by} the set~$G$ of preclosure maps.

\txtskip

The special case of Theorem~\ref{thm:in-dcpo-cl-generd-by-precls-&-induct} where $G=\set{g}$
is of interest on its own.

\thmskip

\begin{cor}\label{cor:cl-gen-by-single-precl-in-dcpo}
If\/ $g$ is a~preclosure map on a~dcpo\/~$P\negtinysp$,
then\/ $\dtinysp\fix(g)$~is a~closure system in\/~$P\negtinysp$,
and the closure operator\/~$\overbar{g}$ on\/~$P$ that has\/ $\dtinysp\fix(\overbar{g})=\fix(g)$
is the least of all closure operators on\/~$P$ that are above\/~$g\tinysp$.
\qed
\end{cor}

\thmskip

The induction principle and the obverse induction principle
of course hold in the special case featuring a~single preclosure map;
there is no need to restate them.

\txtskip

As a~consequence of Theorem~\ref{thm:in-dcpo-cl-generd-by-precls-&-induct},
if~$P$ is a~dcpo, then in the poset $\Cl(\negtinysp P)$ every subset has a~join,
	therefore~$\Cl(\negtinysp P)$ is a~complete lattice.

\thmskip

\begin{cor}\label{cor:Cl(P)-is-complete-lattice}
The poset\/ $\Cl(\negtinysp P)$ of all closure operators on a~dcpo\/~$P$
is a~complete lattice,
and so is the poset\/ $\ClSys(\negtinysp P)$ of all closure systems in\/~$P\negtinysp$.
In\/~$\Cl(\negtinysp P)$, the join of a~set\/~$G$ of closure operators is
the closure operator\/~$\gencl{G}$ generated by\/~$G$,
while~in\/~$\ClSys(\negtinysp P)$, the meet of a~set\/~$\coll{C}$ of closure systems
	is the intersection\/ $\tInters\coll{C}$.
\end{cor}

\interskip

\begin{myproof}
For every $G\subseteq\Cl(\negtinysp P)$,
$\gencl{G}$ is the join of~$G$ in~$\Cl(\leftt G)$
	by Theorem~\ref{thm:in-dcpo-cl-generd-by-precls-&-induct}.
The~map\-ping~$\Cl(\negtinysp P)\to\ClSys(\negtinysp P) : g\mapsto\fix(g)$
	is an~antiisomorphism of~complete lattices,
so for every $G\subseteq\Cl(\negtinysp P)$ we have
    $\tMeet_{g\in G}\fix(g)
	= \fix(\tJoin\negdtinysp G)
	= \fix\bigl(\gencl{G}\dtinysp\bigr)
	= \fix(\leftt G)
	= \tInters_{\tinysp g\in G}\fix(g)$,
therefore
all~meets in $\ClSys(\negtinysp P)$ exist and they are calculated as~intersections.
\end{myproof}

\thmskip

Let $P$ be a~dcpo.
The set $\ClSys(\negtinysp P)$ of all closure systems in~$P$
is a~closure system in the complete lattice $\PowP$.
The corresponding closure operator on~$\PowP$
maps each subset~$X$ of~$P$
	to the closure system~$\dtinysp\clsys(X) = \clsys_P(X)\tinysp$,
which is the least of all closure systems in~$P$ that include~$X$.

\txtskip

Tarski's fixed point theorem, a~version for dcpos,
	easily follows from Theorem~\ref{thm:in-dcpo-cl-generd-by-precls-&-induct}.

\thmskip

\begin{thm}\label{thm:tarski's-fixed-point-thm-for-dcpos}
Let\/ $f$ be an increasing map on a~dcpo\/~$P\negtinysp$.
The subposet\/~$\dtinysp\fix(f)$ of\/~$P$ is~a~dcpo.
For every\/ $x\in P$ on which\/ $f$~ascends
there exists the least fixed point of\/~$f$ above\/~$x$.
If\/~$P$~has a~least element, then\/~$f$~has a~least fixed~point.
\end{thm}

\interskip

\begin{myproof}
Let $A$ be the set of all elements of $P$ on which $f$ ascends;
clearly $A$~contains all fixed points of~$f$.
Since $x\leq f(x)$ implies $f(x)\leq f(f(x))$,
the set $A$ is closed under~$f$.
If $D\subseteq A$ is directed, then the join $\tJoin\!D$ in~$P$ exists,
and for every $d\in D$ we have $d\leq f(d)\leq f(\tJoin\!D)$,
whence $\tJoin\!D\leq f(\tJoin\!D)$,
so~$A$~is closed under directed~joins.
Thus the subposet $A$ is a~dcpo,%
\footnote{\,In fact the subposet $A$ is a~sub$\tinysp$-dcpo of~$P\negtinysp$,
	since the directed joins in~$A$ are inherited from~$P\negtinysp$.}
and the restriction $g\colon A\to A$ of~$f$ is a~preclosure map on~$A\tinysp$.
The closure operator $\gencl{g}$ on~$A$ generated by the preclosure map~$g$ on~$A$
has the fixpoint set $\dtinysp\fix(\gencl{g}) = \fix(g) = \fix(f)$,
which is a~closure system in the dcpo~$A\tinysp$,
and so as a~subposet of~$A\tinysp$, and hence~of~$P\negtinysp$,
	it~is~itself a~dcpo.%
\footnote{\,The subposet $\dtinysp\fix(f)$ of~$P$ is,
	in general, not a~sub$\tinysp$-dcpo of~$P\negtinysp$.}
If~$P$ has a~least element~$\bot$,
then $\bot\in A\tinysp$,
and~$\gencl{g}(\bot)$~is the least element of~$\gencl{g}(A)=\fix(\gencl{g})=\fix(f)$.
\end{myproof}

\thmskip

The last statement of Theorem~\ref{thm:tarski's-fixed-point-thm-for-dcpos}
is the bare-bones Tarski's fixed point theorem for~dcpos;
let us restate it on its own.

\thmskip

\begin{cor}\label{cor:barebones-tarski's-fixed-point-thm}
Every increasing endomap on a~pointed dcpo has a~least fixed point. \qed
\end{cor}

\thmskip

It can be proved, with a~generous help from the~axiom of~choice,
that the bare-bones Tarski's fixed point property, stated in the corollary,
	in fact characterizes pointed dcpos.
See,~for example, Theorem~11
	(after consulting Corollary~2 of Theorem~1) in~\cite{Markowsky-CCP&DSA-1976}.

\txtskip

We round off this section by presenting a~class of curious dcpos.
Each~of~these dcpos is associated with an arbitrary poset in a~rather peculiar way.
The challenging part of the presentation is a~proof that what is offered is in fact a~dcpo.
At~the~crucial point of the proof it is Corollary~\ref{cor:cl-gen-by-single-precl-in-dcpo}
	that gets us over the hurdle, with a~flick of a~finger.

\txtskip

Let $P$ be a~poset.

We shall say that $A$ is a~\notion{directed-complete subposet of~$P$},
or a~\notion{dc-subposet of~$P$},
if $A$ is a~subposet of $P$ that is a~dcpo (with respect to the induced ordering).
We~shall denote by $\DcSpos(\negtinysp P)$
	the poset of all dc-subposets of~$P$ ordered by inclusion,
and by $\DcSpos_\bot(\negtinysp P)$ the subposet of $\DcSpos(\negtinysp P)$
	consisting of all pointed dc-subposets (\ie, dc-subposets that possess a~least element) of~$P$.

If~$A$~and~$B$ are dc-subposets of~$P$ and $A\supseteq B\supseteq Y$ with $Y$ directed,
	then $\stJoin^A\negdtinysp Y \leq \stJoin^B\negdtinysp Y$;
indeed, $\stJoin^B\negdtinysp Y$ is an element of~$A$
	and it is an upper bound of the subset~$Y$ of~$A\tinysp$,
so it is above the least upper bound $\stJoin^A\negdtinysp Y$ of~$Y$ in~$A\tinysp$.

\thmskip
\pagebreak[3]

\begin{prop}\label{prop:Dcsps(P)-is-dcpo}
Let\/~$P$~be a~poset.
If\/~$\coll{F}$ is a~filtered set of dc-subposets of\/~$P$,
	then\/~$\tInters\coll{F}$~is a~dc-subposet of\/~$P$;
moreover, if all dc-subposets in\/~$\coll{F}$ are pointed,
	then the dc-subposet\/~$\tInters\coll{F}$ is pointed.
In~other words,
$\DcSpos(\negtinysp P)^\op$ and $\DcSpos_\bot(\negtinysp P)^\op$
	are sub-dcpos of the dcpo\/ $(\Pow{\negtinysp P})^\op$.
\end{prop}

\interskip

\begin{myproof}
Let $\coll{F}$ be a~filtered subset of~$\DcSpos(\negtinysp P)$.
We shall prove that $\tInters\coll{F}\in\DcSpos(\negtinysp P)$. %

Let $F$ be the product of the sets in~$\coll{F}$.
The elements of $F$ are the functions $f\colon\coll{F}\to P$
	that have $f(A)\in A$ for every $A\in\coll{F}\tinysp$;
that is, they are the choice functions for the set~$\coll{F}$ of~sets.
We~order $F$ by the pointwise ordering.
Endowed with~this ordering,
$F$~becomes a~dcpo whose directed joins are calculated pointwise:
if~$H$~is a~directed subset of~$F$,
then $H(A)$ is a~directed subset of~$A$ for every $A\in\coll{F}$,
and
\begin{equation*}
(\stJoin^F\!\! H)(A) \Eq \stJoin^A\!H(A)~, \qquad A\in\coll{F}\,.
\end{equation*}
Let $E$ be the set of all $f\in F$
	that are increasing on the directed set $\coll{F}^{\tinysp\op}$;
that~is, $f\in E$ iff
	for all $A,\,B\in \coll{F}$ the inclusion $A\supseteq B$ implies $f(A)\leq f(B)$.
The subset $E$ of $F$ is a~sub-dcpo of the dcpo~$F\tinysp$:
if $H$ is a~directed subset of~$E$, and $A\supseteq B$ are in~$\coll{F}$,
then
\begin{equation*}
(\stJoin^F\!\!H)(A)
	\Eq \stJoin^A\!H(A)
	\Leq \stJoin^A\!H(B)
	\Leq \stJoin^B\!H(B)
	\Eq (\stJoin^F\!\!H)(B)~,
\end{equation*}
therefore $\stJoin^F\!\!H\in E$.
The poset~$E$, with the pointwise ordering,
	is a~dcpo in which the directed joins are calculated pointwise.

For each $f\in E$ define $J\negdtinysp f\in F$ by
\begin{equation*}
(\leftt{J\negdtinysp f})(A) \Defeq \stJoin^A \bigset{f(X)\bigsuchthat
					X\in \allbelow{\coll{F}\negdtinysp}{\negdtinysp A}}~,
	\qquad A\in\coll{F}\;;
\end{equation*}
the join exists because
	the set $(\allbelow{\coll{F}\negdtinysp}{\negdtinysp A})^\op$ is directed
	and the map $(\allbelow{\coll{F}\negdtinysp}{\negdtinysp A})^\op \to A : X\mapsto f(X)$
		is increasing,
so the set under the join is a~directed subset of~$A\tinysp$.
We claim that $J\negdtinysp f\in E\tinysp$.
In order to~prove this, let $A\supseteq B$ be sets~in~$\coll{F}$.
First note that~$\stJoin^A\bigset{f(X)\bigsuchthat X\in \allbelow{\coll{F}\negdtinysp}{\negdtinysp B}}
	= \stJoin^A\bigset{f(X)\bigsuchthat X\in \allbelow{\coll{F}\negdtinysp}{\negdtinysp A}}$
because $(\allbelow{\coll{F}\negdtinysp}{\negdtinysp B})^\op$
is a~cofinal subset of~$(\allbelow{\coll{F}\negdtinysp}{\negdtinysp A})^\op$.
Then we have
\begin{align*}
(\leftt{J\negdtinysp f})(A)
	& \Eq \stJoin^A\bigset{f(X)\bigsuchthat X\in \allbelow{\coll{F}\negdtinysp}{\negdtinysp A}}
		\Eq \stJoin^A\bigset{f(X)\bigsuchthat
					X\in \allbelow{\coll{F}\negdtinysp}{\negdtinysp B}} \\ 
	& \Leq \stJoin^B\bigset{f(X)\bigsuchthat X\in \allbelow{\coll{F}\negdtinysp}{\negdtinysp B}}
		\Eq (\leftt{J\negdtinysp f})(B)\,.
\end{align*}
Thus $J$ maps $E$ to~$E$.
It~is clear that the endofunction $J$ on $E$
is ascending ($f\leq J\negdtinysp f$ for every $f\in E$)
and increasing ($f\leq g$ implies $J\negdtinysp f\leq Jg$, for all $f,\,g\in E$).
Since $E$ is a~dcpo,
it~follows by Corollary~\ref{cor:cl-gen-by-single-precl-in-dcpo}
that $G=\fix(J)$~is a~closure system in~$E\tinysp$;
let $\Gamma$ be the closure operator on~$E$ with $\dtinysp\fix(\Gamma)=G$.

Let $g\in G$.
If $A\supseteq B$ are sets in~$\coll{F}$,
then $g(A)\leq g(B)$, but also
\begin{equation*}
g(B) \Leq \stJoin^A\bigset{g(X)\suchthat X\in \allbelow{\coll{F}\negdtinysp}{\negdtinysp A}}
	\Eq (J g)(A)
	\Eq g(A)\,,
\end{equation*}
and therefore $g(A) = g(B)$.
Now, for any two sets $A$ and $B$ in $\coll{F}$
there exists a~set~$C$ in~$\coll{F}$ such that $A\supseteq C$ and $B\supseteq C$,
whence $g(A) = g(C) =\nolinebreak g(B)$.
In~short, $g$~is a~constant function,
with its constant value lying in the intersection $\tInters\coll{F}$.
Conversely, if $u\in\tInters\coll{F}$,
then the function $\widehat{u}\in E$ with the constant value~$u$
is clearly a fixpoint of~$J$ and hence a~fixpoint~of~$\Gamma$.
The mapping $\tInters\coll{F}\to G : u\mapsto\widehat{u}$ is evidently an isomorphism of posets.

Let $S\subseteq G$ be directed.
The set~$S$ has a~join (the pointwise join) $\tJoinS$ in~$E\tinysp$;
but then $\Gamma\tJoinS$ is the join of~$S$ in~$G$.
The poset~$G$ is a~dcpo, and so is then the poset~$\tInters\coll{F}$.

If every $A\in\coll{F}$ has a~least element $\bot_A$,
then the choice function $\bot\colon\coll{F}\to P : A\mapsto\bot_A$
	increases on $\coll{F}^{\tinysp\op}$,
therefore $\bot$ is the least element of~$E$,
and so~$\Gamma\bot$ is the least element of~$G$,
	while the constant value of $\Gamma\bot$ is the least element of~$\tInters\coll{F}$.
\end{myproof}

\section{\Scottcont\ closure operators on dcpos}
\label{sec:Scottcont-clopers-on-dcpo}

To begin with we establish a~general result
about preservation of some special joins
	by the~pointwise join of functions that preserve those special joins.

\txtskip

Recall that for any poset $P$ we denote by $\Joinable\negtinysp P$
	the set of all joinable subsets of~$P$.

Let $P$ and $Q$ be posets,
	and let $\coll{A}$~be a~subset of $\Joinable\negtinysp P$.

We~shall say that a~function $f\colon P\to Q$ \notion{preserves $\coll{A}$-joins}
    if for every $A\in\coll{A}$ we have
	$f(A)\in\Joinable Q$ and $\tJoinfof(A)=f(\tJoinA)\tinysp$.
Let~$F$ be a~set of functions $P\to Q\tinysp$.
We~shall say that $F$~\notion{preserves $\coll{A}$-joins}
	if~every function in $F$ preserves $\coll{A}$-joins.
We~shall say that $F$~is \notion{pointwise-joinable}
	if~$F(x)\in\Joinable Q$ for every $x\in P\negtinysp$.
Whenever $F$ is pointwise\nobreakdash-joinable
	we define the~\notion{pointwise join} $\dottJoinF\colon P\to Q$ of~$F$
by $(\dottJoinF)(x) \defeq \tJoinFof(x)$ for~$x\in P$. %

If $\coll{A}$ contains all subsets $\set{x,y}$ of $P$ with $x\leq y\tinysp$,
then every function $P\to Q$ which preserves $\coll{A}$-joins is increasing.

\thmskip

\begin{prop}\label{prop:preserv-of-joins-by-pntwise-join}
Let\/ $P$ and\/ $Q$ be posets, let\/ $F$~be a~set of functions\/ $P\to Q\tinysp$,
	and let\/ $\coll{A}$~be a subset of\/ $\Joinable\negtinysp P$.
If\/ $F$ preserves\/ $\coll{A}$-joins and is pointwise-joinable,
then the pointwise join\/ $\dottJoinF$ preserves\/ $\coll{A}$-joins;
moreover, for every\/ $A\in\coll{A}$ we have\/ $F(A)\in\Joinable Q$ and
\begin{equation}\label{eq:(pwJoin F)(Join A)=Join(pwJoin f)(A)=Join F(A)}
    (\dottJoinF)(\tJoinA)
	\Eq \tJoin\dtinysp(\dottJoinF)(A)
	\Eq \tJoinFof(A)~.
\end{equation}
\end{prop}

\negdisplayhalfskip
\interskip

\begin{myproof}
Consider any $A\in\coll{A}\tinysp$.
By assumption the join $\tJoinA$ exists
	and the join $\tJoinFof(x)$ exists for every~$x\in P\negtinysp$,
therefore the element
	$(\dottJoinF)(\tJoinA) = \tJoinFof(\tJoinA)$
		of~$Q$
and the subset
	$(\dottJoinF)(A) = \bigset{\tJoinFof(a)\bigsuchthat a\in A}$
		of~$Q$
are well-defined.
Let $y$ be an arbitrary element of $Q$.
The chain of equivalences
\begin{align*}
y \geq (\dottJoinF)(\tJoinA)
    &\wider\Isequiv \text{$y \geq f(\tJoinA)$\, for every \,$f\narrt\in F$} \\
    &\wider\Isequiv \text{$y \geq \tJoinfof(A)$\, for every \,$f\narrt\in F$}
	\quad \text{(since $f(\tJoinA) = \tJoinfof(A)$)}\\
    &\wider\Isequiv \text{$y\geq f(a)$\, for every \,$a\in A$\, and for every\, $f\narrt\in F$} \\
    &\wider\Isequiv \text{$y\geq F(A)$} 
\end{align*}
proves that $F(A)\in\Joinable Q$ and that
	$\tJoin\negdtinysp F(A) = (\dottJoinF)(\tJoinA)$.
Then
\begin{align*}
y \geq (\dottJoinF)(A)
    &\wider\Isequiv \text{$y \geq \tJoinFof(a)$\, for every \,$a\in A$} \\
    &\wider\Isequiv \text{$y \geq f(a)$\, for every \,$f\narrt\in F$\, and for every \,$a\in A$} \\
    &\wider\Isequiv \text{$y\geq F(A)$} 
\end{align*}
clinches the proof ot the equalities~\eqref{eq:(pwJoin F)(Join A)=Join(pwJoin f)(A)=Join F(A)}.
\end{myproof}

\thmskip

By the definition of $\dottJoinF$
we may add $\tJoinFof(\tJoinA) = (\dottJoinF)(\tJoinA)$
to the equalities \eqref{eq:(pwJoin F)(Join A)=Join(pwJoin f)(A)=Join F(A)}.

\txtskip

Proposition \ref{prop:preserv-of-joins-by-pntwise-join} is so general with a~reason:
it makes perfectly clear that
	$F$ preserving $\coll{A}$\nobreakdash-joins and $F$ being pointwise-joinable
are two independent properties of the~set~$F$ of functions;
in a~sense these two properties are orthogonal to each other.
For~example, suppose that $P$ and $Q$ are dcpos
	and that $\coll{A}$ is the set of all directed subsets~of~$P\tinysp$;
in this case $F$ preserving $\coll{A}$-joins means that $F$ preserves directed joins.
The~set~$F$ of functions need not be directed (as long as it is pointwise-joinable).
It~surely helps if~$F$ is directed, since then the sets $F(x)$, $x\in P\negtinysp$,
are directed subsets of~$Q$ and have joins in~$Q$,
thus we know that $F$ is pointwise-joinable \emph{because} it is directed and $Q$~is a~dcpo,
and we can conclude that the pointwise join $\dottJoinF$ preserves directed joins.
But suppose that $Q$ is a~complete lattice, with~$P$ still just any dcpo:
then \emph{every} set~$F$ of functions $P\to Q$ which preserves directed joins
has the~pointwise join $\dottJoinF$ which preserves directed~joins.

\txtskip

A function between posets $f\colon P\to Q$ is said to be \notion{\Scottcont}
if it preserves all existing directed joins.%
\footnote{\,In other words, $f$ is \Scottcont\ iff it preserves $\Joinable\negtinysp\DirSubsets{P}$-joins.}
In detail, $f$ is \Scottcont\ \iff\
for~every joinable directed subset $D$ of $P$
the $f$-image of the join of~$D$ in~$P$ is the join of~$f(D)$ in~$Q$,%
\footnote{\,There are those who prefer the more long-winded
``$\tinysp$the $f$-image of the set $D$ has a~join in $Q$
	which is equal to the $f$-image of the join of~$D$ in~$P\tinysp$''.}
that~is, $f(\tJoinD) = \tJoinfof(D)\tinysp$.
In~particular, a~\Scottcont\ function $f$
	preserves joins of all pairs of comparable elements~of~$P\negtinysp$,
which implies that $f$~is~increasing,
therefore for every joinable directed subset~$D$~of~$P$
	its image~$f(D)$~is a~joinable \emph{directed} subset of~$Q$.
From this it follows that if $f\colon P\to Q$ and $g\colon Q\to R$
are \Scottcont\ functions between posets,
then the composite function $gf\colon P\to R$ is \Scottcont.

For any posets $P$ and $Q$ we let $\Sc(P,Q)$
denote the poset of all \Scottcont\ functions $P\to Q$ with the pointwise ordering;
$\Sc(P,Q)$~is a~subposet of~$\Inc(P,Q)\tinysp$.

For any poset $P$ we denote by $\ScPrecl(\negtinysp P)$ the pointwise-ordered poset
	of all \Scottcont\ preclosure maps on~$P\negtinysp$,
and by $\ScCl(\negtinysp P)$ the pointwise-ordered poset of all \Scottcont\ closure operators on~$P\tinysp$;
$\ScPrecl(\negtinysp P)$~is an~ordered submonoid of~$\Precl(\negtinysp P)\tinysp$.

\txtskip

The following proposition follows from Proposition~\ref{prop:preserv-of-joins-by-pntwise-join}
by specialization.

\thmskip

\begin{prop}\label{prop:dirjoin-of-Scottcont-funs-between-dcpos}
Let\/ $P$ and\/ $Q$ be dcpos.
If\/ $F$ is a~directed subset of\/ $\Sc(P,Q)$,
then the pointwise join\/ $\dottJoinF$ exists and is \Scottcont.
\qed
\end{prop}

\interthmskip

\begin{cor}\label{cor:P-Q-dcpos==>Sc(P,Q)-dcpo-pntwise-joins}
If\/ $P$ and\/ $Q$ are dcpos,
then\/ $\Sc(P,Q)$ is a~dcpo in which directed joins are calculated pointwise.
\qed
\end{cor}

\interthmskip

\begin{cor}\label{cor:P-dcpo==>ScPreCl(P)-dcpo-pntwise-joins}
If $P$ is a~dcpo,
then $\ScPrecl(\negtinysp P)$ is a~pointed dcpo in which directed joins are calculated pointwise.
\qed
\end{cor}

\interskip

\begin{myproof}
If $F$ is a~directed set of \Scottcont\ preclosure maps on~$P\negtinysp$,
then the pointwise join $\dottJoinF$ is a~preclosure map
which is \Scottcont\ by~Proposition~\ref{prop:dirjoin-of-Scottcont-funs-between-dcpos}.
Therefore $\ScPrecl(\negtinysp P)$ is a~dcpo in which directed joins are calculated pointwise,
and it is a~pointed dcpo since the identity map $\widedt{\id_P}$ is its bottom element.
\end{myproof}

\thmskip

Let $P\negtinysp$, $Q$, $R$ be posets.
For any sets of functions $F\subseteq\Inc(P,Q)$ and $G\subseteq\Inc(Q,R)$
we write
	$GF \narrt\defeq \set{\dtinysp gf\negdtinysp\suchthat
					f\narrdt\in F,\, g\narrdt\in G\tinysp}\dtinysp$;
when both $F$ and $G$ are directed, $GF$ is easily seen to be directed
	(recall that composition of increasing maps is increasing in both operands).

\thmskip

\begin{prop}\label{prop:composing-Sc(P,Q)-and-Sc(Q,R)}
Let\/ $P\negtinysp$, $Q$, and\/ $R$ be dcpos.
If\/ $F\subseteq\Sc(P,Q)$ and\/ $G\subseteq\Sc(Q,R)$ are directed,
then\/ $GF \subseteq \Sc(P,R)$ is~directed and
\begin{equation}\label{eq:(Join(G))*(Join(F))=Join(G*F)}
(\tJoin\!G)(\tJoinF) \Eq \tJoin(GF)~.
\end{equation}
\end{prop}

\interskip

\begin{myproof}
Let $x\in P\negtinysp$.
The set $F(x)$ is directed,
every $g\in G$ preserves directed joins,
and~$G$, being directed, is pointwise-joinable,
and so we calculate:
\begin{align*}
\bigl((\tJoinG)(\tJoinF)\bigr)(x)
	&\Eq (\dottJoinG)\bigl((\dottJoinF)(x)\bigr)
	    \Eq (\dottJoinG)\bigl(\tJoinFof(x)\bigr)  \\
	&\Eq \tJoinGof\bigl(F(x)\bigr)
		\qquad\text{(by Proposition~\ref{prop:preserv-of-joins-by-pntwise-join})} \\
	&\Eq \tJoin\dtinysp(GF)(x)
	    \Eq \bigl(\dottJoin(GF)\bigr)(x) \\
	&\Eq \bigl(\tJoin(GF)\bigr)(x)~.
\end{align*}
This proves the equality~\eqref{eq:(Join(G))*(Join(F))=Join(G*F)}.
\end{myproof}

\interthmskip

\begin{cor}\label{cor:dcpo-P==>compos-monoids-Sc(P)-and-ScPreCl(P)}
Let\/ $P$ be a~dcpo.
The dcpo\/ $\Sc(\negtinysp P)$ is an~ordered composition monoid
	in which composition distributes over directed joins:
if\/~$F$ and\/~$G$ are directed subsets of\/~$\Sc(\negtinysp P)\tinysp$,
then\/ $GF$~is a~directed subset of\/~$\Sc(\negtinysp P)$ and
	$(\tJoinG)(\tJoinF) = \tJoin(GF)\tinysp$.%
\footnote{\,On the level of directed subsets of $\Sc(\negtinysp P)$
	the identity $\tJoin(GF)=(\tJoinG)(\tJoinF)$ says that joining distributes over composition.
But on the level of elements of $\Sc(\negtinysp P)$
	it is the composition that distributes over directed joins,
which becomes apparent when we rewrite the identity as
	$(\tJoin_{\!g\in G}g)(\tJoin_{\!f\in F}f)=\tJoin_{\!g\in G,\,f\in F}gf$.}
Likewise the pointed dcpo\/~$\ScPrecl(\negtinysp P)$ is an~ordered composition monoid
	in which composition distributes over directed joins.
\qed
\end{cor}

\thmskip

Let $P$ be a~dcpo.
Every submonoid of the ordered composition monoid $\ScPrecl(\negtinysp P)$
	is a~submonoid of the ordered composition monoid $\Precl(\negtinysp P)$
and is therefore directed.

Recall that for any subset $X$ of a~monoid~$M$
	we denote by $X^*$ the submonoid of~$M$ generated by~$X$.

\thmskip

\begin{thm}\label{thm:in-dcpo-Scottcont-preclmaps-gener-Scottcont-clop}
Let\/ $P$ be a~dcpo, and\/ let $G$ be a~subset of\/~$\ScPrecl(\negtinysp P)$.
The {\rm(}pointwise calculated\/{\rm)}
	directed
	join\/~$h\defeq\tJoinG^*$ in\/~$\ScPrecl(\negtinysp P)$
exists and is a~\Scottcont\ closure operator on\/~$P\negtinysp$.
Moreover, $h$ is the closure operator on the dcpo\/~$P$
	generated by the set\/~$G$ of preclosure maps~on\/~$P\negtinysp$,
therefore\/ $\dtinysp\fix(h) = \fix(\leftt G)\tinysp$.
\end{thm}

\interskip

\begin{myproof}
Since the submonoid $G^*\negdtinysp$ of $\tinysp\ScPrecl(\negtinysp P)$ is directed,
the join
    $h = \tJoinG^*\negdtinysp
	= \dottJoinG^*\negdtinysp$
		in $\ScPrecl(\negtinysp P)$ exists.
According to Corollary~\ref{cor:dcpo-P==>compos-monoids-Sc(P)-and-ScPreCl(P)}
we have
\begin{equation*}
h\tinysp h \Eq (\tJoinG^*)(\tJoinG^*)
	\Eq \tJoin(G^*G^*)
	\Eq \tJoinG^*\negdtinysp
	\Eq h~,
\end{equation*}
thus~$h$~is a~\Scottcont\ closure operator on~$P\negtinysp$, and clearly $h\geq G$.

Now let $k$ be a~closure operator on $P$ and $k\geq G$.
(Note that we are \emph{not} assuming that $k$ is \Scottcont.)
First, $k\geq\id_P\in G^*\negdtinysp$.
Next, if $u\in G^*$ is a~composite of $n\geq 1$ functions in~$G$, then $u\leq k^n = k\tinysp$.
Therefore $k\geq G^*\negdtinysp$, whence
    $k\geq \dottJoinG^*\negdtinysp
	= \tJoinG^*\negdtinysp
	= h\tinysp$.
It~follows that $h$ is the closure operator on~$P$
generated by the set $G$ of preclosure maps on~$P\negtinysp$,
so $h$ is, by~Theorem~\ref{thm:in-dcpo-cl-generd-by-precls-&-induct},
	the closure operator~on~$P$ that has $\dtinysp\fix(h)=\fix(\leftt G)\tinysp$.
\end{myproof}

\thmskip

\begin{cor}\label{cor:single-Scottcont-preclmap-geners-Scottcont-clop}
If\/ $g$ is a~\Scottcont\ preclosure map on a~dcpo\/ $P\negtinysp$,
then\/ $\dtinysp\fix(g)$~is the~fixpoint set of the
	\Scottcont\ closure operator\/ $\dottJoin_{\!k\in\NN} g^{\tinysp k}$,
which is the least of the closure operators on\/~$P$ that are above\/~$g\tinysp$.
\qed
\end{cor}

\thmskip

The nonempty chain $\set{g^k\narrdt\suchthat k\narrt\in\NN}$ appearing in the corollary
	is of course the~composition monoid~$\set{g}^*$
gene\-rated by the~\Scottcont\ preclosure map~$g\tinysp$.

\thmskip

\begin{cor}\label{cor:dcpo-P-joins-in-ScCl(P)}
Let\/ $P$ be a~dcpo, and let\/ $G$ be a~set of \Scottcont\ closure operators on\/~$P\negtinysp$.
The directed pointwise join\/ $\dottJoinG^*\negdtinysp$
is the join\/ $\tJoinG$ of\/ $G$ in\/ $\ScCl(\negtinysp P)\tinysp$,
and~it~is~also the~join of\/ $G$ in\/ $\Cl(\negtinysp P)$ so that\/
    $\dtinysp\fix(\tJoinG) = \fix(\leftt G)\tinysp$.
\qed
\end{cor}

\thmskip

Let $P$ be a~dcpo.
Corollary~\ref{cor:dcpo-P-joins-in-ScCl(P)} tells us that
the the set $\ScCl(\negtinysp P)$ of all \Scottcont\ closure operators
	is~closed under all joins in the complete lattice $\Cl(\negtinysp P)$
and is therefore itself a~complete lattice whose joins are inherited from~$\Cl(\negtinysp P)$.
Correspondingly,
the set $\ScClSys(\negtinysp P)$ of the fixpoint sets of all \Scottcont\ closure operators on~$P$
is closed under all meets in the complete lattice $\ClSys(\negtinysp P)$
	which is a~closure system in~$\PowP$,
and is therefore itself a~closure system in~$\PowP\negtinysp$,
that is, it is closed under arbitrary intersections.
There is a~less roundabout way to see this,
using an explicit characterization of fixpoint sets of the \Scottcont\ closure operators on a~dcpo,
given below in Lemma~\ref{lem:dcpo-P-dirclosed-clsyss}.

\txtskip

But first an~auxiliary lemma, almost trivial, though still worth telling on its own.

\thmskip

\begin{lem}\label{lem:gamma(Join(X))=gamma(gamma(Join(X)))}
Let\/ $P$ be a~poset, $\gamma$ a~closure operator on\/~$P\negtinysp$, and\/ $X$ a~subset of\/~$P\negtinysp$.
If~in\/~$P$ both\/~$\tJoin\!X$ and\/ $\tJoin\negdtinysp\gamma(X)$ exist,
then\/ $\gamma(\tJoin\!X) = \gamma(\tJoin\negdtinysp\gamma(X))\tinysp$.
\end{lem}

\interskip

\begin{myproof}
We get the asserted identity by
applying $\gamma$ to $\tJoin\!X \leq \tJoin\negdtinysp\gamma(X) \leq \gamma(\tJoin\!X)\tinysp$.
\end{myproof}

\thmskip

And here is the promised characterization.

\thmskip

\begin{lem}\label{lem:dcpo-P-dirclosed-clsyss}
A~closure operator\/~$\gamma$ on a~dcpo\/~$P$ is \Scottcont\
	iff\/ $\dtinysp\fix(\gamma)$~is closed under directed joins in\/~$P$.
\end{lem}

\interskip

\begin{myproof}
Write $C\defeq\fix(\gamma)$.

Suppose $\gamma$ is \Scottcont.
If $Y\negdtinysp\subseteq C$ is directed,
then $\gamma(\tinysp\tJoin\negtinysp Y)
	= \tJoin\negdtinysp\gamma(Y)
	=\tJoin\negtinysp Y\negdtinysp$,
therefore $\tJoin\negtinysp Y\negdtinysp\in C$.%
\footnote{\,The necessity part of the lemma holds, mutatis mutandis,
for any \Scottcont\ endomap~$\gamma$ on any poset~$P\dtinysp$:
if~$Y\subseteq\fix(\gamma)$ is joinable directed in~$P$, then $\tJoinY\!\in\fix(\gamma)$.}

Suppose $C$ is closed under directed joins,
and let $Y\negdtinysp\subseteq P$ be directed.
Then $\gamma(Y)$~is a~directed subset of $C$,
thus $\tJoin\negdtinysp\gamma(Y)\in C$,
and $\gamma(\tinysp\tJoin\negtinysp Y)
	= \gamma(\tinysp\tJoin\negdtinysp\gamma(Y))
	= \tJoin\negdtinysp\gamma(Y)$.
\end{myproof}

\thmskip

Let us state yet another~simple\,---\,but handy\,---\,lemma;
	we omit the evident proof.
In~the~lemma we use the following notation:
for any~set~$\coll{S}$ of subsets of a~poset~$P$
we denote by~$\Dc\coll{S}$ the set of all \dircl\ sets belonging to~$\coll{S}$.

\thmskip

\begin{lem}\label{lem:P-poset&C-clsys-in-Pow(P)=>Dc(C)-clsys-in-Pow(P)}
Let\/~$\coll{C}$ be a~set of subsets of a~poset\/~$P$.
If\/~$\coll{C}$ is a~closure system in\/~$\PowP$,
	then\/~$\Dc\coll{C}$ is a~closure system in\/~$\PowP$.
\qed
\end{lem}

\thmskip

Let $P$ be a~dcpo.

We know that the poset $\Cl(\negtinysp P)$ of all closure operators on $P\negtinysp$, ordered pointwise,
and the poset $\ClSys(\negtinysp P)$ of all closure systems on $P\negtinysp$, ordered by inclusion,
are complete lattices,
where $\ClSys(\negtinysp P)$ is a~closure system in the powerset lattice $\PowP\negtinysp$,
meaning that the intersection of any set of closure systems in~$P$
	is a~closure system in~$P\negtinysp$.

Since the set $\ClSys(\negtinysp P)$ of all closure systems in~$P$ is a~closure system in $\PowP$,
it follows by~Lemma~\ref{lem:P-poset&C-clsys-in-Pow(P)=>Dc(C)-clsys-in-Pow(P)}
that the set $\DcClSys(\negtinysp P)$ of all \dircl\ closure systems in~$P$ is a~closure system in $\PowP$,
	and so it is also a~closure system in~$\ClSys(\negtinysp P)$.
Then Lemma~\ref{lem:dcpo-P-dirclosed-clsyss} tells~us
that the set $\DcClSys(\negtinysp P)$
	is the same as
the set $\ScClSys(\negtinysp P)$ of the closure systems
	that are associated with the \Scottcont\ closure operators on~$P$.
The~isomorphism of~complete lattices
	$\Cl(\negtinysp P) \to \ClSys(\negtinysp P)^\op : \gamma \mapsto \fix(\gamma)$
restricts to the isomorphism of~complete lattices
	$\ScCl(\negtinysp P) \to \DcClSys(\negtinysp P)^\op\tinysp$,
and so, in particular, $\ScCl(\negtinysp P)$ is an~interior system in~$\Cl(\negtinysp P)\tinysp$;
the latter we already know,
but now we gained an insight into why it~is inevitable.

Since $\DcClSys(\negtinysp P)$ is a~closure system in $\PowP\negtinysp$,
for every subset~$X$ of~$P$
	there exists the~least of all \dircl\ closure systems that include~$X$,
	which we denote by $\widedt{\dcclsys(X)}\tinysp$.
The endomap $\widedt\dcclsys$ on $\PowP$
	is a~closure operator on $\PowP\negtinysp$,
	and it restricts to a~closure operator on~$\ClSys(\negtinysp P)\tinysp$.
Since $\ScCl(\negtinysp P)$ is an interior system in~$\Cl(\negtinysp P)\tinysp$,
for every closure operator~$\gamma$ on~$P$
    there exists the~greatest of all \Scottcont\ closure operators on~$P$
	that are below $\gamma\tinysp$,
    which we denote by $\widedt{\sccore(\gamma)}$
	and call it the \notion{\Scottcont~core} of the closure operator~$\gamma\tinysp$.
The~endomap~$\widedt\sccore$ on~$\Cl(\negtinysp P)$ is an~inte\-rior~operator on~$\Cl(\negtinysp P)\tinysp$.

Via~the isomorphism $\Cl(\negtinysp P) \to \ClSys(\negtinysp P)^\op : \gamma\mapsto\fix(\gamma)$
	of complete lattices
the interior operator~$\widedt\sccore$ on~$\Cl(\negtinysp P)$
	corresponds to
the restriction of the closure operator $\widedt\dcclsys$ to $\ClSys(\negtinysp P)\dtinysp$:

\thmskip

\begin{lem}\label{lem:Fix(sccore(gamma))=dcclsys(Fix(gamma))}
If\/~$\gamma$ is a~closure operator on a~dcpo\/~$P$
	then\/~$\dtinysp\fix(\sccore(\gamma)) = \dcclsys(\fix(\gamma))\tinysp$.
\qed
\end{lem}

\thmskip

For a~general dcpo $P$ we cannot say much
    about the closure operator~$\widedt\dcclsys$
	on~$\PowP\negtinysp$
	    or its restriction to $\ClSys(\negtinysp P)$,
    or about the interior operator $\widedt\sccore$ on~$\Cl(\negtinysp P)\tinysp$.
However, if~$P$~is a~domain,
then there exist explicit constructions
	of the \Scottcont\ core $\widedt{\sccore(\gamma)}$
		of any~closure operator~$\gamma$ on~$P$,
	and of the \dircl\ closure system $\widedt{\dcclsys(C)}$
		generated by any~closure system~$C$ in~$P\tinysp$;
these two constructions are described in
	Proposition~\ref{prop:sc(gamma)-on-domain} and Proposition~\ref{prop:dcclsys-in-domain}.

\txtskip

The following definitions are from~\cite{Gierz-etal-CL&D-2003}.

Let $P$ be a~poset.

For any elements $x$ and $y$ of~$P$
we say that \notion{$x$ is way below~$y\tinysp$}, and write $x\waybelow y\tinysp$,
if~for~every joinable directed subset~$D$ of~$P$
the inequality $y\leq\tJoin\!D$ implies that $x\leq d$ for some~$d\in\negdtinysp D$.

For~every~$x\in P$ we define the set
	$\ddown x\defeq\set{\widet{u\narrt\in P\suchthat u\narrt\waybelow x}}\tinysp$,
which is a~lower set of~$P$ included in the principal ideal $\ldown x\tinysp$.

The poset $P$ is said to be \notion{continuous}
if it satisfies the \notion{axiom of approximation}:
for every $x\in P$ the set $\ddown x$ is directed
	and has in $P$ the join $\tJoin\negtinysp\ddown x = x\dtinysp$.

A~\notion{domain} is a~continuous dcpo.

\txtskip

We shall silently use the basic properties of the way-below relation.
Besides those we will also need the following two results from~\cite{Gierz-etal-CL&D-2003}.

The first result is the \notion{interpolation property}
	of the way-below relation on a~continuous poset~$P$
	(Theorem~I-1.9(ii)):
{\it for any\/ $x,\tinysp z\in P$ such that\/ $x\waybelow z$
	there exists\/~$y\in P$ so that\/ $x\waybelow y\waybelow z\tinysp$}.

The second result is a~characterization of \Scottcont\ functions between domains
	(Proposition~II-2.1(5)):
{\it a~function\/ $f\colon P\to Q\tinysp$, where\/ $P$ and\/ $Q$ are domains,
is \Scottcont\ iff\/ $f(x) = \tJoinfof(\ddown x)$ for every\/ $x\in P$}
(that is, $f(x)$~is the join of~$f(\ddown x)$ in~$Q$,
or, more long-windedly, the join of $f(\ddown x)$ in $Q$ exists and is~equal~to~$f(x)$).

\txtskip
\pagebreak[3]

Here comes the first of the promised constructions,
namely the construction of the \Scottcont\ core of a~closure operator on a~domain.

\thmskip

\begin{prop}\label{prop:sc(gamma)-on-domain}
Let\/ $\gamma$ be a~closure operator on a~domain\/~$P\negtinysp$.
Then for all\/ $x\in P$ we have\/ $\sccore(\gamma)(x)=\tJoin\negdtinysp\gamma(\ddown x)$.
\end{prop}

\interskip

\begin{myproof}
We define the endomap~$\gamma\inter$ on~$P$
	by $\gamma\inter(x)\defeq\tJoin\negdtinysp\gamma(\ddown x)$ for $x\in P\negtinysp$.

$\gamma\inter$ is ascending:
    $\gamma\inter(x)
	= \tJoin\negdtinysp\gamma(\ddown x)
	\geq \tJoin\!\ddown x
	= x\tinysp$.

$\gamma\inter$ is increasing.
If $x\leq y\tinysp$, then $\ddown x\subseteq\ddown y\tinysp$,
which clearly implies $\gamma\inter(x)\leq\gamma\inter(y)\tinysp$.

$\gamma\inter$ is idempotent.
It suffices to prove that $\gamma\inter(\gamma\inter(x))\leq\gamma\inter(x)\tinysp$,
and to prove this inequality it suffices to prove
that every element of $P$ which is way below the left hand side is below the right hand side.
So~let $u \waybelow \gamma\inter(\gamma\inter(x))
		= \tJoin\negtinysp\gamma(\ddown\gamma\inter(x))\tinysp$.
The join is directed, thus there exists $v\in\ddown\gamma\inter(x)$
	such that $u\leq\gamma(v)\tinysp$.
Now $v\waybelow\gamma\inter(x) = \tJoin\negdtinysp\gamma(\ddown x)\tinysp$,
where the join is directed,
thus $v\leq \gamma(w)$ for some $w\waybelow x\tinysp$.
It follows that
    $u \leq \gamma(v) \leq \gamma(w)
	\leq \tJoin\negdtinysp\gamma(\ddown x)
	= \gamma\inter(x)$.

We have proved that $\gamma\inter$ is a~closure operator on~$P\negtinysp$.

Evidently $\gamma\inter\leq\gamma\tinysp$.

The closure operator $\gamma\inter$ is \Scottcont.
Since $P$ is a~domain, we will prove that $\gamma\inter$~is \Scottcont\
	when we prove that $\tJoin\negdtinysp\gamma\inter(\ddown x)=\gamma\inter(x)\tinysp$
		for every $x\in P$.
It suffices to prove the inequality~$\geq\dtinysp$.
By the definition of $\gamma\inter$ we have
\begin{align}
\gamma\inter(x)
	&\Eq \tJoin_{\!u\waybelow x}\negdtinysp\gamma(u)~,
							\label{eq:gammainter(x)} \\[.25ex]
\tJoin\negdtinysp\gamma\inter(\ddown x)
	&\Eq \tJoin_{\!v\waybelow x}\negdtinysp\gamma\inter(v)
	\Eq \tJoin_{\!v\waybelow x}\!\tJoin_{\!u\waybelow v}\negtinysp\gamma(u)~.
							\label{eq:Join(gammainter(ddown(x)))}
\end{align}
Given any $u$ way below $x$,
there exists, because of the interpolation property,
an~element~$v\in P$ such that $u \waybelow v \waybelow x\tinysp$,
which shows that the term $\gamma(u)$
	in the join in~\eqref{eq:gammainter(x)}
appears also in the double join in~\eqref{eq:Join(gammainter(ddown(x)))}.
This proves the inequality
	$\tJoin\negdtinysp\gamma\inter(\ddown x) \geq \gamma\inter(x)\tinysp$.

Let $\beta$ be a~\Scottcont\ closure operator on~$P$ such that $\beta\leq\gamma\tinysp$.
Then for every $x\in P$ we have
    $\beta(x)
	= \tJoin\!\beta(\ddown x)
	\leq\tJoin\negdtinysp\gamma(\ddown x)
	= \gamma\inter(x)\tinysp$. %

We conclude that $\gamma\inter = \sccore(\gamma)$.
\end{myproof}

\thmskip

For the second construction,
that of a~\dircl\ closure system generated by a~closure system in a~domain,
we have to introduce an operation.

For each subset $X$ of a~dcpo $P$ we let $\widedt{\dj(X)}$ denote
	the set of the joins of all directed subsets of~$X$.
We have $X\subseteq\dj(X)$ because one-element sets are directed,
thus the mapping $\dtinysp\dj\colon\PowP\to\PowP$ is ascending.
It is clear that the mapping $\widedt\dj$ is increasing.
But $\widedt\dj$ is in general not idempotent;
more often than not it is very far from being a~closure operator.
The following proposition thus comes as a~slight surprise.

\thmskip

\begin{prop}\label{prop:dcclsys-in-domain}
If\/ $C$ is a~closure system in a~domain\/~$P\negtinysp$, then\/~$\dcclsys(C)=\dj(C)\dtinysp$;
therefore, if\/~$X$ is any subset of\/~$P\negtinysp$, then\/~$\dcclsys(X)=\dj(\clsys(X))\tinysp$.
\end{prop}

\interskip

\begin{myproof}
Let $P$ be a~domain, $C$ a~closure system in~$P$,
	and $\gamma$ a~closure operator on~$P$ with $\dtinysp\fix(\gamma)=C$.
Since~$\dtinysp\dcclsys(C)$ includes $C$ and is closed under directed joins,
	it~includes~$\dtinysp\dj(C)\tinysp$.
On~the other hand,
$\dtinysp\dcclsys(C)
	= \dcclsys(\fix(\gamma))
	= \fix(\sccore(\gamma))
	= \sccore(\gamma)(\negtinysp P)$
is the set of the closures
    $\dtinysp\sccore(\gamma)(x)
	= \tJoin\negdtinysp\gamma(\ddown x)$
for all $x\in P\negtinysp$.
Since for each $x\in P$ the set $\gamma(\ddown x)$ is a~directed subset of $C$,
its join belongs to $\dj(C)\tinysp$.
This proves the inclusion $\dtinysp\dcclsys(C)\subseteq\dj(C)\tinysp$.

If $X$ is any subset of $P\negtinysp$,
then applying $\widedt\dcclsys$ to $X\subseteq \clsys(X) \subseteq \dcclsys(X)$
we get $\dtinysp\dcclsys(X) = \dcclsys(\clsys(X)) = \dj(\clsys(X))\tinysp$.
\end{myproof}

\thmskip

Proposition~\ref{prop:dcclsys-in-domain} generalizes Theorem~4-1.22 in~\cite{Gratzer-Wehrung-LT:STA2-2016}.
\pagebreak[3]

\section{The frame of nuclei on a preframe}
\label{sec:frame-of-nuclei-on-preframe}

In this section we carry out the project that is only sketched in broad outline
at the end of Section~3 in~\cite{Escardo-JFN-2003}.

\txtskip

A~\notion{preframe}%
\footnote{\,Preframes are also known as meet-continuous semilattices.
See Definition O-4.1 in~\cite{Gierz-etal-CL&D-2003}.}
is a~dcpo $P$ that is also a~meet-semilattice,%
\footnote{A meet-semilattice is a~poset in which any two elements have a~meet,
	or equivalently, in which every nonempty finite subset has a~meet.}%
\,%
\footnote{\,In~\cite{Escardo-JFN-2003},
a~preframe is understood as a~dcpo that is also a~meet-semilattice \emph{with a~top element}
in~which binary meets distribute over directed joins;
quoting almost verbatim,
``a~poset in which there exist finite meets and directed joins,
	and the former distribute over the latter,
is known as a~\knownas{preframe}''.
As a~preframe is defined in the present paper,
it is not required to possess a~top element.
This additional generality seems inconsequential, but isn't.}
in which binary meets distribute over directed joins;
that is, the \notion{directed distributive law} holds:
\begin{equation*}
x\meet\tJoin\negtinysp Y \Eq \Join_{y\in Y} (x\meet y)~,
	\qquad\quad \text{$x\in P$\,, \ directed $Y\subseteq P$}\,.
\end{equation*}
Note that if $x\in P$, and $Y\subseteq P$ is directed,
	then also $\set{x\meet y\suchthat y\in Y}$ is directed.

A~\notion{frame} is a~complete lattice~$L$
in which binary meets distribute over arbitrary joins,
which means that the following \notion{infinite distributive law} holds in~$L\tinysp$:
\begin{equation*}
x\meet\tJoin\negtinysp Y \Eq \Join_{y\in Y}(x\meet y)~,
	\qquad\quad \text{$x\in L$\,, \ $Y\subseteq L$}\,.
\vspace{-1ex}
\end{equation*}

\txtskip

A~closure operator on the meet-semilattice~$P$ that preserves binary meets
	(hence preserves nonempty finite meets%
\footnote{A~``$\tinysp$nonempty finite meet''
	is short for a~``$\tinysp$meet of a~nonempty finite set''.}%
)
is called a~\notion{nucleus~on~$P\negtinysp$}.
A~preclosure map on~$P$ that preserves binary meets
is called a~\notion{prenucleus~on~$P\negtinysp$}.
A~map $\gamma\colon P\to P$ is a~prenucleus
iff~it is ascending and preserves binary meets,
and it is a~nucleus
iff~it is ascending and idempotent and preserves binary meets;
in both cases $\gamma$ is increasing because it preserves binary meets.
If~$\gamma$ is a~prenucleus, and~$P$ possesses a top element~$\top$,
then~$\gamma(\top) = \top$ because $\gamma$ is ascending;
that is, if~the~empty meet%
\footnote{\,When we say ``$\tinysp$the empty meet'' we mean ``$\tinysp$the meet of the empty set''.}
in $P$ exists,
then $\gamma$ preserves~it.

We let $\Prenuc(\negtinysp P)$ and $\Nuc(\negtinysp P)$ denote the pointwise-ordered sets
	of all prenuclei resp.\ of all~nuclei on a~meet-semilattice~$P\negtinysp$.
The fixpoint set of a~nucleus on~$P$ shall be called a~\notion{nuclear system in~$P\negtinysp$},
and the poset of all nuclear systems in~$P\negtinysp$, ordered by inclusion,
shall be denoted by $\NucSys(\negtinysp P)\tinysp$.

\txtskip

We shall prove, among other things, that for any preframe $P$ the poset $\Nuc(\negtinysp P)$ is~a~frame.

\txtskip

For a~while, let $P$ be any~meet-semilattice.

\txtskip

In~the poset $\Fun(\negtinysp P)$ of all endofunctions on~$P$ ordered pointwise,
any~two endofunctions $\gamma$ and $\delta$ have a~meet $\gamma\meet\delta$,
	which is calculated pointwise,
and the following~is~true:
\begin{itemize}[topsep=.5ex,itemsep=.5ex,leftmargin=4.1em]
\item[(i)\:] If $\gamma$ and $\delta$ are ascending, so is $\gamma\meet\delta\tinysp$.
\item[(ii)\:] If $\gamma$ and $\delta$ are increasing, so is $\gamma\meet\delta\tinysp$.
\item[(iii)\:] If $\gamma$ and $\delta$ are closure operators, so is $\gamma\meet\delta\tinysp$.
\item[(iv)\:] If $\gamma$ and $\delta$ preserve binary meets, so does $\gamma\meet\delta\tinysp$.
\end{itemize}
Properties (i) and (ii) are easily verified.
For~(iii), assume $\beta$ and $\gamma$ are closure operators~on~$P$
and put $\beta = \gamma\meet\delta$.
Then $\beta$ is increasing and ascending by~(i)~and~(ii).
Since $\beta\beta \leq \gamma\tinysp\gamma = \gamma$,
and similarly $\beta\beta \leq \delta$,
we have $\beta\beta \leq \gamma\meet\delta = \beta$,
thus $\beta$~is idempotent.
Finally, to~prove~(iv), assume that $\gamma$ and $\delta$ preserve binary meets
	and let~$x,\,y\in P\dtinysp$;
then
\begin{align*}
(\gamma\meet\delta)(x\meet y)
  & \Eq \gamma(x\meet y)\meet\delta(x\meet y)
  	\Eq \gamma(x)\meet\gamma(y)\meet\delta(x)\meet\delta(y)\\ 
  & \Eq (\gamma\meet\delta)(x)\meet(\gamma\meet\delta)(y)\,.
\end{align*}
Because of (i), (ii), and (iii), in the posets $\Precl(\negtinysp P)$
	and $\Cl(\negtinysp P)$ all binary meets exists
and they are calculated pointwise.


From (i)--(iv) above it follows
that the pointwise meet of two prenuclei on the meet-semilattice~$P$ is a~prenucleus on~$P\negtinysp$,
and that the pointwise meet of two nuclei on~$P$ is a~nucleus on~$P\negtinysp$.
Therefore, in the poset~$\Prenuc(\negtinysp P)$
	all binary meets exist and they are calculated pointwise,
and the same is true for the poset $\Nuc(\negtinysp P)$.

\txtskip

From here on let $P$ be a~preframe.

\txtskip

How about the joins of sets of nuclei on the~preframe~$P\dtinysp$?
They~always exist,
and~they are taken in the complete lattice $\Cl(\negtinysp P)$ of all closure operators on~$P\negtinysp$.
We~give a slightly more general~result.

\thmskip

\begin{thm}\label{thm:in-preframe-prenuclei-generate-nucleus}
Let\/ $P$~be a~preframe,
and let\/ $\Gamma$ be a~set of prenuclei on\/~$P\tinysp$;
then the closure operator\/ $\gencl{\Gamma}$,
	generated by the set\/ $\Gamma$ of~preclosure maps,
is a~nucleus.
In~particular, if\/~$\Gamma$ is a~set of nuclei on\/~$P\negtinysp$,
then the join\/~$\tJoin\negtinysp\Gamma$, taken in the complete lattice\/~$\Cl(\negtinysp P)$,
is~a~nucleus.
\end{thm}

\interskip

\begin{myproof}
Put $\delta = \gencl{\Gamma}$.
Let $x,\,y\in P\negtinysp$.
Since~$\delta$ is increasing, we have $\delta(x\meet y) \leq \delta(x)\meet\delta(y)$;
we~must prove that the converse inequality also holds.

We~will first prove the weaker assertion $x\meet \delta(y) \leq \delta(x\meet y)$.%
\footnote{\,In the terminology of~\cite{Simmons-AF-2006},
$\delta$~is a~stable inflator, which turns out to be a~nucleus because it is a~closure operation.
The~discussion in~\cite{Simmons-AF-2006} is restricted to frames,
but closure operations, nuclei, and stable inflators can be defined,
	and they are related as just mentioned,
in an arbitrary meet-semilattice.}
Let~$A$ be the set of all $z\in P$ such that $x\meet z \leq \delta(x\meet y)$.
The~set~$A$ contains~$y\tinysp$, and it is \dircl\ by directed distributivity.
For any $z\in A$ and any $\gamma\in\Gamma$ we have
\begin{equation*}
x\meet\gamma(z)
  \Leq \gamma(x)\meet\gamma(z)
  \Eq \gamma(x\meet z)
  \Leq \gamma\bigl(\delta(x\meet y)\bigr)
  \Eq \delta(x\meet y)\,,
\end{equation*}
hence $\gamma(z)\in A\dtinysp$;
thus $A$ is closed under~$\Gamma$.
The induction principle gives~$\delta(y)\in\nolinebreak A\tinysp$.

Now we substitute $\delta(x)$ for~$x$ in $x\meet \delta(y) \leq \delta(x\meet y)$ and~get
\begin{equation*}
\delta(x)\meet\delta(y)
  \Leq \delta\bigl(\delta(x)\meet y\bigr)
  \Leq \delta\bigl(\delta(x\meet y)\bigr)
  \Eq \delta(x\meet y)\,,
\end{equation*}
where the second inequality holds because $\delta(x)\meet y \leq \delta(x\meet y)$
	and $\delta$ is increasing.
\end{myproof}

\interthmskip

\begin{cor}\label{cor:complete-lattice-of-nuclei}
For any preframe\/~$P\negtinysp$, the subset\/ $\Nuc(\negtinysp P)$ of\/ $\Cl(\negtinysp P)$
is closed under arbitrary joins in the complete lattice\/~$\Cl(\negtinysp P)$,
so it is a~complete lattice.%
\footnote{\,Here we see that not requiring that preframes have top elements
	is not so innocent as it seems.
For~a~preframe~$P$ with a~top element
	it is trivial that there is a~greatest nucleus~on~$P\negtinysp$,
namely the constant map sending every element of~$P$ to its top element.
In~contrast, the existence of a~greatest nucleus
	on a~preframe~$P$ which lacks a~top element
is a~nontrivial matter;
and the greatest nucleus on~$P$,
	the top element of the complete lattice~$\Nuc(\negtinysp P)$,
	does exist.}
Also,~$\Nuc(\negtinysp P)$~is closed under binary meets
	{\rm(}and hence under nonempty finite meets\/{\rm)}
		in\/~$\Cl(\negtinysp P)$,
since binary meets in\/~$\Nuc(\negtinysp P)$ are calculated pointwise,
	same as they are calculated in\/~$\Cl(\negtinysp P)$.
\qed
\end{cor}

\thmskip

If~$P$ possesses a~top element $\top\negtinysp$,
then the top nucleus is the same as the top closure operator,
which is the constant map $P\to P : x\mapsto\top\negtinysp$.
However, if $P$ does not have a~top element,
then the top nucleus might be strictly smaller than the top closure operator.
That~is, though a~subset $\Nuc(\negtinysp P)$ of $\Cl(\negtinysp P)$ is closed under all joins
	and also under all nonempty finite meets,
		both taken in the complete lattice $\Cl(\negtinysp P)\tinysp$,
it might not be closed under the empty meet taken in~$\Cl(\negtinysp P)\tinysp$.
This can already happen in a~finite preframe.
Since every finite directed set has a~greatest element, which is its join,
every finite poset is a~dcpo and every finite meet-semilattice is a~preframe.
The meet-semilattice $P_1$ in the left panel of Figure~\ref{fig:P1-preframe-with-topnuc<topcl}
\begin{figure}[!htp]\centering
\vspace{.5ex}
\includegraphics[draft=false]{mp/clops-on-dcpos-1.mps}
\vspace{-.5ex}
\mycaption{\label{fig:P1-preframe-with-topnuc<topcl}%
	The only nuclear system in the preframe~$P_1$ is the whole~$P_1$
	(represented by the black dot in the right panel),
	so the only nucleus on~$P_1$ is the identity map,
	which is different from the top closure operator on~$P_1$.}
\vspace{-1ex}
\end{figure}
is the simplest possible example:
there are four closure operators corresponding
	to the four closure systems exhibited in the right panel,
while there is only one nucleus, namely the identity map,
	which is the bottom closure operator.

\txtskip

Since $\Nuc(\negtinysp P)$ is an interior system in the complete lattice $\Cl(\negtinysp P)$,
for every closure operator $\gamma$ on the preframe $P$
there exists the largest nucleus below $\gamma\tinysp$,
the \notion{nuclear core} $\widedt{\nuccore(\gamma)}$ of the closure operator~$\gamma\tinysp$.

\thmskip

\begin{thm}\label{thm:frame-of-nuclei-on-preframe}
For any preframe\/~$P$ the complete lattice\/~$\Nuc(\negtinysp P)$ is a~frame.
\end{thm}

\interskip

\begin{myproof}
Let $\beta\in\Nuc(\negtinysp P)$ and $\Gamma\subseteq\Nuc(\negtinysp P)$,
and write $\delta \defeq \tJoin\negtinysp\Gamma$,
	$\deltapr \defeq \tJoin_{\leftdt{\gamma\in\Gamma}}(\beta\meet\gamma)$
(where meets and joins are taken in~$\Nuc(\negtinysp P)$, hence in~$\Cl(\negtinysp P)$).
We~must show that $\beta\meet\delta = \deltapr\negdtinysp$.
The~inequal\-ity~$\beta\meet\delta \geq \deltapr $
holds because $\beta\meet\delta \geq\nolinebreak \beta\meet\nolinebreak\gamma$
for every $\gamma\in\Gamma$.

To~prove the converse inequality,
let~$x\in P$ and put
	$A \defeq \bigset{z\in P\bigsuchthat \beta(x)\meet z \leq \deltapr (x)}\tinysp$.
Evidently $x\in A\tinysp$, and $A$ is \dircl\ by directed distributivity.
In~order to see that $A$~is closed under~$\Gamma$,
consider any $\gamma\in\Gamma$ and any~$z\in A\tinysp$,
exhibit the following chain of equalities and inequalities,
\begin{equation}\label{eq:beta(x)meetgamma(z)...<=...delta'(x)}
\begin{aligned}
\beta(x) \meet \gamma(z)
	&\Eq \beta(x) \widet\meet \gamma\beta(x) \widet\meet \gamma(z)
		\Eq \beta(x) \widet\meet \gamma\bigl(\beta(x)\narrt\meet z\bigr)
		\Leq \beta(x) \widet\meet \gamma\bigl(\deltapr(x)\bigr) \\ 
	&\Leq \beta\bigl(\deltapr(x)\bigr) \widet\meet \gamma\bigl(\deltapr(x)\bigr)
		\Eq (\beta\narrt\meet\gamma)\bigl(\deltapr(x)\bigr) \\ 
	&\Eq \deltapr(x)~,
\end{aligned}
\end{equation}
and then from the inequality between the first and the last expression in the chain
	conclude that $\gamma(z)\in A\tinysp$.
By the induction principle it then follows that $\delta(x)\in A\tinysp$,
that~is, that $(\beta\meet\delta)(x) = \beta(x)\meet\delta(x) \leq \deltapr (x)\tinysp$.
\end{myproof}

\thmskip
\pagebreak[3]

In conclusion to this nuclear-themed section we take a~look at \Scottcont\ nuclei on a~preframe.

\txtskip

Let $P\negtinysp$ be a~preframe.

We denote by $\ScNuc(\negtinysp P)$
	the poset of all \Scottcont\ nuclei on~$P\negtinysp$, ordered pointwise,
and by $\ScNucSys(\negtinysp P)$
	the poset of the fixpoint sets of the \Scottcont\ nuclei, ordered by inclusion.

By Lemma~\ref{lem:dcpo-P-dirclosed-clsyss}
a~subset of~$P$ is the fixpoint set of a~\Scottcont\ nucleus iff it~is a~\dircl\ nuclear system,
which means that $\ScNucSys(\negtinysp P) = \DcNucSys(\negtinysp P)$.
Since $\NucSys(\negtinysp P)$ is a~closure system in $\PowP$,
so is then $\DcNucSys(\negtinysp P)$, by Lemma~\ref{lem:P-poset&C-clsys-in-Pow(P)=>Dc(C)-clsys-in-Pow(P)};
it~follows that $\DcNucSys(\negtinysp P)$~is also a~closure system in~$\NucSys(\negtinysp P)$
	as~well~as in~$\ClSys(\negtinysp P)$ and in~$\DcClSys(\negtinysp P)$.
Passing these observations
	via the antiisomorphism between closure systems and closure operators
we see, at a~glance, that $\ScNuc(\negtinysp P)$~is an~interior system in~$\Nuc(\negtinysp P)$
	as~well~as in~$\Cl(\negtinysp P)$ and in~$\ScCl(\negtinysp P)$.
The following theorem almost does not need a proof, though we do provide one. 

\thmskip

\begin{thm}\label{thm:in-preframe-nuc-generd-by-prenucs}
Let\/ $\Gamma$ be a~set of \Scottcont\ prenuclei on a~preframe\/~$P\negtinysp$.
The pointwise join\/~$\delta\defeq\dottJoin\negtinysp\Gamma^*$
	is a~\Scottcont\ nucleus on\/~$P\negtinysp$.
Moreover, $\delta$ is the closure operator on the dcpo\/~$P$
	that is generated by the set\/~$\Gamma$ of preclosure maps on\/~$P\negtinysp$.

If\/ $\Gamma$ is a~set of \Scottcont\ nuclei on\/~$P\negtinysp$,
then the pointwise join of\/~$\Gamma^*$
is the join of\/~$\Gamma$ in\/~$\ScNuc(\negtinysp P)$
	as well as in\/~$\Nuc(\negtinysp P)$ and in\/~$\Cl(\negtinysp P)\tinysp$.
The~set\/~$\tinysp\ScNuc(\negtinysp P)$ is
	an~interior system in the complete lattice\/~$\Nuc(\negtinysp P)$
		and also in the complete lattice\/~$\Cl(\negtinysp P)$,
and so as~a~poset it~is  a~complete lattice.
\end{thm}

\interskip

\begin{myproof}
According to Theorem~\ref{thm:in-dcpo-Scottcont-preclmaps-gener-Scottcont-clop},
	$\delta$ is a~\Scottcont\ closure operator on~$P\negtinysp$,
	and it is also the closure operator on~$P$ generated by the set $\Gamma$ of preclosure maps.
By Theorem~\ref{thm:in-preframe-prenuclei-generate-nucleus}
	it then follows that $\delta$ is a~nucleus.
\end{myproof}

\interthmskip

\begin{prop}\label{prop:meet-of-two-Scotcont-nuclei-is-Scottcont}
The pointwise meet of two \Scottcont\ nuclei on a~preframe\/~$P$
is a~\Scottcont\ nucleus on\/~$P$.
\end{prop}

\interskip

\begin{myproof}
Let $\gamma$ and $\delta$ be \Scottcont\ nuclei on~$P$.
We know that the pointwise meet $\gamma\meet\delta$ is a~nucleus,
so it remains to prove that the meet is \Scottcont.%
\footnote{The reasoning that follows, suitably adapted, can be made to prove
	 that the pointwise meet of two \Scottcont\ increasing maps from a~dcpo to a~preframe
		is \Scottcont\ (and, of~course, increasing).}

Let $Y$ be a~directed subset of~$P\tinysp$;
we shall prove that
	$(\gamma\narrt\meet\delta)(\tJoin\negtinysp Y) = \tJoin(\gamma\narrt\meet\delta)(Y)\tinysp$.
It~suffices to prove the inequality $\leq\dtinysp$.
We calculate:
\begin{equation}\label{eq:(gamma-meet-delta)(Join(Y))=...}
\begin{aligned}
(\gamma\narrt\meet\delta)(\tJoin\negtinysp Y)
    &\Eq \gamma(\tJoin\negtinysp Y) \meet \delta(\tJoin\negtinysp Y)
	\Eq \tJoin\negdtinysp\gamma(Y) \widedt\meet \tJoin\negdtinysp\delta(Y) \\ 
    &\Eq \tJoin_{\!y_1,y_2\in Y} \bigl(\gamma(y_1)\meet\delta(y_2)\bigr)~;
\end{aligned}
\end{equation}
the last equality holds because of directed distributivity in the preframe~$P$.
Now if $y_1,\tinysp y_2\in Y$,
then there exists $y\in Y$ such that $\set{y_1,y_2}\leq y\tinysp$,
and then
    $\gamma(y_1)\meet\delta(y_2)
	\leq \gamma(y) \meet \delta(y)
	= (\gamma\narrt\meet\delta)(y)
	\leq \tJoin(\gamma\narrt\meet\delta)(Y)\tinysp$.
The last join in~\eqref{eq:(gamma-meet-delta)(Join(Y))=...}
	is therefore $\leq\tJoin(\gamma\narrt\meet\delta)(Y)$,
and we have the desired inequality.
\end{myproof}

\thmskip
\pagebreak[3]

The following is a~straightforward consequence
	of~Theorem~\ref{thm:in-preframe-nuc-generd-by-prenucs}
	and Proposition~\ref{prop:meet-of-two-Scotcont-nuclei-is-Scottcont}.
\thmskip
\begin{cor}\label{cor:ScNuc(preframe)-is-frame}
Let\/ $P$ be a~preframe.
The subposet\/ $\ScNuc(\negtinysp P)$ of the frame\/ $\Nuc(\negtinysp P)$
	is in $\Nuc(\negtinysp P)$ closed under all joins and under binary meets,
therefore it is itself a~frame
	whose arbitrary joins and binary meets are inherited from\/~$\Nuc(\negtinysp P)$.
\qed
\end{cor}

\thmskip

It is an open question, for the author, whether there exists a~preframe~$P$
such that the top element of~$\ScNuc(\negtinysp P)$
	is different from the top element of~$\Nuc(\negtinysp P)$.
If such a~preframe exists, it lacks a~top element,
and it is by necessity infinite
	since every increasing function between finite posets is \Scottcont\
and so, in particular,
	every nucleus on a~finite preframe (that is, on a~finite meet-semilattice) is~\Scottcont.

\section{Nuclei on frames}
\label{sec:nucs-on-frames}

In this section we consider frames, as preframes with special properties.
For starters we specialize Theorem~\ref{thm:frame-of-nuclei-on-preframe} to frames.

\thmskip

\begin{cor}\label{cor:frame-of-nuclei-on-frame}
For any frame\/~$L$ the complete lattice\/~$\Nuc(L)$ is a~frame. \qed
\end{cor}

\thmskip

The concise proof of Theorem~\ref{thm:frame-of-nuclei-on-preframe},
	which also serves (specialized) as a~proof of Corollary~\ref{cor:frame-of-nuclei-on-frame},
compares well with the long-winded proof of
	Proposition~II.2.5 in~\cite{Johnstone-SS-1982}.

\txtskip

From now on, to the end of the section, we let $L$ be a~frame.

\txtskip

Arbitrary joins and finite meets in the complete lattice $\Nuc(L)$
are calculated in the complete lattice $\Cl(L)$, as in any preframe;
but, $L$ being a~frame, more~is~true.

\thmskip

\begin{prop}\label{prop:frame-L-Nuc(L)-sub-(completlatt)-of-Cl(L)}
Let\/ $L$ be a~frame.
Arbitrary meets in the complete lattice\/~$\Nuc(L)$ are calculated pointwise,
same as they are calculated in the complete lattice\/~$\Cl(L)$.
As~arbitrary joins in\/~$\Nuc(L)$ also are calculated in\/~$\Cl(L)$,
the complete lattice\/~$\Nuc(L)$
	is a~sub$\tinysp$-{\rm(}complete lattice\/{\rm)} of the complete lattice\/~$\Cl(L)$.
\end{prop}

\interskip

\begin{myproof}
We have to prove that for any set $\Gamma$ of nuclei
its pointwise meet $\alpha\defeq\nolinebreak\tMeet\negdtinysp\Gamma$,
	which is a~closure operator, 
preserves binary meets;
but this is a~straightforward consequence of the associativity-cum-commutativity of arbitrary meets
	in the complete lattice~$L\tinysp$.
\end{myproof}

\thmskip

We have mentioned, in Section~\ref{sec:overview}, that a~frame is a~complete \Heytalg:
for any two elements $a,\,b\in L$ there exists
	the \notion{relative pseudo-complement} of~$a$ with respect to~$b\tinysp$,
which is the unique element $(a\bimpl b)\in L$ such that
\begin{equation*}
x\meet a \Leq b \wide\Isequiv x\leq (a\bimpl b)
	\qquad\quad \text{for every \ $x\in L$}~;
\end{equation*}
the infinite distributivity implies that
	$(a\bimpl b) \widedt= \tJoin\set{x\suchthat x\meet a\leq b}$.

\txtskip

The following proposition characterizes nuclear systems in a~frame.%
\footnote{A nuclear system in a~frame (alias locale) is also known as a~\knownas{sublocale}.}

\thmskip

\begin{prop}\label{prop:charact-of-nucsyss-in-frames}
Let\/ $L$ be a~frame.%
\footnote{Actually, the proposition holds for any \Heytalg~$L$.
The proof, precisely as it is, remains valid if~$L$~is just a~\Heytalg,
	because it~nowhere requires that $L$~be complete.}
A~subset\/~$C$ of\/~$L$ is a~nuclear system in\/~$L$
\iff\/ $C$ is a closure system in\/~$L\tinysp$,
and\/ $x\in L$, $y\in C$ together always imply\/ $(x\bimpl y)\in C$.
\end{prop}

\interskip

\begin{myproof}
Suppose $C$ is a~nuclear system.
Then $C$ is a~closure system.
Let $\gamma$ be the nucleus with $\dtinysp\fix(\gamma)=C$.
Let $x\in L$ and $y\in C$;
we have to prove that $x\bimpl y$ is a fixed point of~$\gamma\tinysp$.
It~suffices to prove that $\gamma(x\narrt\bimpl y)\leq(x\bimpl y)\tinysp$.
Since $x\meet(x\bimpl y) \leq y$ and $\gamma(y)=y$ we have
\begin{equation*}
x\meet\gamma(x\narrt\bimpl y)
	\widedt\Leq \gamma(x)\meet\gamma(x\narrt\bimpl y)
	\widedt\Eq \gamma\bigl(x\meet(x\narrt\bimpl y)\bigr)
	\widedt\Leq \gamma(y)
	\widedt\Eq y~; 
\end{equation*}
\pagebreak[3]
by the defining property of $x\bimpl y$ it then follows that
	$\gamma(x\narrt\bimpl y)\leq(x\bimpl y)$,
as required.

Conversely, suppose that $C$ is a~closure system,
and that $x\in L$, $y\in C$ always imply $(x\bimpl y)\in C$.
Let $\gamma$ be the closure operator with $\dtinysp\fix(\gamma)=\gamma(L)=C\tinysp$;
we shall prove that $\gamma$~preserves binary meets.
We need only prove that $x\meet\gamma(y)\leq\gamma(x\narrt\meet y)$
	for~all~$x,\tinysp y\in L\tinysp$.%
\footnote{Just as in the proof of Theorem~\ref{thm:in-preframe-prenuclei-generate-nucleus}.}
Let~$z\defeq\bigl(x\bimpl\gamma(x\narrt\meet y)\bigr) \in C$.
From $x\meet y\leq\gamma(x\narrt\meet y)$ 
	we get $y\leq \bigl(x\bimpl\gamma(x\narrt\meet y)\bigr) = z\dtinysp$;
but then $\gamma(y)\leq\gamma(z)=z$,
whence $x\meet\gamma(y)\leq x\meet z \leq\gamma(x\narrt\meet y)\dtinysp$.
\end{myproof}

\thmskip

For any two subsets $A$ and $B$ of a~frame $L$ we shall write
\begin{equation*}
(A\bimpl B) \Defeq \bigset{a\narrt\bimpl b\bigsuchthat a\narrt\in A\tinysp,\, b\narrt\in B}~.
\end{equation*}
Using this notation, Proposition~\ref{prop:charact-of-nucsyss-in-frames} says
    that a~closure system $C$ in the frame~$L$ is a~nuclear system iff $(L\bimpl C)\subseteq C$
	(which implies $(L\bimpl C)=C\tinysp$).

Every subset $X$ of a~complete lattice
generates the smallest closure system $\widedt{\clsys(X)}$ that includes~$X$;
the closure system $\widedt{\clsys(X)}$ consists of the meets of all subsets of~$X$.
Likewise every subset $X$ of the~frame~$L$
generates the smallest nuclear system $\widedt{\nucsys(X)}$ in~$L$ that includes~$X$;
can we somehow construct $\widedt{\nucsys(X)}$?

\thmskip

\begin{prop}\label{prop:in-frame-nucsys(X)=clsys(L=>X)}
If\/ $X$ is any subset of a~frame\/ $L\tinysp$, then
\begin{equation*}
\nucsys(X) \Eq \clsys\dtinysp(L\rightt\bimpl X)\tinysp~.
\end{equation*}
\end{prop}

\negdisplayhalfskip
\interskip

\begin{myproof}
The endofunction $(L\bimpl\anon)\colon X\mapsto(L\narrt\bimpl X)$ on $\PowL$
	is a~closure operator on the powerset lattice~$\PowL\dtinysp$:
it is ascending because $(\top\bimpl x) = x$ for every~$x\in X$;
it is evidently increasing;
and it is idempotent because $(a\bimpl(b\bimpl x)) = ((a\narrt\meet b)\bimpl x)$
	for all $a,\tinysp b\in L$ and every $x\in X$.
In~view of Proposition~\ref{prop:charact-of-nucsyss-in-frames}
the closure operator $\widedt\nucsys$ on $\PowL$
is the join, in the complete lattice $\Cl(\PowL)\tinysp$,
of the closure operators $\widedt\clsys$ and $(L\bimpl\anon)\tinysp$.

Let $X$ be any subset of $L\tinysp$.
We shall prove the inclusion
\begin{equation*}
\bigl(L\bimpl\clsys(X)\bigr) \Subseteq \clsys\dtinysp(L\rightt\bimpl X)~,
\end{equation*}
from which it will follow that
    $\dtinysp\clsys\dtinysp(L\bimpl\anon)
	\widet= \clsys\widet\join\tinysp(L\narrt\bimpl\anon)
	\widet= \nucsys\,$.
Consider a~general element of $\widedt{\clsys(X)}$,
    which is of the form $\tMeet_{\tinysp i\in I}\negtinysp x_i$
	for some elements $x_i$, $i\in I$, of the set~$X$;
also let $y\in L\tinysp$.
Then
\begin{equation*}
\bigl(y\dtinysp\bimpl\tMeet_{\tinysp i\in I}\negtinysp x_i\bigr)
	\widedt\Eq \tMeet_{\tinysp i\in I} (y\bimpl x_i)
	\widedt\In \clsys(L\rightt\bimpl X)~,
\end{equation*}
which proves the asserted inclusion.
\end{myproof}

\thmskip

Mark that we do not obtain a~shortcut when constructing $\widedt{\nucsys(C)}$
	for a~closure system~$C$
since we still have to construct $\widedt{\nucsys(C)}$ as $\widedt{\clsys\dtinysp(L\bimpl C)}$,
	which resists simplification.
But, if~$X$~is any subset of~$L\tinysp$,
then $\dtinysp\clsys\dtinysp(L\bimpl\clsys(X))$
	simplifies to $\widedt{\clsys\dtinysp(L\bimpl X)}$.

\txtskip
\pagebreak[3]

There is an important special case where $\dtinysp\nucsys(X)$ does simplify.

\thmskip

\begin{prop}\label{prop:in-frame-nucsys(x)=(L=>x)}
Let $L$ be a~frame.%
\footnote{\,This is another result that is true for any \Heytalg.
The proof of Proposition~\ref{prop:in-frame-nucsys(x)=(L=>x)}
	relies on the completeness of $L$, but we can do without it.
Here is how.
Let $L$ be a~\Heytalg.
Fix $a\in L$.
Since $x\leq (y\bimpl a)$ iff $x\meet y\leq a$ iff $y\leq(x\bimpl a)$,
we see that the endofunction $\relnot_a\defeq(\anon\bimpl a)$ is a~Galois complementation,
that is, that $\pair{{\relnot_a},{\relnot_a}}$ is a~Galois connection $L\galoisconn L^\op$.
Therefore $\relnot_a$ is decreasing and $\relnot_a\relnot_a\relnot_a=\relnot_a$,
which implies that $\relnot_a\relnot_a$ is a~closure operator
	with the fixpoint set $\relnot_a L = (L\bimpl a)$,
whence $\relnot_a\relnot_a$ is a~nucleus
	because $(L\bimpl(L\bimpl a)) = (L\bimpl a)$.
Also, $L\bimpl a$ is clearly the least nuclear system containing~$a$,
that is, it is $\dtinysp\nucsys(\set{a})$.}
If\/ $x$ is any element of\/~$L\tinysp$, then
\begin{equation*}
\nucsys(\set{x}) \Eq (L\bimpl x)~.
\end{equation*}
\end{prop}

\negdisplayhalfskip
\interskip

\begin{myproof}
It suffices to prove that $L\bimpl x$ is closed under all meets,
so that $\dtinysp\nucsys(\set{x}) = \clsys\dtinysp(L\bimpl x) = (L\bimpl x)\tinysp$.
And indeed, if $y_i$, $i\in I$, are any elements of~$L\tinysp$,
then $\tMeet_{\tinysp i\in I}(y_i\bimpl x)
	= ((\tJoin_{\!i\in I} y_i)\bimpl x)
	\in (L\bimpl x)\tinysp$.
\end{myproof}

\thmskip

We are not yet satisfied.
Now that we have constructed $\widedt{\nucsys(X)}$,
the nuclear system generated by a~subset $X$ of a~frame,
we are curious what the corresponding nucleus $\dtinysp\nuc_X$ looks like,
	the one whose fixpoint set is $\widedt{\nucsys(X)}$.

\thmskip

\begin{prop}\label{prop:in-frame-nucX(y)=Meet(x-in-X)((y=>x)=>x)}
If\/ $X$ is any subset of a frame\/ $L\tinysp$, then
\begin{equation*}
\nuc_X(y) \Eq\! \Meet_{x\in X}\!\bigl((y\bimpl x)\bimpl x\bigr)
			\qquad\quad \text{\rm for\, $y\in L$}\,.
\end{equation*}
\end{prop}

\interskip

\begin{myproof}
For every $y\in L$ we have
\begin{equation*}
\nuc_X(y) \Eq \Meet\negtinysp\bigl(\allabove{\nucsys(X)}{y}\bigr)
	\Eq \Meet\tinysp\bigl(\allabove{\clsys(L\narrt\bimpl\negdtinysp X)}{y}\bigr)
	\Eq \Meet\tinysp\bigl(\allabove{(L\narrt\bimpl\negdtinysp X)}{y}\bigr)~.
\end{equation*}
Fix $x\in X$, and let $u\in L\tinysp$.
Then $y\leq(u\bimpl x)$ iff $u\meet y\leq x$ iff $u\leq (y\bimpl x)\tinysp$,
and for every $u\leq(y\bimpl x)$ we have $(u\bimpl x)\geq((y\bimpl x)\bimpl x)\tinysp$.
The meet of all terms of the form $u\bimpl x$ in
	$\tMeet(\tinysp\allabove{(L\narrt\bimpl\negdtinysp X)}{y})$
is $(y\bimpl x)\bimpl x\tinysp$.
Now we release~$x$ to run through the whole set~$X$
and obtain the formula for~$\widedt{\nuc_X(y)}$ given in the proposition.
\end{myproof}

\thmskip

An important special case of
	Proposition~\ref{prop:in-frame-nucX(y)=Meet(x-in-X)((y=>x)=>x)}
has $X$ consisting of a~single element.

\thmskip

\begin{cor}\label{cor:in-frame-nucx(y)=((y=>x)=>x)}
If\/ $x$ is any element of a~frame\/ $L$, then
\begin{equation*}
\nuc_{\set{x}}(y) \Eq \bigl((y\bimpl x)\bimpl x\bigr)
			\qquad\quad \text{\rm for\, $y\in L$}\,. \tag*{\qed}
\end{equation*}
\end{cor}

\interskip

Let us introduce, for each $x\in L$,
	the nucleus $\regnuc{x} \defeq \nuc_{\set{x}}$~on~$L\dtinysp$;
we shall leave it nameless.%
\footnote{\,In~\cite{Wilson-AT&CAAFT-1994} the nucleus $\regnuc{x}$
	is called \knownas{quasi-closed} and is written~$q(x)$.
In~\cite{Johnstone-SS-1982} only the (nameless) special case
	$\regnuc{0}=\relnot\relnot$ makes a~cameo appearance on page~51.}
By~Proposition~\ref{prop:in-frame-nucsys(x)=(L=>x)}
	the fixpoint set of the nucleus $\regnuc{x}$ is $L\bimpl x\dtinysp$.
The~map\-ping $L \to\nolinebreak \NucSys(L) : x \mapsto \fix(\regnuc{x})$ is injective
	since $x$~is the least element of $\dtinysp\fix(\regnuc{x}) =\nolinebreak (L\bimpl x)$,
and therefore also the mapping $L \to \Nuc(L) : x\mapsto \regnuc{x}$ is~injective.

\txtskip

Here is another special case of proposition~\ref{prop:in-frame-nucX(y)=Meet(x-in-X)((y=>x)=>x)},
with $X$ a~closure system.

\thmskip

\begin{cor}\label{cor:in-frame-nucC(y)=Meet(x-in-X)((y=>gamma(x))=>gamma(x))}
If\/ $\gamma$ is a~closure operator on a~frame\/ $L\tinysp$, then
\begin{equation*}
\nuc_{\tinysp\gamma(L)}(y)
	\Eq\!\! \Meet_{x\in\gamma(L)}\!\! \regnuc{x}(y)
	\Eq\! \Meet_{u\in L} \regnuc{\gamma(u)}(y)~.	\tag*{\qed}
\end{equation*}
\end{cor}

\negdisplayskip
\thmskip

The nucleus
	$\dtinysp\nuccore(\gamma)
	\defeq \nuc_{\tinysp\gamma(L)}
	= \nuc_{\tinysp\fix(\gamma)}\tinysp$
is the nuclear core of the closure operator~$\gamma$,
that~is,
it~is~the~greatest of all nuclei on the frame~$L$ that are below~$\gamma\tinysp$.

\thmskip
\pagebreak[3]

\begin{prop}\label{prop:nu<=rx-iff-x-in-Fix(nu)}
If\/ $\nu$ is a~nucleus on a~frame\/~$L\tinysp$ and\/ $x\in L\tinysp$,
	then\/ $\nu\leq\regnuc{x}$ iff\/ $x\in\fix(\nu)$.
In particular, if\/ $x,\tinysp y\in L\tinysp$,
	then\/ $\regnuc{x} \leq \regnuc{y}$ iff\/ $y\in(L\bimpl x)$.
\end{prop}

\interskip

\begin{myproof}
Let $\nu$ be a~nucleus on a~frame $L$ and $x\in L$.
Then $\nu\leq\regnuc{x}$
	iff $\fix(\regnuc{x})\subseteq\fix(\nu)$,
	iff $(L\bimpl x)\subseteq\fix(\nu)$,
	iff $x\in\fix(\nu)\dtinysp$;
the last equivalence holds since $x=(\top\bimpl x)$~is in~$(L\bimpl x)$,
and because $x\in\fix(\nu)$ implies $(L\bimpl x)\subseteq\fix(\nu)$.
\end{myproof}

\thmskip

According to the second assertion of the proposition,
the set~$L$,
	equipped with the relation $\regnuc{\leq}$
		defined by $x\regnuc{\leq}y$ iff $\regnuc{x}\leq\regnuc{y}$
			{\large(}iff $y\in(L\bimpl x)${\large)},
is a~poset isomorphic to the subposet~$\regnuc{L}$ of~$\Nuc(L)$,
via the isomorphism $L\to\regnuc{L} : x\mapsto\regnuc{x}\tinysp$.%
\footnote{\,The partial order $\regnuc{\leq}$ on~$L$
is in~\cite{Wilson-AT&CAAFT-1994} written $\trianglelefteq$ and called the~\knownas{regularity ordering}.}
Note that $x\regnuc{\leq}y$ implies $x\leq y\tinysp$, for all $x,\tinysp y\in L\tinysp$.

\txtskip

The following is in effect a~rephrasing
	of Proposition~\ref{prop:in-frame-nucX(y)=Meet(x-in-X)((y=>x)=>x)}
in terms of the nuclei~$\regnuc{x}$.

\thmskip

\begin{prop}\label{prop:regnuclei-meet-generate-all-nuclei}
Let\/ $L$ be a~frame.
The set\/ $\regnuc{L}$ of all nuclei $\regnuc{x}$, $x\in L\tinysp$,
	meet-gen\-er\-ates the complete lattice\/~$\Nuc(L)$.
If\/~$\nu$~is a~nucleus on\/~$L$ and\/ $X$~is a~subset of\/~$L\tinysp$,
then\/ $\nu=\nolinebreak\tMeet_{x\in X}\regnuc{x}\tinysp$
	\iff\/ $\dtinysp\fix(\nu)$ is the least nuclear system including the set\/~$X$;
in~particular, $\nu=\tMeet\tinysp\set{\regnuc{x}\narrt\suchthat x\narrt\in\fix(\nu)}$.
\qed
\end{prop}

\thmskip

Since $\regnuc{L}$ meet-generates the frame~$\Nuc(L)$,
it contains all completely meet-irreducible elements of~$\Nuc(L)$,
and in particular it contains all completely meet-prime elements of~$\Nuc(L)$.
We can locate the latter provided
	we know the completely meet-irreducible elements of the frame~$L\tinysp$.

\thmskip

\begin{prop}\label{prop:completely-meet-prime-in-NUC(L)}
Let\/ $L$ be a~frame and\/ $x\in L$.
The nucleus\/ $\regnuc{x}$ is completely meet-prime in\/~$\Nuc(L)$
\iff\/~$x$~is completely meet-irreducible in\/~$L\tinysp$.
\end{prop}

\interskip

\begin{myproof}
Suppose $x$ is completely meet-irreducible in~$L\tinysp$.
Let~$\nu_i\in\Nuc(L)$, $i\in I$, and suppose that $\regnuc{x}\geq\tMeet_{i\in I}\nu_i\tinysp$.
Then $x\in\fix(\tMeet_{i\in I}\nu_i)$, where
\begin{equation*}
\fix\bigl(\tMeet_{i\in I}\nu_i\bigr)
	\Eq \bigset{\tMeet_{i\in I}y_i\bigsuchthat\text{$y_i\in\fix(\nu_i)$ for $i\in I$}}
\end{equation*}
because the meet $\tMeet_{i\in I}\nu_i$ in~$\Nuc(L)$ is taken in~$\Cl(L)$.
It follows that $x=\tMeet_{i\in I}y_i$ for some $y_i\in\fix(\nu_i)$, $i\in I$.
Since $x$ is completely meet-irreducible in~$L$ we have $x=y_{i_0}\in\fix(\nu_{i_0})$ for some $i_0\in I$,
whence $\regnuc{x}\geq\nu_{i_0}$ because $\nu_{i_0}$ is a~nucleus.
This proves that $\regnuc{x}$ is completely meet-prime in $\Nuc(L)$.

Suppose $\regnuc{x}$ is completely meet-prime in $\Nuc(L)$,
and let $y_i\in L$, $i\in I$, be such that $x=\tMeet_{i\in I}y_i\tinysp$.
Then
\begin{align*}
\fix(\regnuc{x}) \Eq (L\bimpl x)
	&\Eq \bigset{\tMeet_{i\in I}(u\bimpl y_i)\bigsuchthat u\in L} \\ 
	&\Subseteq \bigset{\tMeet_{i\in I}z_i\bigsuchthat\text{$z_i\in\fix(\regnuc{y_i})$ for $i\in I$}}
		\Eq \fix\bigl(\tMeet_{i\in I}\regnuc{y_i}\bigr)~.
\end{align*}
That is, $\regnuc{x}\geq\tMeet_{i\in I}\regnuc{y_i}$.
Since $\regnuc{x}$ is completely meet-prime in $\Nuc(L)$,
we have $\regnuc{x}\geq\regnuc{y_{i_0}}$ for some $i_0\in I$, and hence~$x\geq y_{i_0}$.
Since also $x\leq y_{i_0}$, we have $x=y_{i_0}$.
This proves that $x$ is completely meet-irreducible in~$L\tinysp$.
\end{myproof}

\thmskip

The sub$\tinysp$-(complete lattice) $\Nuc(L)$ of~$\Cl(L)$ is not only an interior system in~$\Cl(L)$,
which gives us for every closure operator on~$L$ the largest nucleus below it,
it is also a~closure system in~$\Cl(L)$,
and so for any given closure operator on~$L$ there is the least nucleus above~it.
The following proposition has the details.

\thmskip

\begin{prop}\label{prop:the-least-nuc-above-clop}
Let\/ $\gamma$ be a~closure operator on a~frame\/~$L\tinysp$, and\/~$C=\fix(\gamma)$.
There exists the least nucleus\/~$\nu$ above\/~$\gamma\tinysp$.
The nucleus\/~$\nu$ and its fixpoint set are given by
\begin{equation*}
\nu \Eq \Meet\tinysp\bigset{\regnuc{x}\bigsuchthat x\in C,\, \gamma\leq \regnuc{x}}~,
	\qquad
\fix(\nu) \Eq \bigset{x\in C\bigsuchthat (L\bimpl x)\subseteq C}~,
\end{equation*}
where the meet in the formula for\/~$\nu$ is the pointwise meet.
\end{prop}

\interskip

\begin{proof}
The nucleus $\nu$ is the pointwise meet of all nuclei above~$\gamma\tinysp$.
Since $\Nuc(L)$~is meet-generated by the nuclei~$\regnuc{x}$,
it follows that
    $\nu = \tMeet\tinysp\set{\regnuc{x}\narrt\suchthat
				x\narrt\in L,\,\gamma\narrt\leq\regnuc{x}}\tinysp$.
Let $C\defeq\fix(\gamma)$ and $N\defeq\fix(\nu)$.
Since for every $x\in L$ we have
	$\gamma\leq \regnuc{x}$
	iff $\nu\leq \regnuc{x}$
	iff $x\in N$,
we get $N \defeq \set{x\narrt\in L \suchthat (L\bimpl x)\subseteq C\tinysp}\tinysp$.
Because $N\subseteq C$, we have also
	$N = \set{x\narrt\in C\suchthat (L\bimpl x)\subseteq C\tinysp}$
and~$\nu = \tMeet\tinysp\set{\regnuc{x}\narrt\suchthat x\in N}
	= \tMeet\tinysp\set{\regnuc{x}\narrt\suchthat
				x\narrt\in C,\,\gamma\narrt\leq\regnuc{x}}\tinysp$.
\end{proof}

\section{The \HMJ\ theorem}
\label{sec:HMJ-theorem}

In~\cite{Escardo-JFN-2003} the author demonstrates the utility of join induction
by using it in a~proof of the \HMJ\ theorem.
In this section we use the \HMJ\ theorem as a~training wheel
	on which we try out an application of the obverse induction principle.
The proof of the \HMJ\ theorem is spread through proofs of three lemmas,
	with parts of it reasoned out in the connecting text;
the short concluding reasoning then ties everything together.
The obverse induction principle gets its chance
	in the proof of Lemma~\ref{lem:scott-open-filter-is-nuclear},
where it performs admirably,
simplifying the proofs of the corresponding results
	in~\cite{Johnstone-VL&LSL-1985} and~\cite{Escardo-JFN-2003}
and shortening them to five easy lines of the proof proper
	(after the introductory~line).

\thmskip

\begin{thm}[Johnstone]\label{thm:HMJ-theorem}
The compact fitted quotient frames of any frame
are in order-reversing bijective correspondence%
\footnote{An ``order-reversing bijective correspondence'' means an antiisomorphism of posets,
that is, a~bijection between posets such that both the bijection itself and its inverse
are order-reversing.}
with the \Scottop\ filters of the frame.
\qed
\end{thm}

\thmskip

If this sounds all Greek to you, do not panic; 
everything will be explained below\,---\,slowly and in sickening detail\,---\,%
before we embark on the actual proof of the theorem,
which will be short and quite painless.

\txtskip

Frames we have already defined:
a~frame is a complete lattice in which finite meets distribute over arbitrary joins,
and a~frame morphism is a~mapping from a~frame to a~frame
that preserves finite meets and arbitrary joins.

\txtskip

Henceforward let $L$ be an arbitrary frame.

\txtskip

Consider a~nucleus $\gamma$ on~$L$.
The subposet $\gamma(L) = \fix(\gamma)$ of~$L$ is a~complete lattice
in which meets are calculated in~$L$
and the join of a~subset $S$ of $\gamma(L)$ is
$\tJoin^{\tinysp\gamma}\!S =\nolinebreak \gamma\bigl(\tJoinS\bigr)$.
Moreover, the infinite distributive law holds in the~complete lattice~$\gamma(L)$,
so~it~is in fact a~frame:
given any $x\in \gamma(L)$ and any $Y\subseteq\gamma(L)$, we have
\begin{equation*}
x\meet\tJoin^{\tinysp\gamma}\negtinysp Y
    \Eq \gamma(x)\meet\gamma\bigl(\tJoin\negtinysp Y\bigr)
    \Eq \gamma\bigl(x\meet\tJoin\negtinysp Y\bigr)
    \Eq \gamma\bigl(\tJoin_{\negdtinysp y\in Y}(x\meet\nolinebreak y)\bigr)
    \Eq \tJoin^{\tinysp\gamma}_{\negdtinysp y\in Y}(x\meet y)~.
\end{equation*}
The restriction $\gammapr\colon L\to\gamma(L)$ of~$\gamma$
preserves finite meets because $\gamma$ preserves~them,
and it preserves joins because $\gamma$ is a~closure~operator;
thus $\gammapr$~is a~surjective frame morphism.
This~is~why the~nuclear system~$\gamma(L)$
	is~also~called a~\knownas{quotient frame~of~$L\tinysp$}.

Let~$f\colon L\to K$ be a~morphism of frames.
Since~$f$ preserves all joins,
it has a~right adjoint $g\colon K\to L$, which preserves all meets,
thus the closure operator $\gamma \defeq gf$~on~$L$ preserves finite meets,
that is,~it is a~nucleus on~$L\tinysp$.
Now suppose that $f$ is surjective,
	and~hence $g$~is injective and~$fg=\id_K$.
Denoting by~$\gammapr\colon L\to\gamma(L)=g(K)$ the restriction of~$\gamma$
and by $h\colon g(K)\to f(L) = K$ the restriction of~$f$,
we have an~isomorphism~$h$ of posets and hence of frames
	such that~$h\tinysp\gammapr = f\negtinysp gf = f$.
Therefore, every surjective frame morphism from~$L$ is isomorphic
to an `inner' surjective frame morphism from~$L$ associated with a~nucleus on~$L\tinysp$.

\txtskip

A~frame is said to be \knownas{compact} if its top element is inaccessible by directed joins.
Spelled out: a frame $K$, with a~top element $\top$, is compact
\iff\ every directed subset~$S$ of~$K$ whose join is~$\top$ already contains~$\top$.

\txtskip

So we now know what is a~compact quotient frame.
``Fitted'' comes next.

\txtskip

Let $a\in L$.
The~principal ideal $\ldown a$ is a~frame,
the map $f_a\colon L\to\ldown a : x\mapsto x\meet a$ is a~surjective frame morphism,
and the defining property of $(\anon\bimpl\anon)$ shows that
the right adjoint of $f_a$ is the map $g_a\colon \ldown a\to L : y \mapsto (a\bimpl y)$;
the nucleus $a\open := g_a f_a$ on~$L$
maps $x\in L$ to $a\open(x) = \bigl(a\bimpl(x\meet a)\bigr) = (a\bimpl x)$.%
\footnote{\,We can verify directly that $a\open = (a\bimpl\anon)$ is a~nucleus:
it is ascending and increasing;
it is idempotent, $\bigl(a\bimpl(a\bimpl x)\bigr) = \bigl((a\meet a)\bimpl x) = (a\bimpl x)$;
it preserves binary meets, $\bigl(a\bimpl(x\meet y)\bigr) = (a\bimpl x)\meet(a\bimpl y)$.}
The nucleus $a\open$ is called the~\knownas{open nucleus} associated with~$a$;
the corresponding nuclear system is $a\open(L)=(a\bimpl L)$.
The restriction of the mapping $f_a$ to $(a\bimpl L)\to\ldown a$
is an isomorphism of posets and hence of frames.

A~nucleus $\gamma$ on~$L$,
	and the corresponding nuclear system~$\gamma(L)=\fix(\gamma)$, 
are said to be~\knownas{fitted},
if~$\gamma$ is a~join of open nuclei
	(with the join taken in the complete lattice~$\Nuc(L)$).\linebreak[3]
We shall denote by $\FitNuc\tinysp(L)$ the set of all fitted nuclei on~$L$,
and by $\FitNucSys(L)$ the set of all fitted nuclear systems
(that is,~fitted quotient frames)~on~$L$.
Subposet $\FitNuc\tinysp(L)$ of $\Nuc(L)$ is a~complete lattice
because it is evidently closed under joins~in~$\Nuc(L)$.%
\footnote{Actually, $\FitNuc\tinysp(L)$ is a~subframe of $\Nuc(L)$,
that is,~it is also closed under finite meets,
since $\bot\!\open\colon x\mapsto\top$ is the greatest nucleus,
and $a\open\!\meet b\open = (a\join b)\open$ for all $a,\,b\in L$.}\linebreak[3]
Correspondingly, the subposet $\FitNucSys(L)$ of $\NucSys(L)$ is a~complete lattice;
it~is closed~under meets in~$\NucSys(L)$,
and since meets in~$\NucSys(L)$ are intersections,
$\FitNucSys(L)$ is a~closure system in~$\PowL$.

Below any nucleus $\gamma\in\Nuc(L)$
there exists the greatest fitted nucleus ${\gamma\tinysp}^\fitit\in\FitNuc\tinysp(L)$,
which is the join of all open nuclei below~$\gamma\tinysp$.
The~mapping $\gamma\mapsto{\gamma\tinysp}^\fitit$ is an interior operator on~$\Nuc(L)$;
it is fittingly called the~\knownas{fitting} of nuclei on~$L\tinysp$.

\txtskip

Given a~nucleus $\gamma$ on~$L$, which open nuclei on~$L$ are below~$\gamma\tinysp$?
Lemma~\ref{lem:open-nucs-below-a-nuc} has the answer.
In the proof of this lemma we are going to use the following
inequality satisfied by an endomap~$f$ on the frame~$L$
	that preserves binary meets, and hence is increasing:
for all $x,\tinysp y\in L$, $f(x\bimpl y)\leq \bigl(f(x)\bimpl\nolinebreak f(y)\bigr)$.
The inequality follows from the inequality between the first and the last expressions in
$f(x)\meet f(x\bimpl y) = f\bigl(x\meet(x\bimpl y)\bigr) \leq f(y)$.

\thmskip

\begin{lem}\label{lem:open-nucs-below-a-nuc}
Let\/~$L$ be a~frame.
If\/~$a\in L$ and\/~$\gamma\in\Nuc(L)$, then\/~$a\open\leq\gamma$ iff\/ $\gamma(a)=\top$.
\end{lem}

\interskip

\begin{myproof}
If~$a\open\leq\gamma$, then $\top = (a\bimpl a) = a\open(a) \leq \gamma(a)$.
Conversely, if $\gamma(a) = \top$,
then for every $x\in L$,
$a\open(x) = (a\bimpl x) \leq \gamma(a\bimpl x) \leq \bigl(\gamma(a)\bimpl\gamma(x)\bigr)
= \bigl(\top\bimpl\gamma(x)\bigr) = \gamma(x)$.
\end{myproof}

\thmskip

We can now write down the following formula for the fitting of a nucleus:
\begin{equation}\label{eq:fitting-formula}
{\gamma\tinysp}^\fitit \Eq \tJoin\set{a\open\suchthat\gamma(a)=\top}
	\qquad\quad \text{for every \ $\gamma\in\Nuc(L)$}~.
\end{equation}

\negdisplayshortskip
\thmskip

So far we completely understand one side of the bijection mentioned in Theorem~\ref{thm:HMJ-theorem}.
There is not much left to understand on the other side.

A~\knownas{filter} of a~poset $P$ is a~downward directed upper set of~$P$.
In~our frame~$L$ a~filter is an upper set closed under finite meets
(including the empty meet, that is,~a~filter always contains~$\top$).
Every~filter $V\negdtinysp$ of~$L$
obeys the \knownas{modus~ponens} rule:
for~all~$a,\,b\in L$,
	if $a\in V\negdtinysp$ and $(a\bimpl b)\in V\negdtinysp$,
	then $b\in V\negdtinysp$ because $b \geq a\meet(a\bimpl b) \in V\negdtinysp$.

For any nucleus $\gamma$ on~$L$,
the set $\gamma^{-1}(\top)$ is a~filter of~$L$;
we~shall call filters of this form \notion{nuclear filters of~$L$},
and will denote by $\NucFilt(L)$
the poset of all nuclear filters of~$L$ ordered by~inclusion.
Since
\begin{equation*}
\bigl(\tMeet\negdtinysp\Gamma\bigr)^{-1}(\top)
    \widedt\Eq \tInters\dtinysp\set{\gamma^{-1}(\top)\suchthat \gamma\in\Gamma}
    \qquad\quad \text{for \ $\Gamma\subseteq\Nuc(L)$}
\end{equation*}
(recall that all meets of nuclei are calculated pointwise),
it follows that $\NucFilt(L)$ is a~closure system in~$\PowL\tinysp$.
Indeed, given a~set $\coll{V}$ of nuclear filters,
let $\Gamma$ be the set of all nuclei $\gamma$ such that $\gamma^{-1}(\top)\in\coll{V}$.
Since every filter in $\coll{V}$ is of the form $\gamma^{-1}(\top)$
for some nucleus $\gamma$ in $\Gamma$,
the intersection
$\tInters\negdtinysp\coll{V}
	= \tInters_{\gamma\in\Gamma}\gamma^{-1}(\top)
	= (\tMeet\negdtinysp\Gamma)^{-1}(\top)$
is a~nuclear filter.

By definition, a~\Scottop\ subset of a~poset is an upper set inaccessible by directed joins.
Since every~filter is an upper set by definition,
a~filter is \Scottop\ iff it is inaccessible by directed joins.

\txtskip

We have everything ready to relate compact fitted quotient frames to \Scottop\ filters.
The following lemma is Lemma~4.4 in~\cite{Escardo-JFN-2003},
which in turn is Lemma~3.4(i) in~\cite{Johnstone-VL&LSL-1985};
its proof is almost verbatim as in~\cite{Escardo-JFN-2003},
which in turn is lifted from~\cite{Johnstone-VL&LSL-1985}.
Anyway, this lemma is not very deep,
it is an~immediate consequence
	of the relationship between joins in a~frame and joins in a~quotient frame of~the~frame.

\thmskip

\begin{lem}\label{lem:compact-qr-vs-scottopen-nf}
Let\/ $\gamma$ be a nucleus on a frame\/~$L$.
Then the quotient frame\/ $\gamma(L)$ is compact
\iff\ the nuclear filter\/ $\gamma^{-1}(\top)$ is \Scottop.
\end{lem}

\interskip

\begin{myproof}
($\implies$)\, Suppose $\gamma(L)$ is compact,
and let $S\subseteq L$ be directed with $\tJoinS \in \gamma^{-1}(\top)$.
Since $\tJoin^{\tinysp\gamma}\!\gamma(S) = \gamma\bigl(\tJoinS\bigr) = \top$ and $\gamma(L)$ is compact,
there is some $s\in S$ with $\gamma(s) = \top$, that is,~with $s\in\gamma^{-1}(\top)$.

($\isimplied$)\, Suppose $\gamma^{-1}(\top)$ is \Scottop,
and let $S\subseteq\gamma(L)$ be directed with $\tJoin^{\tinysp\gamma}\!S = \top$.
Since $\tJoin^{\tinysp\gamma}\!S = \gamma\bigl(\tJoinS\bigr)$, we have $\tJoinS\in\gamma^{-1}(\top)$,
and since $\gamma^{-1}(\top)$ is \Scottop,
there is some $s\in S$ such that $s\in\gamma^{-1}(\top)$,
that is, such that $s = \gamma(s) = \top$.
\end{myproof}

The fitting of nuclei, the interior operator $\gamma\mapsto\gamma^\fitit$ on $\Nuc(L)$,
is a~counit of a~certain (covariant) Galois connection,%
\footnote{\,Strictly speaking, the counit is the family of relationships
	$\gamma^\fitit\leq\gamma$ for $\gamma\in\Nuc(L)$.
Mind that in a~poset the relationships $x\leq y$
	are the morphisms $x\to y$ of the poset regarded as a~category.}
which we now proceed to describe.

For~every $\gamma\in\Nuc(L)$ put $\oneker\negtinysp\gamma := \gamma^{-1}(\top)$,
and for every $S\in\PowL$ put $\fitnuc\tinysp S := \tJoin_{\!s\in S}s\open$.
For any $\gamma\in\Nuc(L)$ and any $S\in\PowL$
the chain of equivalences
\begin{align*}
\qquad\qquad\qquad
\fitnuc\tinysp S \Leq \gamma
    & \wide\Isequiv (\forall s\narr\in S)(s\open\leq\gamma) \\ 
    & \wide\Isequiv (\forall s\narr\in S)(\gamma(s) = \top)
		&& \text{(by Lemma~\ref{lem:open-nucs-below-a-nuc})}
							\qquad\qquad \\ 
    & \wide\Isequiv S\Subseteq\oneker\negtinysp\gamma
\end{align*}
shows that $\pair{\fitnuc,\negdtinysp\oneker}$
	is a~Galois connection $\PowL\galoisconn\Nuc(L)\tinysp$.
By our definitions, $\fitnuc\tinysp\PowL$
	is the~set $\FitNuc\tinysp(L)$ of all fitted nuclei on~$L$,
while $\oneker\negtinysp\Nuc(L)$ is the set $\NucFilt(L)$ of all nuclear filters of~$L$.
From the general properties of Galois connections it follows that
$\fit\defeq\fitnuc\negtinysp\oneker$ is an interior operator on~$\Nuc(L)$
and that for any~nucleus~$\gamma$,
$\fit(\gamma) =\nolinebreak \gamma^\fitit$~(see~\eqref{eq:fitting-formula})
	is~the~greatest fitted nucleus below~$\gamma$
		(all~of~which we already~know),
while on the other side,
$\dtinysp\nucfilt \defeq \oneker\negtinysp\fitnuc$ 
	is a~closure operator on~$\PowL\tinysp$,
where for any subset~$S$~of~$L$,
$\dtinysp\nucfilt(S)$ is the least nuclear filter of~$L$ that includes~$S$.
We~have the identities $\oneker\negtinysp\fitnuc\negtinysp\oneker = \oneker$
and $\fitnuc\negtinysp\oneker\negtinysp\fitnuc = \fitnuc$,
meaning, respectively, that
$\bigl(\gamma^\fitit\bigr)^{-1}(\top) = \gamma^{-1}(\top)$ for every $\gamma\in\Nuc(L)$,
and that
$\tJoin\dtinysp\set{s\open\negtinysp\narrdt\suchthat s\narrt\in\nucfilt(S)}
	= \tJoin\dtinysp\set{s\open\negtinysp\narrdt\suchthat s\narrt\in S}$
for every $S\subseteq L$.
And, restricting $\fitnuc$ to $\NucFilt(L)\to\FitNuc\tinysp(L)$
and $\oneker$ to $\FitNuc\tinysp(L)\to\NucFilt(L)$,
we obtain two isomorphisms of complete lattices which are inverses to each other.

\txtskip

At last, here comes the punch line\,---\,or should it be the~punch lemma?

The following lemma is Lemma 3.4(ii) in~\cite{Johnstone-VL&LSL-1985},
reappearing as Lemma~4.3(2) in~\cite{Escardo-JFN-2003}.
We~give a~short proof which uses the obverse induction principle
	instead of transfinite induction in~\cite{Johnstone-VL&LSL-1985},
	and instead of the join induction in~\cite{Escardo-JFN-2003}.

\thmskip

\begin{lem}\label{lem:scott-open-filter-is-nuclear}
Every \Scottop\ filter of a~frame\/ $L$ is nuclear.
\end{lem}

\interskip

\begin{myproof}
Let $V\negdtinysp$ be a~\Scottop\ filter of a~frame $L$,
and let
	$\gamma = \fitnuc\negtinysp V\negdtinysp 
		= \tJoin\dtinysp\set{v\open\negdtinysp\suchthat v\in V}\tinysp$.
The filter $V$ is \dirinacc.
If $v\in V\negdtinysp$ and $x\in L$, and $v\open(x) = (v\bimpl x) \in V$,
then~$x\in V$ by modus ponens, meaning that~$V\negdtinysp$
is inversely closed under~$\set{v\open\suchthat v\in V}$.
Invoking~the~obverse induction principle we find that $V$~is~inversely closed under~$\gamma$,
so~certainly
$\nucfilt(V)
	= \oneker\negtinysp\gamma
	= \gamma^{-1}(\top) \subseteq V\negdtinysp$.
Since also $V\negdtinysp\subseteq\nucfilt(V)$,
we conclude that $V\negdtinysp=\nucfilt(V)$ is a~nuclear~filter.
\end{myproof}

\thmskip

After all the preparations,
Theorem~\ref{thm:HMJ-theorem} is easy to prove.

\thmskip

\begin{myproof}[Proof of Theorem~\ref{thm:HMJ-theorem}.]
Let~$\coll{F}$ be the poset of all \Scottop\ filters of~$L$ ordered by inclusion,
let $\coll{Q}$ be the subposet of $\FitNucSys(L)$
	consisting of all compact fitted quotient frames on~$L$,
and let $\coll{G}$ be the subposet of $\FitNuc(L)$
	consisting of all nuclei $\gamma$ on~$L$ such that $\gamma(\negtinysp P)\in\coll{Q}\dtinysp$.
Now Lemma~\ref{lem:scott-open-filter-is-nuclear} and Lemma~\ref{lem:compact-qr-vs-scottopen-nf}
tell us that $\coll{F}$ is a~subposet of $\NucFilt(L)$
and that the~isomorphisms of complete lattices
\begin{equation*}
\NucFilt(L) \longto \FitNuc(L) \longto \FitNucSys(L)^\op
	\wide: V \longmapsto \fitnuc\negtinysp V \longmapsto \fix(\fitnuc\negtinysp V)
\end{equation*}
restrict to isomorphisms of posets $\coll{F} \to \coll{G} \to \coll{Q}^\op$.
\end{myproof}

\section{Doing it with maximal elements}
\label{sec:do-it-with-max-elems}

Let $P$ be a~dcpo.
Since the poset $\ClSys(\negtinysp P)$ of all closure systems in~$P$
	is a~closure system in $\PowP\negtinysp$,
it is determined by a~set of closure rules on~$P\negtinysp$.
One such set of closure rules is, of~course, the full-fledged closure theory
consisting of all closure rules obeyed by~$\ClSys(\negtinysp P)\tinysp$.
But this closure theory is too large;
we want some smaller set of closure rules
	that determines the closure system $\ClSys(\negtinysp P)\tinysp$,
and moreover,
we want a~set of closure rules
	which can be described in terms of the structure of the dcpo~$P\negtinysp$.

We shall obtain a~suitable set of closure rules using the approach
	in the paper~\cite{Ranzato-CCPOFCL-1999}.
We~will not follow the exposition in the paper;
our treatment will be more streamlined,	and we will obtain some results that are not in the paper.

\txtskip

A~\notion{default closure rule} on a~poset~$P$ 
is a~closure rule $B\negdtinysp\adjoins{}c$ on the set~$P$
	(that is, $B\subseteq P$ and $c\in P\tinysp$)
where $c$ is a~maximal lower bound of~$B$.
We shall denote by $\Rulesdflt(\negtinysp P)$
	the set of all default closure rules on a~poset~$P\negtinysp$.
We shall write $B\negdtinysp\adjoins{\dflt}c$
	to mean that $B\negdtinysp\adjoins{}c$~is a~default closure rule,
that is, that $\Rulesdflt(\negtinysp P)\colon B\negdtinysp\adjoins{}c\dtinysp$.

Default closure rules on a~poset~$P$
generalize the default closure rules $B\adjoins{}\tMeet\!B$ on a~complete lattice $L$,
where~$B$ is an arbitrary subset of~$L\tinysp$.

\txtskip

Let $P$ be a~poset.
A~default closure rule $B\negdtinysp\adjoins{}c$ on~$P$ can be reflexive,
	which means that it has~$c\in B\dtinysp$;
it is reflexive iff $c$ is the least element of~$B$.
In~more detail: let~$B\negdtinysp\adjoins{}c$~be a~closure rule on~$P\tinysp$;
if~$B\negdtinysp\adjoins{\dflt}c$ and~$c\in B$,
	then $c$~is the least element of~$B\tinysp$;
conversely, if~$B$~has a~least element~$c\tinysp$,
	then $B\negdtinysp\adjoins{}c$ is the unique default closure rule with~the~body~$B$.

\thmskip

\begin{lem}\label{lem:f-preclosure-on-P=>Fix(f)-obeys-Rdf(P)}
Let\/ $P$ be a~poset.
If\/ $f$ is a~preclosure map on\/~$P$
	then\/ $\dtinysp\fix(f)$ obeys\/~$\Rulesdflt(\negtinysp P)$.
\end{lem}

\interskip

\begin{myproof}
Let $f$ be a~preclosure map,
and suppose that $B\negdtinysp\adjoins{\dflt}c$ with $B\subseteq\fix(f)\tinysp$.
For~any~$b\in B$ we have $f(c)\leq f(b) = b$, thus $f(c)$~is a~lower bound of~$B$.
Since $c\leq f(c)$ and $c$ is a~maximal lower bound of~$B$,
it~follows that~$f(c)=c\in\fix(f)\tinysp$.
\end{myproof}

\thmskip

Let $P$ be a~poset.

We~shall say that $P$ \notion{has a~ceiling}
if for every element~$x$ of~$P$ there exists a~maximal element~$y$ of~$P$
	such that~$x\leq y\tinysp$.%
\footnote{\,When $P$ has a~ceiling,
the set of all maximal elements of~$P$ is \emph{the} ceiling we have in mind here.}
Mark that the empty poset has a~ceiling.
We shall say that a~subset~$A$ of~$P$ has a~ceiling
	if the subposet~$A$ of~$P$ has a~ceiling.

We shall say that $P$ is \notion{default-enabled}
if for every subset $X$ of $P$ the set of all lower bounds of $X$ in $P$ has a~ceiling
	(that is, every lower bound of $X$ is below some maximal lower bound of~$X\tinysp$).%
\footnote{A~default-enabled poset is in~\cite{Ranzato-CCPOFCL-1999}
called a~\knownas{relatively maximal lower bound complete} poset,
which is rather a~mouthful, so Ranzato shortens it to \knownas{rmlb$\tinysp$-complete} poset,
which is not very mnemonic.}
If $P$ is default-enabled, then in particular the set $P$ itself has a~ceiling,
	since $P$ is the set of all lower bounds of the empty subset.

\txtskip

The following lemma tells us that a~default-enabled poset has enough default closure rules
to determine the closure systems in the poset.

\thmskip

\begin{lem}\label{lem:dflt-enabled-P-has-enough-dflt-clrules}
Let\/ $P$ be a~default-enabled poset.
If a~subset\/~$C$ of\/~$P$ obeys\/~$\Rulesdflt(\negtinysp P)$,
	then\/ $C$~is a~closure system in\/~$P\negtinysp$.
\end{lem}

\interskip

\begin{myproof}
Suppose $C\subseteq P$ obeys $\Rulesdflt(\negtinysp P)$.
Let $x$ be any element of~$P\tinysp$;
we shall prove that the set~$B\defeq\allabove{C}{x}$ has a~least element.
The~ele\-ment~$x$ is a~lower bound of~$B$,
thus $x\leq u$ for some maximal lower bound~$u$ of~$B$ because $P$~is default-enabled.
Then $B\negdtinysp\adjoins{\dflt}u\dtinysp$,
therefore $u\in C$ because $C$~obeys~$\Rulesdflt(\negtinysp P)\tinysp$,
whence $u\in\allabove{C}{x} = B$ is the least element of~$B$.
\end{myproof}

\thmskip

The proofs of Lemma~\ref{lem:f-preclosure-on-P=>Fix(f)-obeys-Rdf(P)}
	and Lemma~\ref{lem:dflt-enabled-P-has-enough-dflt-clrules}
correspond to the two parts of the proof of Theorem~4.4 in~\cite{Ranzato-CCPOFCL-1999}
(where Lemma~\ref{lem:f-preclosure-on-P=>Fix(f)-obeys-Rdf(P)}
	is slightly more general than the first part of Theorem~4.4).

\txtskip

We have the following consequence of
	Lemma~\ref{lem:f-preclosure-on-P=>Fix(f)-obeys-Rdf(P)}
		and Lemma~\ref{lem:dflt-enabled-P-has-enough-dflt-clrules}:

\thmskip

\begin{prop}\label{prop:P-dflt-enabled=>ClSys(P)-clsys-in-Pow(P)-etc}
Let\/ $P$ be a~default-enabled poset.

If\/ $g$ is a~preclosure map on\/~$P$,
then\/ $\dtinysp\fix(g)$ is a~closure system in\/~$P$
and the closure operator\/~$h$ on\/~$P$ with\/ $\fix(h)=\fix(g)$
	is the least closure operator above\/~$g$.

A~subset of\/ $P$ is a~closure system~in\/~$P$ \iff\ it obeys\/~$\Rulesdflt(\negtinysp P)$.
Consequently,
the set\/~$\ClSys(\negtinysp P)$ of all closure systems in\/~$P$
is a~closure system in\/~$\PowP\negtinysp$,
so~it~is a~complete lattice in which the meets are intersections.
Also the poset\/~$\Cl(\negtinysp P)$ of all closure operators on\/~$P\negtinysp$,
being antiisomorphic to the poset\/~$\ClSys(\negtinysp P)$,
	is a~complete lattice.
\end{prop}

\interskip

\begin{myproof}
Only the assertion about the preclosure map $g$ needs any proving.
In any poset, if $g$ is a~preclosure map,
then a~closure operator $h$ is the least closure operator above~$g$
iff~$\dtinysp\fix(h)$~is the largest of all closure systems included in~$\dtinysp\fix(g)$.
If~$g$~is a~preclosure map on the default-enabled poset~$P$,
then $\dtinysp\fix(g)$, which obeys~$\Rulesdflt(\negtinysp P)$,
is itself the largest closure system in~$P$ included in~$\dtinysp\fix(g)$,
and the assertion follows.
\end{myproof}

\thmskip

We~shall cook up a~theorem for default-enabled posets
	that will resemble Theorem~\ref{thm:in-dcpo-cl-generd-by-precls-&-induct} for dcpos.
With this aim in mind we introduce the following notion:

\txtskip

Let us say that a~subset $A$ of a~poset $P$ is \notion{default-enabled within~$P$}
if it satisfies the following two conditions:
\begin{itemize}[topsep=.5ex,itemsep=.5ex,leftmargin=4em]
\item[{\rm(i)}\:] the~subposet~$A$ is default-enabled;
\item[{\rm(ii)}\:] for~every~$x\in P$ the set $\allbelow{A}{x}$ has a~ceiling.
\end{itemize}
The condition~(i) is a~property of the structure of the subposet $A$ alone,
	independent of the rest of the structure of the ambient poset~$P\negtinysp$,
while the condition~(ii) prescribes how the subposet $A$ has to be situated inside the poset~$P\negtinysp$.

\txtskip

And here is the theorem mimicking Theorem~\ref{thm:in-dcpo-cl-generd-by-precls-&-induct}
	(minus the obverse induction principle);
with it we wander a~little way beyond~\cite{Ranzato-CCPOFCL-1999}.

\thmskip

\begin{thm}\label{thm:in-dfltenab-cl-generd-by-precls-&-induct}
Let\/ $P$ be a~default-enabled poset,
and let\/ $G$ be a~set of preclosure maps on\/~$P\negtinysp$.
The set\/~$\dtinysp\fix(\leftt G)$ is a~closure system in\/~$P\negtinysp$,
and the closure operator\/~$\overbar{G}$ on\/~$P$
	which has\/~$\dtinysp\fix\bigl(\overbar{G}\dtinysp\bigr)=\fix(\leftt G)$
is the least closure operator on $P$ that is above\/~$G$.

The following {\bfseries induction principle} holds:
if a~subset\/~$A$ of\/~$P$ is default-enabled within\/~$P$ and is closed under\/~$G$,
then it is closed under\/~$\overbar{G}\tinysp$.
\end{thm}

\interskip

\begin{myproof}
For every $g\in G$ the fixpoint set $\dtinysp\fix(g)$ obeys $\Rulesdflt(\negtinysp P)$
	by Lemma~\ref{lem:f-preclosure-on-P=>Fix(f)-obeys-Rdf(P)},
therefore $\dtinysp\fix(g)\in\ClSys(\negtinysp P)$ 
	by Lemma~\ref{lem:dflt-enabled-P-has-enough-dflt-clrules}.
Since $\ClSys(\negtinysp P)$ is closed under arbitrary intersections,
	by~Proposition~\ref{prop:P-dflt-enabled=>ClSys(P)-clsys-in-Pow(P)-etc},
the set $\dtinysp\fix(\leftt G)=\tInters_{g\in G}\fix(g)$ is~a~closure system;
let $h$ be the closure operator on~$P$ which has $\dtinysp\fix(h)=\fix(\leftt G)$.
If $g\in G$, then $\dtinysp\fix(h)\subseteq\fix(g)$,
	whence $h\geq g\dtinysp$; that is, $h\geq G$.
Let $k$ be a~closure operator above~$G$.
Then $\dtinysp\fix(k)\subseteq\fix(g)$ for every $g\in G$,
	thus $\dtinysp\fix(k)\subseteq\fix(\leftt G)=\fix(h)$,
	and so $k\geq h$.
The closure operator $\overbar{G}\defeq h$
	has the properties stated in the theorem.

\emph{The induction principle.}

Assume that $A\subseteq P$ is default-enabled within~$P$ and closed under~$G$.

Let $G_{\negdtinysp A}$ be the set of restrictions
	$g_{\negtinysp A}\colon A\to A$ of the maps~$g\in G$.
We obtained a~set~$G_{\negdtinysp A}$ of preclosure maps on a~default-enabled subposet~$A\tinysp$,
thus there is (by the first part of the proof above, applied to the poset~$A\tinysp$)
a~closure operator $\hpr$ on~$A$ such that $\dtinysp\fix(\hpr)=\fix(G_{\negdtinysp A})$,
where $\dtinysp\fix(G_{\negdtinysp A}) = A\inters\fix(\leftt G) = A\inters\fix(h)\subseteq\fix(h)$.
If~$a\in A\tinysp$, then~$a\leq\hpr(a)\in\fix(\hpr)=\fix(G_{\negdtinysp A})\subseteq\fix(h)$,
thus $h(a)\leq h(\hpr(a))=\hpr(a)$.

We shall show that for any $a\in A$ also $h(a)\geq\hpr(a)$,
and therefore $h(a)=\hpr(a)\in A\tinysp$.

So let $a\in A\tinysp$.
By assumption $\allbelow{A}{h(a)}$ has a~ceiling.
Since $a\in \allbelow{A}{h(a)}\tinysp$,
there exists in $\allbelow{A}{h(a)}$ a~maximal element $\apr$ such that $a\leq\apr$.
For every $g\in G$ we have $g(\apr)\leq g(h(a))=h(a)$ and $g(\apr)\in A\tinysp$,
thus $g(\apr)\in \allbelow{A}{h(a)}\tinysp$;
now since $\apr\leq g(\apr)$ and $\apr$~is maximal in $\allbelow{A}{h(a)}\tinysp$,
it follows that $g(\apr)=\apr$.
Thus we have $a\leq\apr$ in~$A\tinysp$,
where $\apr$~is fixed by~$g_{\negtinysp A}$ for every $g\in G$,
therefore $\apr$ is fixed by~$\hpr$,
and it follows that $\hpr(a)\leq\hpr(\apr)=\apr\leq h(a)$.
\end{myproof}

\thmskip

Every dcpo is a~default-enabled poset.
First, every nonempty dcpo has a~maximal element, by Zorn's lemma.
Next, if~$P$~is a~dcpo and $x$~is any element of~$P\negtinysp$,
then the principal filter~$\lup x$ is a~sub-dcpo of~$P$
and hence has a~maximal element which is also a~maximal element of~$P\dtinysp$;
it~follows that $P$~has a~ceiling.
Finally, if $X$ is any subset of a~dcpo~$P\negtinysp$,
then the set of all lower bounds of~$X$ in~$P$ is a~sub-dcpo of~$P\negtinysp$,
thus it has a~ceiling,
and we see that $P$~is default-enabled.

Suppose that a~subset $A$ of a~dcpo $P$ is closed under directed joins in~$P\negtinysp$.
If~$x$~is any element of~$A$,
then $\allbelow{A}{x} = A \inters \ldown x$
	is the intersection of two sub-dcpos of~$P\negtinysp$,
thus it is a~sub-dcpo of~$P\negtinysp$,
so~it~has a~ceiling.
The sub-dcpo~$A$ is default-enabled within~$P\negtinysp$.

Theorem~\ref{thm:in-dfltenab-cl-generd-by-precls-&-induct}
thus specializes to Theorem~\ref{thm:in-dcpo-cl-generd-by-precls-&-induct}
	(minus the obverse induction principle),
but we need the~axiom of~choice to do~it.
In fact we cannot do the specialization without involving the~axiom of~choice,
since it is easy to prove that the assertion
	that every nonempty~dcpo has a~maximal element
implies, in~the~theory of sets without~the~axiom of~choice,
	the Hausdorff's maximal chain condition
	(to see this, consider the pointed dcpo of all chains in a~poset).

Luckily we do not need the help of the~axiom of~choice
	in order to specialize Theorem~\ref{thm:in-dfltenab-cl-generd-by-precls-&-induct}
		to Theorem~\ref{thm:in-dcpo-cl-generd-by-precls-&-induct},
since we already proved the latter theorem on its own.

\txtskip

The obvious question to ask at this point is
whether the class of default-enabled posets is strictly larger than the class of dcpos.
The answer is yes, it is strictly larger:
the~poset~$P_2$ in Figure~\ref{fig:dfltenabd-post-not-dcpo}
\begin{figure}[!htp]\centering
\vspace{.5ex}
\includegraphics[draft=false]{mp/clops-on-dcpos-2.mps}
\vspace{-.5ex}
\mycaption{\label{fig:dfltenabd-post-not-dcpo}%
	A default-enabled poset which is not a dcpo.}
\vspace{-1ex}
\end{figure}
	(reproduced from~\cite{Ranzato-CCPOFCL-1999})
is default-enabled while it is not a~dcpo.
This poset~$P_2\tinysp$, though it answers the question in the affirmative,
is not very exciting, since the only closure operator on it is the identity map.
Here is a~challenge: describe a~class of \emph{interesting} default-enabled posets
that are far~from~being~dcpos%
\footnote{\,Which means no cheap tricks.
For example, we can place the poset $P_2$ on top of any dcpo
and obtain a~default-enabled poset which is not a~dcpo\,---\,%
but such a~poset is as uninteresting as the poset~$P_2\tinysp$.}
and whose complete lattices of closure operators are quite nontrivial.

\txtskip
\pagebreak[3]

Now our travels will carry us beyond the horizon of~\cite{Ranzato-CCPOFCL-1999}.

Let $P$ be a~meet-semilattice.

For any two elements $a$ and $b$ of $P$
we define the set $(a\astbimpl b) \defeq \set{x\in P\suchthat x\narrt\meet a\leq b}$,
and then define the set $(a\dotbimpl b)$
	as the set of all maximal elements of the set~$(a\astbimpl b)$.
Note that the set $(a\astbimpl b)$ is always nonempty as it contains the element~$b\dtinysp$;
however, $(a\astbimpl b)$ may not have any maximal elements,
	so it is possible that the set $(a\dotbimpl b)$ is empty.

A~\notion{nuclear closure rule} on~$P$ is a~unary closure rule $b\adjoins{}c\tinysp$,
where $c\in(a\dotbimpl b)$ for some~$a\in P$.
The set of all nuclear closure rules on $P$ shall be denoted by~$\Rulesnuc(\negtinysp P)$.
We~shall write $b\adjoins{\nuc}c$
	to mean that the closure rule $b\adjoins{}c$ is nuclear;
that~is, $b\adjoins{\nuc}c$ is synonymous with $\Rulesnuc(\negtinysp P)\colon b\adjoins{}c$.
A~subset~$X$ of~$P$ obeys~$\Rulesnuc(\negtinysp P)$
	\iff\ $(a\dotbimpl x)\subseteq X$ for all~$a\in P$ and all~$x\in X$.

\thmskip

\begin{lem}\label{lem:Fix(prenuc)-obeys-nuclear-rules}
If\/ $\gamma$ is a~prenucleus on a~meet-semilattice\/~$P$,
	then\/ $\dtinysp\fix(\gamma)$ obeys\/~$\Rulesnuc(\negtinysp P)$.
\end{lem}

\interskip

\begin{myproof}
Suppose that $b\adjoins{\nuc}c$ with $b\in\fix(\gamma)\dtinysp$;
	we shall prove that $c\in\fix(\gamma)$.
There exists $a\in P$ such that $c\in(a\dotbimpl b)$.
Since $c\meet a\leq b$ and $\gamma(b)=b\tinysp$, we have
\begin{equation*}
\gamma(c)\meet a
	\Leq \gamma(c)\meet\gamma(a)
	\Eq \gamma(c\narrt\meet a)
	\Leq \gamma(b)
	\Eq b~,
\end{equation*}
thus $\gamma(c)\in(a\astbimpl b)$.
Since $c\leq\gamma(c)$ and $c$ is maximal in $(a\astbimpl b)$, we have $\gamma(c)=c\tinysp$.
\end{myproof}

\thmskip

A~sort of converse of Lemma~\ref{lem:Fix(prenuc)-obeys-nuclear-rules} holds
if for all elements $a$, $b$ of a~meet-semi\-lattice~$P$
	the set $(a\astbimpl b)$ has a~ceiling.

\thmskip

\begin{lem}\label{lem:enough-nucrules-clop-obeys-nucrules=>clop-is-nuc}
Let\/ $P$ be a~meet-semilattice in which every set\/ $(a\astbimpl b)$ with\/~$a,\tinysp b\in P$
	has~a~ceiling.
Let\/~$\gamma$~be a~closure operator on\/~$P$.
If\/~$(a\dotbimpl b)\subseteq\fix(\gamma)$ for~all\/~$a\in P$ and all\/~$b\in\fix(\gamma)$,
then\/ $\gamma$ preserves binary meets, that is, it is a~nucleus.
\end{lem}

\interskip

\begin{myproof}
Let $a,\tinysp b\in P$.
The inequality $\gamma(a)\meet\gamma(b)\geq\gamma(a\narrt\meet b)$ holds
	since $\gamma$ is increasing.
For the converse inequality
it suffices to prove that $a\meet\gamma(b)\leq\gamma(a\narrt\meet b)$.
Since $a\meet b\leq \gamma(a\narrt\meet b)$,
	$b$ lies in $\bigl(a\astbimpl\gamma(a\narrt\meet b)\bigr)$.
Since $\bigl(a\astbimpl\gamma(a\narrt\meet b)\bigr)$ has a~ceiling,
there exists $c\in\negdtinysp \bigl(a\dotbimpl\gamma(a\narrt\meet b)\bigr)$
	such that $b\leq c\tinysp$.
Since $\gamma(a\narrt\meet b)\in\fix(\gamma)$,
it follows from our assumption about $\dtinysp\fix(\gamma)$ that $c\in\fix(\gamma)$,
therefore $\gamma(b) \leq \gamma(c) = c\tinysp$,
and we conclude that $a\meet\gamma(b) \leq a\meet c \leq \gamma(a\narrt\meet b)$.
\end{myproof}

\thmskip

Let us say that a meet-semilattice $P$ is \notion{nuclear-enabled}
if it is default-enabled and every set $(a\astbimpl b)$ with $a,\tinysp b\in P$ has a~ceiling.
The following proposition is a~consequence of Lemma~\ref{lem:Fix(prenuc)-obeys-nuclear-rules}
and Lemma~\ref{lem:enough-nucrules-clop-obeys-nucrules=>clop-is-nuc}.

\thmskip

\begin{prop}\label{prop:nuc-enab-meet-semilatt}
If\/ $P$ is a~nuclear-enabled meet-semilattice,
then\/ $\NucSys(\negtinysp P)$ is a~closure system in\/ $\PowP$
	determined by the set of closure rules\/
		 $\Rulesdflt(\negtinysp P)\union\Rulesnuc(\negtinysp P)\tinysp$.
\qed
\end{prop}

\thmskip

Let $P$ be a~nuclear-enabled meet-semilattice.

The poset $\NucSys(\negtinysp P)$, 
	being a~closure system in the complete lattice $\PowP$,
is a~complete lattice; it is also a~closure system in the complete lattice $\ClSys(\negtinysp P)$.
The meets in $\NucSys(\negtinysp P)$, as well as in $\ClSys(\negtinysp P)$, are intersections.
Correspondingly, $\Nuc(\negtinysp P)$ is an interior system in~$\Cl(\negtinysp P)$,
	and hence is a~complete lattice with joins inherited
		from the complete lattice~$\Cl(\negtinysp P)\dtinysp$:
for every subset $\Gamma$ of $\Nuc(\negtinysp P)$ the join $\tJoin\negtinysp\Gamma$,
	taken in $\Cl(\negtinysp P)$,
is~a~nucleus, thus it is the~join of~$\Gamma$ in~$\Nuc(\negtinysp P)\tinysp$;
also, $\fix(\tinysp\tJoin\negtinysp\Gamma)=\tInters_{\tinysp\gamma\in\Gamma}\fix(\gamma)$
by the~antiisomorphism between $\Cl(\negtinysp P)$ and~$\ClSys(\negtinysp P)$.

Let $\Gamma$ be a~subset of $\Nuc(\negtinysp P)$,
	and set $\coll{C}\defeq\set{\rightdt{\fix(\gamma)}\suchthat\leftdt{\gamma\narrt\in\Gamma}}$.
The join of $\coll{C}$ in $\NucSys(\negtinysp P)$ is $B\defeq\nucsys(\tUnion\coll{C})$,
where $\nucsys=\nucsys_P$ is the closure operator on $\PowP$
	determined by the closure rules $\Rulesdflt(\negtinysp P)\union\Rulesnuc(\negtinysp P)\tinysp$.
If $\beta$ is the meet of $\Gamma$ in $\Nuc(\negtinysp P)$, then $\dtinysp\fix(\beta)\narrt=B$.

The nonempty finite meets in $\Cl(\negtinysp P)$ as well as in $\Nuc(\negtinysp P)$ are calculated pointwise.
If~$P$~has a~top element~$\top$,
then the constant map $P\to P : x\mapsto \top$
is the top element of both $\Cl(\negtinysp P)$ and $\Nuc(\negtinysp P)$.
If $P$ does not have a~top element,
then the top element of $\Nuc(\negtinysp P)$ may be different (thus strictly smaller)
	than the top element of $\Cl(\negtinysp P)$.

\thmskip

\begin{prop}\label{prop:nuc-enab-meet-semilatt--nuc-generd-by-prenucs}
Let\/ $P$ be a~nuclear-enabled meet-semilattice, and let\/~$\Gamma$ be a~set of prenuclei on~$P$.
The closure operator\/~$\chi$ on\/~$P$ which has\/~$\dtinysp\fix(\chi)=\fix(\Gamma)\tinysp$,
	the least of the closure operators on\/ $P$ that are above\/~$\Gamma$,
is a~nucleus.
\end{prop}

\interskip

\begin{myproof}
For every $\gamma\in\Gamma$ the set $\dtinysp\fix(\gamma)$ obeys $\Rulesnuc(\negtinysp P)$
	by Lemma~\ref{lem:Fix(prenuc)-obeys-nuclear-rules},
therefore $\dtinysp\fix(\Gamma)=\tInters_{\tinysp\gamma\in\Gamma}\fix(\gamma)$ obeys $\Rulesnuc(\negtinysp P)$,
and besides that, $\dtinysp\fix(\Gamma)$ is a~closure system
	by Theorem~\ref{thm:in-dfltenab-cl-generd-by-precls-&-induct}.
By Lemma~\ref{lem:enough-nucrules-clop-obeys-nucrules=>clop-is-nuc}
	the closure operator $\chi$ preserves binary meets, that is, its is a~nucleus.
By Theorem~\ref{thm:in-dfltenab-cl-generd-by-precls-&-induct}
	the closure operator $\chi$ is the least closure operator on $P$ that is above~$\Gamma$.
\end{myproof}

\thmskip

The following theorem is a~do-it-by-maximal-elements analogue
	of Theorem~\ref{thm:frame-of-nuclei-on-preframe}.

\thmskip

\begin{thm}\label{thm:when-Nuc(meet-semilatt)-is-a-frame}
If\/ $P$ is a~default-enabled meet-semilattice,
and for all\/ $a,\tinysp b\in P$
	the set\/~$(a\astbimpl b)$ is default-enabled within\/~$P$
		{\rm(}so it certainly has a~ceiling\/{\rm)},
then the complete lattice\/~$\Nuc(\negtinysp P)$ is a~frame.
\end{thm}

\interskip

\begin{myproof}
We shall prove,
    for all $\beta\in\Nuc(\negtinysp P)$ and all $\Gamma\subseteq\Nuc(\negtinysp P)\tinysp$,
the following identity:
\begin{equation*}
\beta\meet\tJoin\Gamma \Eq \tJoin_{\!\gamma\in\Gamma}(\beta\meet\gamma)~.
\end{equation*}
The inequality $\geq$ is clear, so it remains to prove the converse inequality~$\leq\dtinysp$.
We~write $\delta\defeq\tJoin\Gamma$
	and $\deltapr\defeq\tJoin_{\!\gamma\in\Gamma}(\beta\meet\gamma)\tinysp$.
We have to prove that
\begin{equation*}
(\beta\narrt\meet\delta)(x) \Eq \beta(x)\meet\delta(x) \Leq \deltapr(x)
	\qquad\quad \text{for every $x\in P$}
\end{equation*}
(recall that finite meets in $\Nuc(\negtinysp P)$ are calculated pointwise).
Let $A \defeq \bigl(\beta(x)\astbimpl\deltapr(x)\bigr)
	= \bigset{z\narrdt\in P\bigsuchthat \beta(x)\meet z\leq\deltapr(x)}\tinysp$.
Clearly $x\in A$, and $A$ is by assumption default-enabled within~$P$.
We~prove that the set $A$ is closed under $\Gamma$
precisely as we did in the proof of Theorem~\ref{thm:frame-of-nuclei-on-preframe}.
By the induction principle,
	formulated in Theorem~\ref{thm:in-dfltenab-cl-generd-by-precls-&-induct},
it follows that $A$ is closed under $\delta$, and hence that $\delta(x)\in A$,
which means that $\beta(x)\meet\delta(x)\leq\deltapr(x)\tinysp$.
\end{myproof}

\thmskip

The last two propositions above specialize to propositions about preframes
since in a~preframe every set of the form $(a\astbimpl b)$ is a~sub-dcpo
	(in~fact it is a~\Scottcl\ subset)
and as such it is default-enabled within the preframe.
The act of specialization requires the use of the~axiom of~choice,
so we are lucky, again, that we have already proved the specialized propositions.

\txtskip

Notice that we somehow managed to prove the lemmas, the propositions, and the theorems of this section
without ever using the~axiom of~choice, or even the~law of~excluded middle,
	whatever this might be good for.

\section{Two convex geometries associated with a dcpo}
\label{sec:conv-geoms-assoc-with-dcpo}

Let $E$ be a~set and $\gamma$ a~closure operator on~$\PowE$.
The closure operator $\gamma$ is said to be \notion{convex}
	if it satisfies the following \notion{anti-exchange axiom}:
\begin{itemize}[topsep=1ex,itemsep=1ex,leftmargin=5.3em]
\item[(AE)\:] For every subset $A$ of $E$ and all elements $x$, $y$ of $E\tinysp$,\\
	if $x,\tinysp y\notin \gamma(A)$ and $x\neq y$
		and $x\in\gamma(A\narrt\union\set{y})\tinysp$,
	then $y\notin\gamma(A\narrt\union\set{x})\tinysp$.
\end{itemize}%

The anti-exchange axiom is equivalent to the following condition:
\begin{itemize}[topsep=1ex,itemsep=1ex,leftmargin=5.3em]
\item[(CAS)\:] For every $\gamma$-closed subset $C$ of $E$ and all elements $x$, $y$ of $E\tinysp$,\\
	if $x,\tinysp y\notin C$
		and $\gamma(C\narrt\union\set{y})=\gamma(C\narrt\union\set{x})$,
	then~$x=y\tinysp$.
\end{itemize}%
Let $\gamma$ be an arbitrary~closure operator on $\PowE\tinysp$.
For every subset~$A$ of~$E$,
the closure operator $\gamma$ induces the preorder  $\leq_A$ on the set $E\setdiff A$,
	where $x\narrt{\leq_A}y$ iff $x\in\gamma(A\narrt\union\set{y})$.
The~condition~(CAS) requires that for every $\gamma$-closed subset~$C$ of~$E$
	the preorder $\leq_C$ on $E\setdiff C$ is~anti\-symmetric,
	that is, that it is a~partial order.

A~\notion{convex geometry} is a~structure $\pair{E,\gamma}$
where $E$ is a~set and $\gamma$ is a~convex closure operator on~$\PowE$.

\txtskip

It is clear from the form of (AE), or of (CAS),
that handling of convex geometries will require invocations of the~law of~excluded middle.
However, we will get by without any help from the~axiom of~choice.

\txtskip

The following proposition is the main result of this section.
It generalizes Proposition~5\nobreakdash-5.1 in~\cite{Gratzer-Wehrung-LT:STA2-2016}.
It will be proved in due time.

\thmskip

\begin{prop}\label{prop:dcpo-P==>(P,clsysP)-(P,dcclsysP)-convex-geoms}
For a~dcpo\/ $P$,
$\pair{P,\dtinysp\clsys_P}$ and\/ $\pair{P,\dtinysp\dcclsys_P}$ are convex geometries.
\qed
\end{prop}

\thmskip

We start with some very general observations.

\thmskip

\begin{lem}\label{lem:clsys-C-lowset-A==>clsys-CunionA}
Let\/ $P$ be a~poset.
If\/~$C$ is a~closure system in\/~$P$ and\/~$A$ is a~lower set of\/~$P$,
then\/~$C\union A$ is a~closure system in\/~$P$.
\end{lem}

\interskip

\begin{myproof}
Let $x\in P$.
We have to prove that $\allabove{(C\narrdt\union A)}{x}$ has a~least element~$u\tinysp$.
If $x\in  A$ then $u=x\tinysp$.
If $x\notin A$, then $\allabove{(C\narrdt\union A)}{x}=\allabove{C}{x}$ since $A$ is a~lower set,
	and $u=\cl_C(x)$.
\end{myproof}

\thmskip

Notice the application of the~law of~excluded middle in the proof,
	where we consider the pair of cases $x\in A$ and $x\notin A\tinysp$ as being exhaustive.
More applications of the law of excluded middle lay ahead,
	but we will no longer raise alarums over them.

\thmskip

\begin{lem}\label{lem:clsys-dirclsd-C-fingen-lowset-A=>clsys-dirclsd-CunionA}
Let\/ $P$ be a~poset.
If\/ $C$~is a~closure system in\/~$P$
	that is closed under existing directed joins in\/~$P$,
and\/ $A$~is a~finitely generated lower set of\/~$P$,
then\/ $C\union A$~is a~closure system in\/~$P$
	that is closed under existing directed joins in\/~$P$.
\qed
\end{lem}

\thmskip

The set $C\union A$ is a~closure system
	by Lemma~\ref{lem:clsys-C-lowset-A==>clsys-CunionA}.
It remains to prove that $C\union A$ is closed under existing directed joins.
Since $A$ is a~union of finitely many principal ideals,
and every principal ideal is closed under all existing joins
	hence under all existing directed joins,
the desired result is a~consequence of the following lemma:

\thmskip

\begin{lem}\label{lem:union-of-two-diclsd-sets-is-dirclsd}
If subsets\/ $A$ and\/ $B$ of a~poset\/~$P$
	are closed under existing directed joins in\/~$P$,
then the subset\/ $A\union B$ is closed under existing directed joins in\/~$P$.
\end{lem}

\interskip

\begin{myproof}
Let a~directed subset $Y$ of $A\union B$ have a~join~$u$ in~$P$.
We consider two cases.

Case~1: $Y\inters A$ is a~cofinal subset of~$Y$.
The set $Y\inters A$ has the same upper bounds in~$P$ as the set $Y$,
thus the join $u$ of $Y$ in $P$ is also the join of $Y\inters A$ in $P$.
Since $Y\inters A$ is a~directed subset of $A$ and $A$ is \dircl,
we have $u\in A\tinysp$.

Case~2: $Y\inters A$ is not a~cofinal subset of~$Y$.
There exists $b\in Y$ such that the set $\allabove{Y\negdtinysp}{b}$ is disjoint with $A$
and is therefore included in~$B$.
Since $\allabove{Y\negdtinysp}{b}$ is a~cofinal subset of~$Y$, so~is~$Y\inters B$,
and by Case~1, with $A$ and $B$ exchanged, it follows that $u\in B$.
\end{myproof}

\thmskip

Let $P$ be a~poset.
We denote by $\DcClSys(\negtinysp P)$ the subposet of $\PowP$
consisting of all closure systems~in~$P$ that are closed under existing directed joins in~$P$.

Suppose that a~poset~$P$ has the property that $\ClSys(\negtinysp P)$ is a~closure system in $\PowP$,
that~is, that it is closed under all intersections.
Then for any subset~$X$ of~$P$
	there is the least set in $\ClSys(\negtinysp P)$ that includes~$X$,
which we denote by~$\clsys_P(X)$.
But then, by Lemma~\ref{lem:P-poset&C-clsys-in-Pow(P)=>Dc(C)-clsys-in-Pow(P)},
the set $\DcClSys(\negtinysp P)$ is likewise a~closure system in~$\PowP$,
and so for any subset $X\subseteq P$
	there is the least set in $\DcClSys(\negtinysp P)$ that includes~$X$,
which we denote by $\dcclsys_P(X)$.

\thmskip

\begin{prop}\label{prop:ClSys(P)-clsys-in-PowP=>clsysP-satisfies-(CAS)}
Let\/ $P$ be a~poset.
If\/ $\ClSys(\negtinysp P)$ is a~closure system in\/~$\PowP$,
then\/ $\pair{P,\dtinysp\clsys_P}$ is a~convex geometry.
\end{prop}

\interskip

\begin{myproof}
Let $C$ be a~closure system in $P$,
and suppose that $x,\tinysp y\in P$ are not in $C$ and that
    $\clsys_P(C\narrt\union\set{y})
	= \clsys_P(C\narrt\union\set{x})\tinysp$.
The set $C\union\ldown y$ is,
	according to Lemma~\ref{lem:clsys-C-lowset-A==>clsys-CunionA},
		a~closure system,
	and it includes $C\union\set{y}\tinysp$,
thus it includes $\widedt{\clsys_P(C\narrt\union\set{y})}$.
Now~from
    $x \in \clsys_P(C\narrt\union\set{x})
	= \clsys_P(C\narrt\union\set{y})
	\subseteq C\union\ldown y$
and $x\notin C$ it follows that $x\in\ldown y$,
	that~is, that~$x\leq y\tinysp$.
Likewise we see that $y\leq x$, and we conclude that $x=y\tinysp$.
The closure operator $\widedt{\clsys_P}$ on $\PowP$ satisfies the condition~(CAS).
\end{myproof}

\interthmskip

\begin{prop}\label{prop:ClSys(P)-clsys-in-PowP=>dcclsysP-satisfies-(CAS)}
For a~poset~$P$\!,
if\/ $\ClSys(\negtinysp P)$ is a~closure system in\/~$\PowP\!$,
then\/ $\pair{P,\dtinysp\dcclsys_P}$ is a~convex geometry.
\qed
\end{prop}

\thmskip

The proof is the same
	as that of Proposition~\ref{prop:ClSys(P)-clsys-in-PowP=>clsysP-satisfies-(CAS)}
except that it uses Lemma~\ref{lem:clsys-dirclsd-C-fingen-lowset-A=>clsys-dirclsd-CunionA}
	instead of Lemma~\ref{lem:clsys-C-lowset-A==>clsys-CunionA}.

\thmskip

\begin{myproof}[\bf Proof of Proposition~\ref{prop:dcpo-P==>(P,clsysP)-(P,dcclsysP)-convex-geoms}.]
As $P$ is a~dcpo, $\ClSys(\negtinysp P)$ is a~closure system in $\PowP$.
Now~apply Proposition~\ref{prop:ClSys(P)-clsys-in-PowP=>clsysP-satisfies-(CAS)}
	and Proposition~\ref{prop:ClSys(P)-clsys-in-PowP=>dcclsysP-satisfies-(CAS)}.
\end{myproof}

\thmskip

The following proposition has essentially the same proof
	as Proposition~\ref{prop:dcpo-P==>(P,clsysP)-(P,dcclsysP)-convex-geoms}.

\thmskip

\begin{prop}\label{prop:P-dfltenab-==>(P,clsysP)-(P,dcclsysP)-convex-geoms}
If\/ $P$ is a~default-enabled poset,
then\/ $\pair{P,\dtinysp\clsys_P}$ and\/ $\pair{P,\dtinysp\dcclsys_P}$ are convex geometries.
\qed
\end{prop}

\thmskip

We are not done yet.
For a~dcpo $P\negtinysp$, the closure operator $\widedt{\clsys_P}$ is convex for a~reason,
the reason being that this closure operator is \emph{acyclic}.
Below we give the definition of acyclic closure operators,
but only after the definition of a~funnel for a~closure operator.

\txtskip

Let $E$ be a~set and $\gamma$ a~closure operator on $\PowE\tinysp$.

A~\notion{funnel} for the closure operator~$\gamma$
is a~preorder~$\leq$ on~$E$ that has the following property:
to every $X\subseteq E$ and every $y\in\gamma(X)$
there is a~subset~$Z$ of~$X$ such that $y\leq Z$ and $y\in\gamma(Z)\tinysp$.
Mark that a~preorder~$\leq$ on~$E$ is a~funnel for~$\gamma$
iff for every $X\subseteq E$ and every $y\in\gamma(X)$
	it follows that $y\in\gamma(\allabove{X\negdtinysp}{y})\tinysp$.
If~$\leq$ is a~funnel for $\gamma$,
then we also say that $\gamma$ \notion{has a~funnel}~$\leq\,$.

We shall say that the closure operator $\gamma$ is \notion{acyclic}
if it has an~antisymmetric funnel, that is, a~funnel which is a~partial order on~$E\tinysp$.

\thmskip

\begin{prop}\label{prop:properts-equiv-to-<=-beig-funnel-for-clop}
Let\/ $E$ be a~set,
	$\gamma$ a~closure operator on\/~$\PowE\tinysp$,
	and\/ $\leq$ a~preorder on\/~$E\tinysp$.
The following are equivalent:

\begin{itemize}[topsep=.5ex,itemsep=.5ex,leftmargin=3.7em]
\item[{\rm(1)}\:] $\leq$~is a~funnel for\/~$\gamma\dtinysp$;
\item[{\rm(2)}\:] for all\/ $X,\tinysp U\negdtinysp\subseteq E\tinysp$,
	if\/ $U$~is an upper set of~$\pair{E,{\leq}\tinysp}\tinysp$,
	then\/ $\gamma(X)\inters U\subseteq\gamma(X\narrt\inters U)\dtinysp$;
\item[{\rm(3)}\:] for all $X\subseteq E$ and all $y\in E$ we have
	$\allabove{\gamma(X)}{y} \subseteq \gamma(\allabove{X}{y})\tinysp$.
\end{itemize}%
\end{prop}

\negdisplayhalfskip
\interskip

\begin{myproof}
(1)$\Implies$(2).\,
Assume~(1), and let $X,\tinysp U\subseteq E$ with~$U$ an upper set of~$\pair{E,{\leq}\tinysp}\tinysp$.
Let~$u$~be an element of~$\gamma(X)\inters U$ and write $Z\defeq\allabove{X\negdtinysp}{u}\tinysp$.
Then $Z\subseteq U$ because $U$ is an upper set,
	and~$u\in\gamma(Z)$ since $\leq$ is a~funnel for~$\gamma\tinysp$,
thus $Z\subseteq X\inters U$ and $u\in\gamma(Z)\subseteq\gamma(X\narrt\inters U)\tinysp$.

(2)$\Implies$(3) holds by specialization ($\widet{U\leftt=\lup y}$).

(3)$\Implies$(1).\,
Assuming (3), suppose that $y\in\gamma(X)$;
then $y\in \allabove{\gamma(X)}{y} \widet\subseteq \gamma(\allabove{X\negdtinysp}{y})$.
\end{myproof}

\thmskip
\pagebreak[3]

And why are the acyclic closure operators of such interest to us?
This is why:

\thmskip

\begin{prop}\label{prop:generd-Lemma8-3.23-in-LT-STA-2}
Let\/ $E$ be a~set, $\leq$ a~preorder on\/~$E\tinysp$,
	and\/ $\gamma$ a~closure operator on\/~$\PowE\tinysp$.
If\/~$\leq$ is a~funnel for\/~$\gamma\tinysp$, then the following statements are true:
\begin{itemize}[topsep=.5ex,itemsep=.5ex,leftmargin=3.7em]
\item[{\rm(i)}\:] For all\/ $A\subseteq E$ and all\/ $x,\tinysp y\in E$,
	if\/ $x\notin\gamma(A)$ and\/ $x\in \gamma(A\narrt\union\set{y})$,
	then\/ $x\leq y\tinysp$.
\item[{\rm(ii)}\:] If\/ $\leq$ is a~partial order, then the closure operator\/ $\gamma$ is convex.
\end{itemize}%
\end{prop}

\negdisplayhalfskip
\interskip

\begin{myproof}
(i)\, Assume that $A$, $x$, $y$ satisfy the premises.
Since $x\in\gamma(A\narrt\union\set{y})$ and~$\leq$~is a~fun\-nel for~$\gamma\tinysp$,
it follows that $x\in\gamma\bigl(\allabove{(A\narrt\union\set{y})}{x}\bigr)\tinysp$.
Now the set $\allabove{(A\narrt\union\set{y})}{x}$ must contain~$y$ 
since~otherwise we would have $\allabove{(A\narrt\union\set{y})}{x} = \allabove{A}{x}$
	and $x\in\gamma(\allabove{A}{x})\subseteq\gamma(A)\tinysp$,
contrary to assumptions.
That~is,~we~have~$x\leq y\tinysp$.

(ii)\, For every $\gamma$-closed $C\subseteq P$
	and for all $x,\tinysp y\in P\setdiff C$,
if $\gamma(C\narrt\union\set{x})=\gamma(C\narrt\union\set{y})\tinysp$,
then by part (i) it follows that $x\leq y$ and $y\leq x$, whence $x=y\tinysp$.
The closure operator $\gamma$ satisfies the condition~(CAS).
\end{myproof}

\thmskip

Proposition~\ref{prop:generd-Lemma8-3.23-in-LT-STA-2} generalizes
	Lemma~8-3.23 in~\cite{Gratzer-Wehrung-LT:STA2-2016},
from algebraic closure operators of poset type to arbitrary acyclic closure operators.%
\footnote{\,Mark that the partial order in Proposition~\ref{prop:generd-Lemma8-3.23-in-LT-STA-2}
is the converse of the partial order in Lemma~8-3.23.}
The proof of the proposition is \emph{not} completely modeled after the proof of Lemma~8-3.23,
since the latter proof uses Lemma~8-3.2
	which provides a~useful consequence of algebraicity of the closure operator,
and the proof above has no use 
	for such a~lemma.

\thmskip

\begin{prop}\label{prop:ClSys(P)-clsys-in-PowP=>clsysP-is-acyclic}
Let\/ $P$ be a~poset.
If\/ $\ClSys(\negtinysp P)$ is a~closure system in\/~$\PowP$,
then the partial order of\/~$P$ is a~funnel for the closure operator\/~$\clsys_P$,
which is therefore~acyclic.
\end{prop}

\interskip

\begin{myproof}
Let $X\subseteq P$ and $y\in\clsys_P(X)\dtinysp$;
we have to prove that $y\in\clsys_P(\allabove{X\negdtinysp}{y})\tinysp$.

The set $\dtinysp\clsys_P(\allabove{X\negdtinysp}{y})
				\union(P\setdiff\narrt\lup y)$
    is by Lemma~\ref{lem:clsys-C-lowset-A==>clsys-CunionA}
	a~closure system in~$P\tinysp$;
it~includes the set~$X$,
so it includes the closure system $\widedt{\clsys_P(X)}$
	and hence contains the element~$y\dtinysp$;
since $y\notin P\setdiff\lup y$,
we conclude that $y\in\clsys_P(\allabove{X\negdtinysp}{y})\tinysp$.
\end{myproof}

\interthmskip

\begin{cor}\label{cor:P-dcpo-or-dfltenab==>clsysP-acyclic}
If\/ $P$ is a~dcpo or a~default-enabled poset,
then\/ $\widedt{\clsys_P}$ is acyclic.%
\footnote{\,The only property of a~dcpo or a~default-enabled poset which we need here
is that every intersection of its closure systems is a~closure system.
We proved this property separately for dcpos and for default-enabled posets,
even though every dcpo is default-enabled;
we did this because the specialization from default-enabled posets to dcpos
	by necessity involves the axiom of choice.
This~is also the reason why the corollary mentions both types of posets.}
\qed
\end{cor}

\thmskip

Therefore, if $P$ is a~dcpo or a~default-enabled poset,
then the closure operator $\widedt{\clsys_P}$ on $\PowP$
is convex \emph{because} it is acyclic,
	in view of Proposition~\ref{prop:generd-Lemma8-3.23-in-LT-STA-2}.
This proves again the first halves of
	Proposition~\ref{prop:dcpo-P==>(P,clsysP)-(P,dcclsysP)-convex-geoms}
	and Proposition~\ref{prop:P-dfltenab-==>(P,clsysP)-(P,dcclsysP)-convex-geoms};
but the original direct proofs of those halves were markedly simpler,
so one can be excused for not seeing the point
	of the new proofs that go the roundabout way through acyclicity.
However, Corollary~\ref{cor:P-dcpo-or-dfltenab==>clsysP-acyclic} is of independent interest.
For example, the dcpo part of the corollary,
	partnered with Proposition~\ref{prop:properts-equiv-to-<=-beig-funnel-for-clop},
is a~special case of Lemma~4.3~in~\cite{Ranzato-CCPOFCL-1999},
and this special case is then used in the proof of~Theorem~5.2 in~\cite{Ranzato-CCPOFCL-1999}
which is about a~poset that satisfies the ascending chain condition and is therefore (trivially) a~dcpo.


\bibliographystyle{alpha}

\bibliography{clops-on-dcpos-references}


\end{document}